\documentclass{amsart}

\usepackage{amsmath,amsopn,amssymb,amsthm,mathtools,mathrsfs,yhmath,tensor,graphicx,geometry,stmaryrd,tabulary,booktabs,leftindex,dsfont}
\usepackage[euler]{textgreek}

\usepackage{xcolor}

	\definecolor{UCIB}{HTML}{0064A4}
	\definecolor{UCSDB}{HTML}{00629B}
	\definecolor{OSUR}{HTML}{ba0c2f} % R 186   G 12   B 47  

\usepackage{accents}

\usepackage{upgreek,textgreek}
	
\usepackage{hyperref}
\hypersetup{backref, urlcolor = UCIB, citecolor=blue, colorlinks=true}

\newcommand{\bhyperlink}[3][black]{\hyperlink{#2}{\color{#1}{#3}}}%

\usepackage{ifthen}

\newcommand{\mybinom}[3][0.8]{\scalebox{#1}{$\dbinom{#2}{#3}$}}
	
\usepackage{lipsum} % just for the example

\usepackage{tikz}
\usetikzlibrary{calc}
\usetikzlibrary{intersections}
\usetikzlibrary{graphs,graphs.standard,quotes}

\usepackage{enumitem}

\usepackage{tikz}

\usepackage{pgfplots}
\pgfplotsset{compat=1.17}
\usetikzlibrary{quotes,angles}

\usepackage{booktabs}

\usepackage{nicematrix}

\usepackage{xcolor}
\DeclareMathOperator\osc{osc}

\DeclareMathOperator\Id{Id}

\DeclareMathOperator\CY{CY}

\theoremstyle{plain}% default
\newtheorem{fthm}{Theorem}[section]
\newtheorem*{fthm*}{Theorem}
\newtheorem{flemma}{Lemma}[section]
\newtheorem*{flemma*}{Lemma}
\newtheorem{fprop}{Proposition}[section]
\newtheorem*{fprop*}{Proposition}
\newtheorem{fcor}{Corollary}[section]
\newtheorem*{fcor*}{Corollary}

\theoremstyle{definition}
\newtheorem{fdefi}{Definition}[section]
\newtheorem*{fdefi*}{Definition}

\newtheorem*{fexmp*}{Example}

\theoremstyle{remark}
\newtheorem{frmk}{Remark}[section]
\newtheorem*{frmk*}{Remark}

\newtheorem*{fconj*}{Conjecture}

\newtheorem*{fclaim*}{Claim}

\newtheorem*{fquest*}{Question}

\newcommand{\Address}{{% additional braces for segregating \footnotesize
  \bigskip
  \footnotesize

  Chao-Ming Lin, \textsc{Department of Mathematics, Ohio State University, OH}\par\nopagebreak
  \textit{E-mail address:} \href{lin.4579@osu.edu}{lin.4579@osu.edu}\par\nopagebreak
  \textit{Personal Website:} \href{https://chaominl.github.io}{https://chaominl.github.io}

}}
\newcommand{\Acknow}{{% additional braces for segregating \footnotesize
  \bigskip
\textbf{Acknowledgements:} The author is grateful to Zhiqin~Lu, Xiangwen~Zhang, and Bo Guan for giving me enlightening help.
The author would like to thank Richard Schoen, Tristan Collins, and Duong Hong Phong for their interest in this work and support. 
The author would also like to thank 
Chung-Jun Tsai, Mu-Tao Wang, Adam Jacob, Mao-Pei Tsui, Yng-Ing Lee, Man Chun Lee, Yu-Shen Lin, Martin Man Chun Li, and Ziming Ma
for the helpful discussions when the author was visiting the National Taiwan University, the National Center for Theoretical Sciences, and the Chinese University of Hong Kong during the Summer of 2023.
}}

\begin{document}

\title{On the Solvability of General Inverse $\sigma_k$ Equations}
\author{Chao-Ming Lin}
%\address{Department of Mathematics\\ Ohio State University, Columbus, OH}
%\date{}

\begin{abstract}
We prove that if there exists a $C$-subsolution to a constant coefficients strictly $\Upsilon$-stable general inverse $\sigma_k$ equation, then there exists a unique solution. As a consequence, this result covers all the analytical results of the classical strictly $\Upsilon$-stable general inverse $\sigma_k$ equations, for example, the complex Monge--Ampère equation, the complex Hessian equation, the J-equation, the deformed Hermitian--Yang--Mills equation, etc. Hence, we confirm an analytical conjecture by Collins--Jacob--Yau \cite{collins20201} of the solvability of the deformed Hermitian--Yang--Mills equation. Their conjecture states that the existence of a $C$-subsolution to a supercritical phase deformed Hermitian--Yang--Mills equation gives the solvability.
\end{abstract}
\maketitle
\vspace{-0.6cm}

\section{Introduction}

Let $(M, \omega)$ be a compact connected Kähler manifold of complex dimension $n$ with a Kähler form $\omega$ and $[\chi_0] \in H^{1,1}(M; \mathbb{R})$, where $H^{1,1}(M; \mathbb{R})$ is the $(1, 1)$-Dolbeault cohomology group. The study of the solvability of the following equation is widely considered:
\begin{align*}
\label{eq:1.1}
\chi^n = c_{n-1} \mybinom[0.8]{n}{n-1} \chi^{n-1}  \wedge \omega^{1} + \cdots + c_{1} \mybinom[0.8]{n}{1} \chi^{1}  \wedge \omega^{n-1} + c_{0} \mybinom[0.8]{n}{0}   \omega^{n}    =  \sum_{k = 0}^{n-1} c_k \mybinom[0.8]{n}{k} \chi^{k}  \wedge \omega^{n-k}, \tag{1.1}
\end{align*}where $c_k$ are real functions on $M$ and $\chi \in [\chi_0]$ is a real smooth, closed $(1, 1)$-form. We call an equation having the same format as equation (\ref{eq:1.1}) a degree $n$ \emph{general inverse $\sigma_k$ equation}. \smallskip 

A general inverse $\sigma_k$ equation (\ref{eq:1.1}) is very likely to be ill-posed. For example, by letting $[\chi_0]$ be a Kähler class and $c_k \leq 0$ for all $k \in \{0, \cdots, n-1\}$, then it is straightforward to see that equation (\ref{eq:1.1}) is not solvable in this case. 
However, some special combinations of the coefficients raise some famous equations. For instance, by letting $[\chi_0]$ be a Kähler class, $c_k = 0$ for all $k \in \{1, \cdots, n-1\}$, and $c_0$ be a positive function, equation (\ref{eq:1.1}) becomes the well-known complex Monge--Ampère equation in the Calabi conjecture \cite{calabi1954kahler, calabi1957kahler}, which was solved by Yau \cite{yau1978ricci}. For more examples, the J-equation, the complex Hessian equation, the deformed Hermitian--Yang--Mills equation, etc. \smallskip

The solvability of these equations has deep connections with the geometric properties of the underlying manifold. For example, 
inspired by the study of the Hermitian--Yang--Mills connections by Donaldson \cite{donaldson1985anti} and Uhlenbeck--Yau \cite{uhlenbeck1986existence}, Donaldson \cite{donaldson1999moment} discovered the J-equation when studying the Mabuchi energy \cite{mabuchi1986k} using the moment maps. For the analytical part, the J-equation was studied extensively by Chen \cite{chen2000lower}, Collins--Székelyhidi \cite{collins2017convergence}, and Song--Weinkove \cite{song2008convergence}. For the numerical part, Lejmi--Székelyhidi \cite{lejmi2015j} conjectured that the existence of the solution to the J-equation is equivalent to a certain numerical stability condition. Chen \cite{chen2021j} and Song \cite{song2020nakai} studied the solvability and the stability of the J-equation using a Nakai--Moishezon type criterion inspired by the work of Demailly--Păun \cite{demailly2004numerical} and proved this conjecture by Lejmi--Székelyhidi \cite{lejmi2015j}. \smallskip

Another crucial example will be the {deformed Hermitian--Yang--Mills equation}, which will be abbreviated to the dHYM equation from now on.
Motivated by mirror symmetry in string theory, the dHYM equation was discovered around the same time by Mariño--Minasian--Moore--Strominger \cite{marino2000nonlinear} and Leung--Yau--Zaslow \cite{leung2000special} using different viewpoints. The dHYM equation was initiated by Jacob--Yau \cite{jacob2017special} and can be formulated as follows:
\begin{align}
\label{eq:1.2}
\Im \bigl ( \omega + \sqrt{-1} \chi \bigr )^n &= \tan \bigl (   {\theta}   \bigr ) \cdot \Re \bigl ( \omega + \sqrt{-1} \chi \bigr )^n. \tag{1.2}
\end{align}Here $\Im$ and $\Re$ are the imaginary and real parts, respectively, and $ {\theta}$ is a topological constant determined by the cohomology classes $[\omega]$ and $[\chi_0]$.  For the analytical part, Collins--Jacob--Yau \cite{collins20201} showed that if there exists a slightly stronger $C$-subsolution, then the supercritical phase dHYM equation is solvable. Collins--Jacob--Yau conjectured that this slightly stronger condition can be removed, that is, the existence of a $C$-subsolution will lead to the solvability. This analytical conjecture was confirmed by Pingali \cite{pingali2019deformed} when the complex dimension equals three. By introducing new real algebraic geometry techniques, the author \cite{lin2023d} confirmed this analytical conjecture when the complex dimension equals four. For the numerical part, Collins--Jacob--Yau conjectured that the existence of the solution to the dHYM equation (\ref{eq:1.2}) is equivalent to a certain stability condition for all analytic subvarieties and confirmed this numerical conjecture for complex surfaces. 
Chen \cite{chen2021j} proved a Nakai--Moishezon type criterion for the dHYM equation when the phase is supercritical under a slightly stronger condition that these holomorphic intersection numbers have a uniform lower bound independent of analytic subvarieties. Chu--Lee--Takahashi \cite{chu2021nakai} improved the result by Chen \cite{chen2021j} without assuming a uniform lower bound for these holomorphic intersection numbers on projective manifolds. \smallskip

There are also other significant works in general inverse $\sigma_k$ equations recently. In \cite{lin2023c}, the author gave more detailed descriptions of many significant works. The interested reader is referred to \cite{chu2022hypercritical, collins2020stability, collins2021moment, collins2018deformed, datar2021numerical, fang2013convergence, fang2011class, fang2023,  jacob2019weak, jacob2020deformed, lin2020, phong2017fu, phong2018anomaly, phong2019estimates, phong2021fu, schlitzer2021deformed, yuan2006global} and the references therein. \bigskip

At a pointwise level, if we write equation (\ref{eq:1.1}) in terms of the eigenvalues of the Hermitian endomorphism $\Lambda = \omega^{-1} \chi$ at a point, then we can rewrite equation (\ref{eq:1.1}) as
\begin{align*}
\label{eq:1.3}
\lambda_1 \cdots \lambda_n = \sum_{k = 0}^{n-1} c_k \sigma_k(\lambda_1, \cdots, \lambda_n) = \sum_{k = 0}^{n-1} c_k \sigma_k(\lambda). \tag{1.3}
\end{align*}Here, $\lambda_i$ are the eigenvalues of $\Lambda$, $\sigma_k(\lambda_1, \cdots, \lambda_n)$ is the $k$-th elementary symmetric polynomial of $\{\lambda_1, \cdots, \lambda_n\}$, and we denote $\sigma_k(\lambda_1, \cdots, \lambda_n)$ by $\sigma_k(\lambda)$ for convenience. We call a multilinear polynomial having the format $\lambda_1 \cdots \lambda_n - \sum_{k = 0}^{n-1} c_k \sigma_k(\lambda)$ a degree $n$ \emph{general inverse $\sigma_k$ multilinear polynomial}.
Assume that $[\chi_0]$ is a Kähler class, if one wishes to find a Kähler representative in $[\chi_0]$ solving equation (\ref{eq:1.1}), then pointwise the level set of $\lambda_1 \cdots \lambda_n - \sum_{k = 0}^{n-1} c_k \sigma_k(\lambda) = 0$ should be contained in the positive orthant. Naturally, we have the following definition. 
\hypertarget{D:1.1}{\begin{fdefi}[\texorpdfstring{$\Upsilon$}{}-stableness]}
Let $f(\lambda) \coloneqq   \lambda_1 \cdots \lambda_n - \sum_{k = 0}^{n-1} c_k \sigma_k(\lambda)$ be a general inverse $\sigma_k$ multilinear polynomial and $\Gamma_f^{n}$ be a connected component of $\{f(\lambda) > 0\}$. We say that this connected component $\Gamma^n_f$ of $f(\lambda)$ is $\Upsilon$-stable if
\begin{align*}
\Gamma^n_f \subseteq  q +  \Gamma_n \text{ for some } q \in \mathbb{R}^n, 
\end{align*}where $\Gamma_n$ is the positive orthant of $\mathbb{R}^n$. We say that this connected component $\Gamma^n_f$ is strictly $\Upsilon$-stable if it is $\Upsilon$-stable and the boundary $\partial \Gamma^n_f$ is contained in the $\Upsilon_1$-cone. 
\end{fdefi} 
The $\Upsilon_k$-cones will be defined later in Section~\ref{sec:2.1} for $k \in \{1, \cdots, n-1\}$. In particular, the $\Upsilon_1$-cone is the $C$-subsolution cone introduced by Székelyhidi \cite{szekelyhidi2018fully} and Guan \cite{guan2014second}. \smallskip

Pointwise, if a level set of equation (\ref{eq:1.1}) is strictly $\Upsilon$-stable for every point on $M$, then we call equation (\ref{eq:1.1}) a \emph{strictly $\Upsilon$-stable general inverse $\sigma_k$ equation}. 

\hypertarget{R:1.1}{\begin{frmk}}
The following general inverse $\sigma_k$ equations are all strictly $\Upsilon$-stable:
\begin{itemize}[leftmargin = 2cm]
\item Complex Monge--Ampère equation.
\item J-equation. 
\item Complex Hessian equation. 
\item Deformed Hermitian--Yang--Mills equation with supercritical phase.
\item Special Lagrangian equation with supercritical phase.
\item General inverse $\sigma_k$ equation with non-negative $c_k$ for $k \in \{0, \cdots, n-1\}$.
\end{itemize}
\end{frmk}

Throughout all these works in these equations, the convexity of either the equation itself or the level set plays a crucial role. To be more precise, analytically, to get a priori estimates, we highly rely on convexity. Since all these classical equations are strictly $\Upsilon$-stable, it is natural to ask the following question: Does the strictly $\Upsilon$-stability condition give some kind of convexity? In \cite{lin2023c}, the author confirmed this question and showed that if a connected component of a general inverse $\sigma_k$ multilinear polynomial is strictly $\Upsilon$-stable, then this connected component is convex. The next natural question will be the following: Given a $C$-subsolution to a strictly $\Upsilon$-stable general inverse $\sigma_k$ equation, does there exist a unique solution?
In this work, we show that if there exists a $C$-subsolution to a constant coefficients strictly $\Upsilon$-stable general inverse $\sigma_k$ equation, then there exists a unique solution. Our next goal is to explore the parabolic versions of these strictly $\Upsilon$-stable general inverse $\sigma_k$ equations which are interesting and still very wide open.
\smallskip

We also want to remark that equation (\ref{eq:1.1}) can be generalized further and plays an important role in reals. For example, Caffarelli--Nirenberg--Spruck \cite{caffarelli1985dirichlet}, Krylov \cite{krylov1993lectures, krylov1997fully}, Trudinger \cite{trudinger1995dirichlet}, and Cheng--Yau \cite{cheng1976regularity} on the study of the Dirichlet problem for the Hessian equations on various settings. Also, in convex geometry, the Christopher--Minkowski problem was studied extensively by Guan--Ma \cite{guan2003christoffel}, Guan--Lin--Ma \cite{guan2006christoffel}, and Guan--Zhang \cite{guan2019class}.
In the author's thesis \cite{lin2023thesis}, the author showed that we can extend these new real algebraic geometry techniques in \cite{lin2023c} to various settings.
Since we mainly focus on the general inverse $\sigma_k$ equations on Kähler manifolds in this work and because of the space limitations, the interested reader is referred to \cite{caffarelli2000priori, firey1967determination, firey1968christoffel, gilbarg2015elliptic, guan1999dirichlet, hou2010second, joyce2011existence, lewy1938differential, lu2022dirichlet, minkowski1897allgemeine, nirenberg1953weyl, pogorelov1978minkowski, siu2012lectures, wang2013singular} and the references therein.
\bigskip

Let us state some of our settings, definitions, and results now. First, for simplicity, by doing a substitution $X \coloneqq \chi - c_{n-1}\omega$ into equation (\ref{eq:1.1}), we get the following equation instead:
\begin{align*}
\label{eq:1.4}
X^n = d_{n-2} \mybinom[0.8]{n}{n-2}X^{n-2}  \wedge \omega^{2} + \cdots + d_{1} \mybinom[0.8]{n}{1} X \wedge \omega^{n-1} + d_{0} \mybinom[0.8]{n}{0}   \omega^{n}    =  \sum_{k = 0}^{n-2} d_k \mybinom[0.8]{n}{k} X^{k}  \wedge \omega^{n-k}, \tag{1.4}
\end{align*}where $d_k$ are real functions on $M$ obtained from this substitution. Be aware that $X$ might not be $d$-closed anymore by doing this substitution because $c_{n-1}$ is a real function on $M$. In this paper, we will only consider the case that $c_{n-1}$ is a constant, so $X$ is still $d$-closed. \smallskip

By collecting all strictly $\Upsilon$-stable general inverse $\sigma_k$ multilinear polynomial, we get the following subset of $\mathbb{R}^{n-1}$.
\hypertarget{D:1.2}{\begin{fdefi} 
Consider the following set $\tilde{\mathscr{C}}_n \subset \mathbb{R}^{n-1}$, which is defined by
\begin{align*}
\tilde{\mathscr{C}}_n \coloneqq \Bigl \{ (c_{n-2},  \cdots, c_{1}, c_{0}) \in \mathbb{R}^{n-1} \colon  \lambda_1 \cdots \lambda_n - \sum_{k = 0}^{n-2} c_k \sigma_k(\lambda) \text{ is strictly } \Upsilon\text{-stable}   \Bigr \}.
\end{align*}
\end{fdefi}}

If equation (\ref{eq:1.4}) is strictly $\Upsilon$-stable, then the coefficients $(d_{n-2}(z), \cdots, d_1(z), d_0(z)) \in \tilde{\mathscr{C}}_n$ for all $z \in M$. From now on, for convenience, we will always do this substitution to make all strictly $\Upsilon$-stable general inverse $\sigma_k$ equations look like format (\ref{eq:1.4}). Under this setting, we can view any strictly $\Upsilon$-stable general inverse $\sigma_k$ equation (\ref{eq:1.4}) as a map from $M$ to $\tilde{\mathscr{C}}_n$. That is to say, we have
\begin{align*}
\label{eq:1.5}
d \colon M \rightarrow \tilde{\mathscr{C}}_n \text{ with }  d(z) = (d_{n-2}(z), \cdots, d_1(z), d_0(z)) \in \tilde{\mathscr{C}}_n. \tag{1.5}
\end{align*}We are interested in whether there exists a representative in the cohomology class $[X]$ such that equation (\ref{eq:1.4}) (equivalently the map (\ref{eq:1.5})) is solvable. \smallskip

The complex Monge--Ampère equation can be reformulated into our settings, we view it as a map from $M$ to a subset of $\tilde{\mathscr{C}}_n$. The solvability of the complex Monge--Ampère equation, which is the Calabi conjecture proposed by Eugenio Calabi in \cite{calabi1954kahler, calabi1957kahler}, was solved by Shing-Tung~Yau in \cite{yau1978ricci}. To honor their remarkable contributions, we define the following subset of $\mathbb{R}^{n-1}$.

\hypertarget{D:1.6}{\begin{fdefi}
We define the Calabi--Yau set $\CY_n$ to be the following subset of $\tilde{\mathscr{C}}_n$:
\begin{align*}
\CY_n \coloneqq \bigl \{  (0, \cdots, 0, c_0) : c_0 > 0   \bigr\} \subset \tilde{\mathscr{C}}_n \subsetneq \mathbb{R}^{n-1}.
\end{align*}
\end{fdefi}}

\begin{fthm}[Reformulation of the complex Monge--Ampère equation, Yau \cite{yau1978ricci}]
Let $(M, \omega)$ be a compact connected Kähler manifold with Kähler form $\omega$ and $[X_0]$ be a $(1, 1)$-Dolbeault class. Given a map $d \colon M \rightarrow \CY_n \subset \tilde{\mathscr{C}}_n$ which is defined by
\begin{align*}
d \colon M \longrightarrow \CY_n \subset \tilde{\mathscr{C}}_n; \quad 
z \longmapsto (0, \cdots, 0, d_0(z))  
\end{align*}and satisfies the integrability condition $\int_M X_0^n = \int_M d_0 \omega^n$.
If there exists a $C$-subsolution (Kähler form in this case) to $d$ in $[X_0]$, then there exists a unique representative $X \in [X_0]$ such that 
\begin{align*}
X^n =  d_0 \omega^n.
\end{align*}
\end{fthm} 

\hypertarget{R:1.2}{\begin{frmk}
If $X_{\underline{u}} \coloneqq X + \sqrt{-1} \partial \bar{\partial} \underline{u}$ is a $C$-subsolution to $d \colon M \rightarrow \tilde{\mathscr{C}}_n$, then $X_{\underline{u}} > 0$ as a $(1, 1)$-form. That is, $X_{\underline{u}}$ is a Kähler form.
\end{frmk}}

For any $d \in \tilde{\mathscr{C}}_n$, we introduce the following subset $\tilde{\mathcal{P}}^d$ of $\mathbb{R}^{n-1}$ which plays an important role in obtaining a priori estimates. A more detailed explanation will be in Section~\ref{sec:2.3}.

\hypertarget{D:1.4}{\begin{fdefi}
Let $d  \in \tilde{\mathscr{C}}_{n}$, then we define the following polyhedron in $\tilde{\mathscr{C}}_{n}$:
    \begin{align*}
        \tilde{\mathcal{P}}^d \coloneqq \Bigl \{  {c} \in   \tilde{\mathscr{C}}_{n} \colon     x_l^{n-l}(d)  - \sum_{k=l}^{n-2} {c}_k\mybinom[0.8]{n-l}{k-l} x_l^{k-l}(d)  \geq 0, \quad \forall l \in \{1, \cdots, n-2\}   \Bigr \},  
    \end{align*}where $x_l(d)$ is the largest real root of the $l$-th derivative of $x^n - \sum_{k=0}^{n-2} d_k \binom{n}{k} x^k$. 
\end{fdefi}}

Given a strictly $\Upsilon$-stable general inverse $\sigma_k$ equation $d =  (d_{n-2}(z), \cdots, d_1(z), d_0(z) ) \colon M \rightarrow \tilde{\mathscr{C}}_n$.
For any strictly $\Upsilon$-stable general inverse $\sigma_k$ equation $c =  (c_{n-2}(z), \cdots, c_1(z), c_0(z) ) \colon M \rightarrow \tilde{\mathscr{C}}_n$, pointwise, if $(c_{n-2}(z), \cdots, c_1(z), c_0(z) ) \in \tilde{\mathcal{P}}^{d(z)}$ for any $z \in M$, then we say $c \in \tilde{\mathcal{P}}^d$ for convenience. We have the following important a priori estimates.

\hypertarget{T:1.2}{\begin{fthm}
Suppose $S$ is a compact subset of $\tilde{\mathscr{C}}_n$, $d \colon M^n \rightarrow  \tilde{\mathscr{C}}_n$ is a strictly $\Upsilon$-stable general inverse $\sigma_k$ equation with $d_1, \cdots, d_{n-2}$ constants and range in $S$, and $X$ is a $C$-subsolution to $d$. For any $c \in \tilde{\mathcal{P}}^d$ with $c_1, \cdots, c_{n-2}$ constants and range in $S$, if $u \colon M   \rightarrow \mathbb{R}$ is a smooth function solving $c \colon M^n \rightarrow \tilde{\mathscr{C}}_n$, then for any $\alpha \in (0, 1)$, there exists a constant $C$ such that 
\begin{align*}
\|\partial \bar{\partial} u \|_{C^{2, \alpha}} \leq C.
\end{align*}Here, $C = C(M, X, S, d, \omega, \alpha, \|c_0\|_{C^2})$ is a constant.
\end{fthm}}

In this paper, we study the solvability of the following equation:
\begin{align*}
\label{eq:1.6}
X^n = d_{n-2} \mybinom[0.8]{n}{n-2} X^{n-2}  \wedge \omega^{2} + \cdots   + d_1  \mybinom[0.8]{n}{1} X  \wedge \omega^{n-1} + d_{0}(z) \mybinom[0.8]{n}{0}   \omega^{n}, \tag{1.6}
\end{align*}where $d_1, \cdots, d_{n-2}$ are constants (not necessary non-negative), $(d_{n-2}, \cdots, d_1, d_0(z)) \in \tilde{\mathscr{C}}_n$ for any $z \in M$, and $\int_M X^n = \sum_{k=1}^{n-2}d_k \binom{n}{k} \int_M X^k \wedge \omega^{n-k} + \int_M d_0(z) \omega^n$. To obtain the solvability, we use the method of continuity. With Theorem~\hyperlink{T:1.2}{1.2}, for $t \in [0, 1]$, we consider continuity paths of the following format
\begin{align*}
\label{eq:1.7.t}
X^n = d_{n-2}(t) \mybinom[0.8]{n}{n-2} X^{n-2}  \wedge \omega^{2} + \cdots   + d_1(t) \mybinom[0.8]{n}{1} X  \wedge \omega^{n-1} + d_{0}(z, t) \mybinom[0.8]{n}{0}   \omega^{n}, \tag{1.7.$t$}
\end{align*}satisfying the following four constraints.
\begin{enumerate}[leftmargin=5.5cm]
	 \setlength\itemsep{0.6em}
\item[Topological constraint:] $\Omega_0 - \sum_{k=1}^{n-2} d_k(t) \Omega_{n-k} - \int_M d_0(z, t) \omega^n = 0$.
\item[Boundary constraints:]  $(d_{n-2}(1), \cdots, d_1(1), d_0(z, 1)) = (d_{n-2}, \cdots d_1, d_0(z))$  \\    
						and $(d_{n-2}(0), \cdots, d_1(0), d_0(z, 0)) \in \CY_n$ for any $z \in M$.
\item[Positivstellensatz constraint:] $(d_{n-2}(t), \cdots, d_1(t), d_0(z, t)) \in  \tilde{\mathscr{C}}_n$.
\item[$\Upsilon$-dominance constraint:] $(d_{n-2}(t), \cdots, d_1(t), d_0(z, t)) \in \tilde{\mathcal{P}}^d$.
\end{enumerate}Here, we denote $\Omega_i \coloneqq \int_M X^{n-i} \wedge \omega^i$. \smallskip

The idea of these constraints follows from the author's work \cite{lin2023d}. The first two constraints are not surprising.
First, the topological constraint is natural, along the continuity path (\ref{eq:1.7.t}), the integration over $M$ should be equal. 
Second, the boundary constraint is to make sure that one endpoint is the unsolved equation (\ref{eq:1.6}) and another endpoint is a well-known solvable one. So we wish to connect the unsolved equation to the complex Monge--Ampère equation, which is solvable by Yau \cite{yau1978ricci}. The last two constraints are not obvious and are related to the algebraic properties of $\Upsilon$-stable general inverse $\sigma_k$ polynomials. \smallskip

The Positivstellensatz constraint guarantees that these equations will not degenerate on the continuity path. From a real algebraic geometry viewpoint, if a general inverse $\sigma_k$ multilinear polynomial is $\Upsilon$-stable but not strictly $\Upsilon$-stable, then the level set degenerates. It is interesting to see whether it is possible to define a solution in this case and the regularity of this solution if such a solution can be defined.
If a general inverse $\sigma_k$ multilinear polynomial is not even $\Upsilon$-stable, then it is not likely that any arbitrary level set is convex and to have these fruitful algebraic properties coming from the $\Upsilon$-stability. \smallskip

The $\Upsilon$-dominance constraint ensures that we can apply a priori estimates on this continuity path (\ref{eq:1.7.t}). If $(d_{n-2}(t), \cdots, d_1(t), d_0(z, t)) \in \tilde{\mathcal{P}}^d$, then we will show in Section~\ref{sec:2} that $X$ is also a $C$-subsolution to $(d_{n-2}(t), \cdots, d_1(t), d_0(z, t))$. Hence, by Theorem~\hyperlink{T:1.2}{1.2}, we have a priori estimates for the solution to $(d_{n-2}(t), \cdots, d_1(t), d_0(z, t))$. If we have a priori estimates on this whole continuity path, then we can apply the Arzelà--Ascoli theorem and standard point-set topology argument and say that (\ref{eq:1.7.t}) is solvable for any $t \in [0, 1]$ provided that (\ref{eq:1.7.t}) is solvable when $t = 0$.\smallskip

In Section~\ref{sec:4}, we explicitly find a continuity path. That is to say, we have the following.
\hypertarget{T:1.3}{\begin{fthm}
There exists a continuity path $(d_{n-2}(t), \cdots, d_1(t), d_0(z, t))$ satisfying all the four constraints on (\ref{eq:1.7.t}).
\end{fthm}}

Due to the purpose of the four constraints on (\ref{eq:1.7.t}), if such a continuity path exists, then equation (\ref{eq:1.4}) is solvable. Hence, we finish the proof and obtain the following.

\begin{fthm}[Solvability]
Suppose that $[X]$ and $[\omega]$ satisfies the following integrability condition: 
\begin{align*}
\int_M X^n = d_{n-2} \int_M X^{n-2} \wedge \omega^2 + \cdots + d_1 \int_M X \wedge \omega^{n-1} + \int_M d_0(z) \omega^n.
\end{align*}with $d = (d_{n-2}, \cdots, d_1, d_0(z)) \colon M \rightarrow \tilde{\mathscr{C}}_n$ for any $z \in M$.
If there exists a $C$-subsolution to $d$, then there exists a unique representative $X_u \coloneqq X + \sqrt{-1} \partial \bar{\partial}u \in [X]$ such that
\begin{align*}
X_u^n = d_{n-2}  X_u^{n-2} \wedge \omega^2 + \cdots + d_1  X_u \wedge \omega^{n-1} +  d_0(z) \omega^n.
\end{align*}
\end{fthm}

As an application, we prove the analytical conjecture by Collins--Jacob--Yau in \cite{collins20201}.

\begin{fthm}[Solvability of the dHYM equation]
The analytical conjecture by Collins--Jacob--Yau in \cite{collins20201} is confirmed for all dimensions. That is, if there exists a $C$-subsolution to equation (\ref{eq:1.2}), then there exists a unique representative in the same cohomology class solving equation (\ref{eq:1.2}).
\end{fthm}

% Layout
The layout of this paper is as follows: 
in Section~\ref{sec:2}, we discuss some background materials. 
In Section~\ref{sec:2.1}, we introduce a special class of multivariate polynomials called $\Upsilon$-stable general inverse $\sigma_k$ multilinear polynomials.
We extend strictly $\Upsilon$-stable general inverse $\sigma_k$ multilinear polynomials to fully nonlinear elliptic partial differential equations on Kähler manifolds. We state some backgrounds and notations of strictly $\Upsilon$-stable general inverse $\sigma_k$ equations in this subsection. 
In Section~\ref{sec:2.2}, we introduce a special class of univariate polynomials called right-Noetherian polynomials and we state some new results of right-Noetherian polynomials. 
In \cite{lin2023c}, the author showed the correspondence between $\Upsilon$-stable general inverse $\sigma_k$ multilinear polynomials and right-Noetherian polynomials. Hence, new properties of right-Noetherian polynomials will lead to new properties of $\Upsilon$-stable general inverse $\sigma_k$ multilinear polynomials.
In Section~\ref{sec:2.3}, we revisit $\Upsilon$-stable general inverse $\sigma_k$ multilinear polynomials. Since we benefit from a better understanding of right-Noetherian polynomials in Section~\ref{sec:2.2}, we get some new properties of $\Upsilon$-stable general inverse $\sigma_k$ multilinear polynomials. For example, the space for the method of continuity. Roughly speaking, we know under which conditions a strictly $\Upsilon$-stable general inverse $\sigma_k$ equation has a priori estimates to its solution.
In Section~\ref{sec:3}, from a better understanding of $\Upsilon$-stable general inverse $\sigma_k$ multilinear polynomials in Section~\ref{sec:2}, we find a unified approach to get a priori estimates for strictly $\Upsilon$-stable general inverse $\sigma_k$ equations for all dimensions. 
In Section~\ref{sec:4}, we construct an explicit continuity path satisfying all four constraints and show that the existence of a $C$-subsolution to a strictly $\Upsilon$-stable general inverse $\sigma_k$ equation provides the solvability.

\Acknow

\section{Preliminaries}
\label{sec:2}
\subsection{Background and Notations of \texorpdfstring{$\Upsilon$}{}-stable General Inverse \texorpdfstring{$\sigma_k$}{} Equations}
\label{sec:2.1}
In \cite{lin2023c, lin2023thesis}, the author introduced strictly $\Upsilon$-stable general inverse $\sigma_k$ multilinear polynomials and showed that the level set of a strictly $\Upsilon$-stable general inverse $\sigma_k$ multilinear polynomial is convex. We state some definitions and important results in \cite{lin2023c, lin2023thesis}, and then we use these to state some basic settings and terminologies. This framework of strictly $\Upsilon$-stable general inverse $\sigma_k$ equations is useful and important for future study in analyzing equations with complicated mixed signs but a nice algebraic structure.\smallskip

First, let us state some conventions, see Spruck \cite{spruck2005geometric} for more details. For an $n$-tuple $\lambda = \{ \lambda_1, \cdots, \lambda_n\}$ and for $k \in \{1, \cdots, n \}$, the $k$-th elementary symmetric polynomial $\sigma_k(\lambda)$ of $\lambda$ will be
\begin{align*}
\sigma_k(\lambda) \coloneqq \sum_{1 \leq i_1 < \cdots < i_k \leq n} \lambda_{i_1} \cdots \lambda_{i_k}.
\end{align*}We also define $\sigma_0(\lambda) \coloneqq 1$ for convenience. For $l \in \{1, \cdots, n \}$ and pairwise distinct indices $i_1, \cdots, i_l$, where $i_j \in \{1, \cdots, n\}$ for all $j \in \{1, \cdots, l\}$, we denote the set $\lambda - \{ \lambda_{i_1}, \cdots, \lambda_{i_l} \}$ by $\lambda_{; i_1, \cdots, i_l}$. We call a multilinear polynomial having the following format
\begin{align*}
\lambda_1 \cdots \lambda_n -   \sum_{k = 0}^{n-1} c_k \sigma_k(\lambda).
\end{align*}a degree $n$ \emph{general inverse $\sigma_k$ multilinear polynomial}.\smallskip

Inspired by the work of Trudinger \cite{trudinger1995dirichlet} on the Dirichlet problem (over the reals) for equations of the eigenvalues of the Hessian, the results of Caffarelli--Nirenberg--Spruck \cite{caffarelli1985dirichlet}, and the results of Collins--Székelyhidi \cite{collins2017convergence}. In \cite{lin2023d}, the author introduced $\Upsilon$-cones to keep track of the information of the original equation as much as possible. We abstractly define the following sets.
\hypertarget{D:2.1}{\begin{fdefi}[\texorpdfstring{$\Upsilon$}{}-cones. Lin \cite{lin2023d, lin2023c}]
Let $f(\lambda) \coloneqq  \lambda_1 \cdots \lambda_n - \sum_{k = 0}^{n-1} c_k \sigma_k(\lambda)$ be a degree $n$ general inverse $\sigma_k$ multilinear polynomial and $\Gamma_f^{n}$ be a connected component of $\{f(\lambda) > 0\}$, we denote by $\Gamma^{n-1}_f$ the projection of $\Gamma^{n}_f$ onto $\mathbb{R}^{n-1}$ by dropping the last entry. We define  
\begin{align*}
\Upsilon_1 \coloneqq \bigl \{  \mu \in \mathbb{R}^n  \colon    \bigl(\mu_{s(1)}, \cdots, \mu_{s(n-1)}\bigr) \in \Gamma^{n-1}_f,\quad \forall s \in S_n \bigr \},  
\end{align*}where $S_n$ is the symmetric group. For $n-1 \geq k \geq 2$, we define the following $\Upsilon$-cones
\begin{align*}
\Upsilon_k  \coloneqq \bigl \{  \mu \in \mathbb{R}^n  \colon    \bigl(\mu_{s(1)}, \cdots, \mu_{s(n-k)}\bigr) \in \Gamma^{n-k}_f,\quad \forall s \in S_n  \bigr \},  
\end{align*}where we define $ \Gamma^{n-k}_f$ inductively by the projection of $ \Gamma^{n+1-k}_f$ onto $\mathbb{R}^{n-k}$ by dropping the last entry. Here, $k$ increases from $2$ to $n-1$.
\end{fdefi}}

Based on the classical equations, we are interested in whether there exists a connected component of $\{f(\lambda) > 0\}$ contained in the positive orthant $\Gamma_n$ after translation. We define the following stableness condition for general inverse $\sigma_k$ multilinear polynomials.
\hypertarget{D:2.2}{\begin{fdefi}[\texorpdfstring{$\Upsilon$}{}-stableness, L.~\cite{lin2023c}]
Let $f(\lambda) \coloneqq   \lambda_1 \cdots \lambda_n - \sum_{k = 0}^{n-1} c_k \sigma_k(\lambda)$ be a general inverse $\sigma_k$ multilinear polynomial and $\Gamma_f^{n}$ be a connected component of $\{f(\lambda) > 0\}$. We say that this connected component $\Gamma^n_f$ of $f(\lambda)$ is $\Upsilon$-stable if
\begin{align*}
\Gamma^n_f \subseteq  q +  \Gamma_n \text{ for some } q \in \mathbb{R}^n, 
\end{align*}where $\Gamma_n$ is the positive orthant of $\mathbb{R}^n$. We say that this connected component $\Gamma^n_f$ is strictly $\Upsilon$-stable if it is $\Upsilon$-stable and the boundary $\partial \Gamma^n_f$ is contained in the $\Upsilon_1$-cone. 
\end{fdefi}} 

Throughout this paper, if a general inverse $\sigma_k$ multilinear polynomial has a (strictly) $\Upsilon$-stable connected component, then we also call this polynomial a (strictly) $\Upsilon$-stable general inverse $\sigma_k$ multilinear polynomial for convenience. We hope that no confusion will result.  

\hypertarget{R:2.1}{\begin{frmk}
Let $f(\lambda) \coloneqq   \lambda_1 \cdots \lambda_n - \sum_{k = 0}^{n-1} c_k \sigma_k(\lambda)$ be a strictly $\Upsilon$-stable general inverse $\sigma_k$ multilinear polynomial. Then the $\Upsilon_1$-cone of $f(\lambda)$ is exactly the tangent cone at infinity introduced by Guan \cite{guan2014second} and the $C$-subsolution cone introduced by Székelyhidi \cite{szekelyhidi2018fully}.
\end{frmk}}

\hypertarget{L:2.1}{\begin{flemma}[L.~\cite{lin2023c}]
Let $f(\lambda) \coloneqq  \lambda_1 \cdots \lambda_n - \sum_{k = 0}^{n-1} c_k \sigma_k(\lambda)$ be an $\Upsilon$-stable general inverse $\sigma_k$ multilinear polynomial and say $\Gamma_f^{n}  \subseteq  q +  \Gamma_n$ for some $q = (q_1, \cdots, q_n)$, then 
\begin{align*}
\label{eq:2.1}
\Gamma^n_f \subseteq  \Upsilon_1 \subseteq \Upsilon_2  \subseteq \cdots \subseteq \Upsilon_{n-1} = (c_{n-1}, \cdots, c_{n-1}) + \Gamma_n, \tag{2.1}
\end{align*}$c_{n-1} \geq q_i$ for all $i \in \{1, \cdots, n\}$. For any $l \in \{1, \cdots, n-1\}$, we have $\Upsilon_l$ is open, connected, and
%\begin{align*}
%\Upsilon_l  =      \bigcap_{1 \leq i_1 < \cdots < i_l \leq n} \Gamma^n_{\frac{\partial^l f}{\partial \lambda_{i_1} \cdots \partial \lambda_{i_l}}} =  \bigcap_{1 \leq i_1 < \cdots < i_l \leq n} \Bigl \{   \sigma_{n-l} (\lambda_{; i_1, \cdots, i_l})   - \sum_{k=l}^{n-2} c_k \sigma_{k - l}(\lambda_{;i_1, \cdots, i_l})  > 0  \Bigr \}.
%\end{align*}
\begin{align*}
\label{eq:2.2}
\Upsilon_l  =      \bigcap_{1 \leq i_1 < \cdots < i_l \leq n} \Gamma^n_{ f_{i_1 \cdots i_l} } =  \bigcap_{1 \leq i_1 < \cdots < i_l \leq n} \Bigl \{   \sigma_{n-l} (\lambda_{; i_1, \cdots, i_l})   - \sum_{k=l}^{n-1} c_k \sigma_{k - l}(\lambda_{;i_1, \cdots, i_l})  > 0  \Bigr \}. \tag{2.2}
\end{align*}Here, we write $f_{i_1 \cdots i_l}$ as the $l$-th partial derivative $\frac{\partial^l f}{\partial \lambda_{i_1} \cdots \partial \lambda_{i_l}}$ for $l \in \{1, \cdots, n-1\}$.
\end{flemma}}

One of the most important results in \cite{lin2023c} is the strictly $\Upsilon$-stableness implies the level set convexity (of course need to specify the $\Upsilon$-stable connected component), which is the following.

\hypertarget{T:2.1}{\begin{fthm}[Convexity of a general inverse $\sigma_k$ multilinear polynomial, L.~\cite{lin2023c}]
If $f(\lambda)$ is a strictly $\Upsilon$-stable general inverse $\sigma_k$ multilinear polynomial, then the level set $\{ f = 0\}$ is convex.
\end{fthm}}

In the author's thesis \cite{lin2023thesis}, by collecting all strictly $\Upsilon$-stable general inverse $\sigma_k$ multilinear polynomials, we get the following algebraic set in the Euclidean space (the author also calls it the moduli space of strictly $\Upsilon$-stable general inverse $\sigma_k$ multilinear polynomials).
%We can also reformulate all classical strictly $\Upsilon$-stable general inverse $\sigma_k$ equations into our new setting.

\hypertarget{D:2.3}{\begin{fdefi}[L.~\cite{lin2023thesis}]
Let $\lambda = \{\lambda_1, \cdots, \lambda_n\}$. We define the following subsets in the Euclidean space endowed with the standard Euclidean topology:
\begin{align*}
\label{eq:2.3}
\mathscr{C}_n &\coloneqq \Bigl \{   (c_{n-1}, \cdots, c_1, c_0) \in \mathbb{R}^n \colon \lambda_1 \cdots \lambda_n - \sum_{k = 0}^{n-1} c_k \sigma_k(\lambda) \text{ is strictly } \Upsilon\text{-stable}    \Bigr \}; \tag{2.3} \\
\label{eq:2.4}
\tilde{\mathscr{C}}_n &\coloneqq \Bigl \{   (c_{n-2}, \cdots, c_1, c_0) \in \mathbb{R}^{n-1} \colon \lambda_1 \cdots \lambda_n - \sum_{k = 0}^{n-2} c_k \sigma_k(\lambda) \text{ is strictly } \Upsilon\text{-stable}    \Bigr \}. \tag{2.4}
\end{align*}Here, we let $\mathscr{C}_n$ and $\tilde{\mathscr{C}}_n$ be topological spaces using the subspace topology induced from the standard Euclidean topology of the Euclidean space.
\end{fdefi}}

Let $(M, \omega)$ be a compact connected Kähler manifold of complex dimension $n$ with a Kähler form $\omega$ and $[\chi_0] \in H^{1,1}(M; \mathbb{R})$, where $H^{1,1}(M; \mathbb{R})$ is the $(1, 1)$-Dolbeault cohomology group. The study of the solvability of the following equation is widely considered:
\begin{align*}
\label{eq:2.5}
\chi^n = c_{n-1} \mybinom[0.8]{n}{n-1} \chi^{n-1}  \wedge \omega  + \cdots + c_{1} \mybinom[0.8]{n}{1} \chi   \wedge \omega^{n-1} + c_{0} \mybinom[0.8]{n}{0}   \omega^{n}    =  \sum_{k = 0}^{n-1} c_k \mybinom[0.8]{n}{k} \chi^{k}  \wedge \omega^{n-k}, \tag{2.5}
\end{align*}where $c_k$ are real functions on $M$ and $\chi \in [\chi_0]$ is a real smooth, closed $(1, 1)$-form. We call an equation having the same format as equation (\ref{eq:2.5}) a degree $n$ \emph{general inverse $\sigma_k$ equation}. At $z \in M$, if we write equation (\ref{eq:2.5}) in terms of the eigenvalues of the Hermitian endomorphism $\Lambda = \omega^{-1} \chi$, then we can rewrite equation (\ref{eq:2.5}) at point $z$ as
\begin{align*}
\lambda_1 \cdots \lambda_n = \sum_{k = 0}^{n-1} c_k(z) \sigma_k(\lambda).
\end{align*}The above expression corresponds to the general inverse $\sigma_k$ multilinear polynomial $\lambda_1 \cdots \lambda_n - \sum_{k = 0}^{n-1} c_k(z) \sigma_k(\lambda)$. If the coefficients $(c_{n-1}(z), \cdots, c_1(z), c_0(z)) \in \mathscr{C}_n$ for all $z \in M$, then we call equation (\ref{eq:2.5}) a \emph{strictly $\Upsilon$-stable general inverse $\sigma_k$ equation}. \smallskip

Similarly, by doing a substitution $X \coloneqq \chi - c_{n-1}\omega$, we get the following equation instead:
\begin{align*}
\label{eq:2.6}
X^n = d_{n-2} \mybinom[0.8]{n}{n-2}X^{n-2}  \wedge \omega^{2} + \cdots + d_{1} \mybinom[0.8]{n}{1} X \wedge \omega^{n-1} + d_{0} \mybinom[0.8]{n}{0}   \omega^{n}    =  \sum_{k = 0}^{n-2} d_k \mybinom[0.8]{n}{k} X^{k}  \wedge \omega^{n-k}, \tag{2.6}
\end{align*}where $d_k$ are real functions on $M$ obtained from the substitution. Be aware that $X$ might not be $d$-closed anymore by doing this substitution because $c_{n-1}$ is a real function on $M$. In this paper, we will only consider the case that $c_{n-1}$ is a constant, so $X$ is still in a $(1, 1)$-Dolbeault cohomology class. If equation (\ref{eq:2.6}) is strictly $\Upsilon$-stable, then the coefficients $(d_{n-2}(z), \cdots, d_1(z), d_0(z)) \in \tilde{\mathscr{C}}_n$ for all $z \in M$. From now on, for convenience, we will always do this substitution to make all equations look like format (\ref{eq:2.6}). So, under this setting, we can view any strictly $\Upsilon$-stable general inverse $\sigma_k$ equation (\ref{eq:2.6}) as a map from $M$ to $\tilde{\mathscr{C}}_n$. That is to say, we have
\begin{align*}
\label{eq:2.7}
d \colon M \rightarrow \tilde{\mathscr{C}}_n \text{ with }  d(z) = (d_{n-2}(z), \cdots, d_1(z), d_0(z)) \in \tilde{\mathscr{C}}_n. \tag{2.7}
\end{align*}And we are interested in whether there exists a representative in the cohomology class $[X]$ such that equation (\ref{eq:2.6}) (equivalently the map (\ref{eq:2.7})) is solvable. \smallskip

The complex Monge--Ampère equation can also be reformulated into our settings, we can view it as a map from $M$ to a subset of $\tilde{\mathscr{C}}_n$. The solvability of the complex Monge--Ampère equation, which is the Calabi conjecture proposed by Eugenio Calabi in \cite{calabi1954kahler, calabi1957kahler}, was solved by Shing-Tung~Yau in \cite{yau1978ricci}. To honor their remarkable contributions, we define the following set.

\hypertarget{D:2.4}{\begin{fdefi}
We define the Calabi--Yau set $\CY_n$ to be the following subset of $\tilde{\mathscr{C}}_n$:
\begin{align*}
\CY_n \coloneqq \bigl \{  (0, \cdots, 0, c_0) : c_0 > 0   \bigr\} \subset \tilde{\mathscr{C}}_n \subsetneq \mathbb{R}^{n-1}.
\end{align*}
\end{fdefi}}

We state an alternative definition of $C$-subsolution to equation (\ref{eq:2.6}) (or map (\ref{eq:2.7})), which was introduced by Guan \cite{guan2014second} and Székelyhidi \cite{szekelyhidi2018fully}, so it is simpler for us to operate in our settings.
\hypertarget{D:2.5}{\begin{fdefi}[$C$-subsolution]
Let $(M, \omega)$ be a compact connected Kähler manifold of complex dimension $n$ with a Kähler form $\omega$ and $[X_0] \in H^{1,1}(M; \mathbb{R})$. We say a smooth function $\underline{u} \colon M \rightarrow \mathbb{R}$ is a $C$-subsolution to equation (\ref{eq:2.6}) (or map (\ref{eq:2.7})) if for any $z \in M$, we always have
\begin{align*}
X_{\underline{u}}(z) \coloneqq \bigl( X_0 + \sqrt{-1} \partial \bar{\partial} \underline{u} \bigr) (z) \in \Upsilon_1^c(z),
\end{align*}where $\Upsilon_1^c(z)$ is the $\Upsilon_1$-cone of $\lambda_1 \cdots \lambda_n - \sum_{k=0}^{n-2} c_k(z) \sigma_k(\lambda)$.
\end{fdefi}}

Throughout this paper, for notational conventions, if $(c_{n-2}, \cdots, c_1, c_0) \in \tilde{\mathscr{C}}_n$, then we denote $(c_{n-2}, \cdots, c_1, c_0)$ by $c$. If $\underline{u} \colon M \rightarrow \mathbb{R}$ is a $C$-subsolution to equation (\ref{eq:2.6}) (or map (\ref{eq:2.7})), then we also call the $(1, 1)$-form $X_{\underline{u}} \coloneqq  X_0 + \sqrt{-1} \partial \bar{\partial} \underline{u}$ a $C$-subsolution to equation (\ref{eq:2.6}) (or map (\ref{eq:2.7})). We hope that there will be no confusion.

\subsection{Notations and Properties of Right-Noetherian Polynomials}
\label{sec:2.2}
In Section~\ref{sec:2.1}, we have seen the importance of $\Upsilon$-stable general inverse $\sigma_k$ multilinear polynomials. However, in practice, it is difficult to verify whether a level set of a general inverse $\sigma_k$ multilinear polynomial after translation is contained in the positive orthant. In \cite{lin2023c}, the author found a way to reduce this difficult problem related to special multivariate polynomials to a simpler problem related to special univariate polynomials. We state some definitions and important results of these special univariate polynomials in \cite{lin2023c, lin2023thesis}, and then we use these to obtain some new properties. 

\hypertarget{D:2.6}{\begin{fdefi}[Noetherian polynomial, L.~\cite{lin2023c}]
We say a degree $n$ real univariate polynomial $p(x)$ is right-Noetherian if for all $l \in \{0, \cdots, n-2 \}$, there exists a real root of $p^{(l)}(x)$, the $l$-th derivative of $p$, which is greater than or equal to the largest real root of $p^{(l+1)}(x)$, the $(l+1)$-th derivative of $p$. We say a right-Noetherian polynomial $p(x)$ is strictly right-Noetherian if the largest real root of $p(x)$ is strictly greater than the largest real root of $p'(x)$, the first derivative of $p$.
\end{fdefi}}

\hypertarget{P:2.1}{\begin{fprop}[L.~\cite{lin2023c}]
Let $p(x)$ be a right-Noetherian polynomial of degree $n$. Then for any $l \in \{0, \cdots, n-2\}$, there exists a unique (ignoring multiplicity) real root of $p^{(l)}(x)$ which is greater than or equal to the largest real root of $p^{(l+1)}(x)$. Moreover, this real root is the largest real root of $p^{(l)}(x)$. In particular, if we denote $x_l$ the largest real root of $p^{(l)}(x)$ for $l \in \{0, \cdots, n-2\}$, then
\begin{align*}
x_0 \geq x_1 \geq \cdots \geq x_{n-1}.
\end{align*}
\end{fprop}}

\hypertarget{T:2.2}{\begin{fthm}[Positivstellensatz, L.~\cite{lin2023c}]
Let $f(\lambda) \coloneqq  \lambda_1 \cdots \lambda_n - \sum_{k = 0}^{n-1} c_k \sigma_k(\lambda)$ be a general inverse $\sigma_k$ multilinear polynomial. Then $f(\lambda)$ is $\Upsilon$-stable if and only if the diagonal restriction $r_f(x)$ of $f(\lambda)$, which is defined by the following  
\begin{align*}
r_f(x) \coloneqq f(x, \cdots, x) = x^n - \sum_{k=0}^{n-1} c_k \mybinom[0.8]{n}{k} x^k,  
\end{align*} is right-Noetherian. Moreover, $\Gamma^n_f$ is strictly $\Upsilon$-stable if and only if $r_f$ is strictly right-Noetherian.
\end{fthm}}

A quick consequence of Theorem~\hyperlink{T:2.2}{2.2} will be the explicit description of $\tilde{\mathscr{C}}_n$ when the degree $n$ is low. The interested reader is referred to the author's previous works \cite{lin2023d, lin2023c, lin2023thesis}. \smallskip

Given $c \in \tilde{\mathscr{C}}_n$, which corresponds to the following strictly $\Upsilon$-stable general inverse $\sigma_k$ multilinear polynomial $\lambda_1 \cdots \lambda_n - \sum_{k=0}^{n-2} c_k \sigma_k(\lambda)$. By Theorem~\hyperlink{T:2.2}{2.2}, equivalently the diagonal restriction $x^n - \sum_{k=0}^{n-2} c_k \binom{n}{k} x^k$ is strictly right-Noetherian. If we denote $x_l(c)$ the largest real root of the $l$-th derivative of $x^n - \sum_{k=0}^{n-2} c_k \binom{n}{k} x^k$ for $l \in \{0, \cdots, n-1\}$, then by Proposition~\hyperlink{P:2.1}{2.1} we have
\begin{align*}
x_0(c) > x_1(c) \geq x_2(c) \geq \cdots \geq x_{n-2}(c) \geq x_{n-1}(c) = 0.
\end{align*}
So, any $c \in \tilde{\mathscr{C}}_n$ gives us an $(n-1)$-tuple $(x_{n-2}(c), \cdots, x_1(c), x_0(c))$ in the following polyhedron.

\hypertarget{D:2.7}{\begin{fdefi}
Let $\tilde{\mathcal{X}}_n$ to be the following polyhedron in $\mathbb{R}^{n-1}$:
\begin{align*}
%\tilde{\mathcal{X}}_n \coloneqq \bigl \{   x_0 > x_1 \geq x_2 \geq \cdots \geq x_{n-2} \geq 0    \bigr\} \subsetneq \mathbb{R}^{n-1}.
\tilde{\mathcal{X}}_n \coloneqq \bigl \{  (x_{n-2}, \cdots, x_1, x_0) \colon    x_0 > x_1 \geq x_2 \geq \cdots \geq x_{n-2} \geq 0     \bigr\} \subsetneq \mathbb{R}^{n-1}.
\end{align*}Here, we let $\tilde{\mathcal{X}}_n$ be a topological space using the subspace topology induced from the standard Euclidean topology of the Euclidean space.
\end{fdefi}}

We can treat the above correspondence as a map $\varphi$ from $\tilde{\mathscr{C}}_n$ to $\tilde{\mathcal{X}}_n$ sending $c = (c_{n-2}, \cdots, c_1, c_0) \in \tilde{\mathscr{C}}_n$ to $(x_{n-2}(c), \cdots, x_1(c), x_0(c)) \in \tilde{\mathcal{X}}_n$. By the argument principle, we know that the roots of a polynomial depend continuously on its coefficients. But here, since we only focus on the largest real roots of all the derivatives of the diagonal restriction, it is not obvious whether these largest roots will still depend continuously on the coefficients. Indeed, in the following Lemma, we will show that $\varphi$ is a continuous map and furthermore determines a homeomorphism from $\tilde{\mathscr{C}}_n$ to $\tilde{\mathcal{X}}_n$. 
\hypertarget{L:2.2}{\begin{flemma}
The map $\varphi \colon \tilde{\mathcal{X}}_n \rightarrow \tilde{\mathscr{C}}_n$ is a homeomorphism.
\end{flemma}}
\begin{proof}
First, we consider the following map $\varphi \colon \tilde{\mathscr{C}}_n \rightarrow \tilde{\mathcal{X}}_n$ defined by sending $c = (c_{n-2}, \cdots, c_1, c_0) \in  \tilde{\mathscr{C}}_n$ to $(x_{n-2}(c), \cdots, x_1(c), x_0(c)) \in \tilde{\mathcal{X}}_n$, where $x_l$ is the largest real root of the $l$-th derivative of $x^n - \sum_{k=0}^{n-2} c_k \binom{n}{k} x^k$ for $l \in \{0, \cdots, n-2\}$. To check whether $\varphi$ is injective, suppose $\varphi(c) = \varphi(\tilde{c})$ for $c, \tilde{c} \in  \tilde{\mathscr{C}}_n$. We immediately get $c_{n-2} = x_{n-2}^2(c) = x_{n-2}^2(\tilde{c}) = \tilde{c}_{n-2}$. Similarly, we get $c_{n-3} = x_{n-3}^3(c) - 3 c_{n-2} x_{n-3}(c) = x_{n-3}^3(\tilde{c}) - 3 \tilde{c}_{n-2} x_{n-3}(\tilde{c}) = \tilde{c}_{n-3}$. It should be straightforward to obtain that $c = \tilde{c}$ by mathematical induction. Hence, the map $\varphi$ is injective. \smallskip

Second, we consider the following map  $\psi \colon  \tilde{\mathcal{X}}_n \rightarrow \mathbb{R}^{n-1}$ sending
$(x_{n-2}, \cdots, x_1, x_0) \in \tilde{\mathcal{X}}_n$ to $(c_{n-2}, \cdots, c_1, c_0)$, where $c_l$ for $l \in \{0, 1, \cdots, n-2\}$ is defined recursively by
\begin{align*}
\label{eq:2.8.l}
c_l \coloneqq x_l^{n-l} - \sum_{k = l +1}^{n-2} c_k \mybinom[0.8]{n-l}{k-l} x_l^{k-l} \tag{2.8.$l$}
\end{align*}from $n-2$ back to $0$. To show that $(c_{n-2}, \cdots, c_1, c_0) \in \tilde{\mathscr{C}}_n$, we have $x_{n-2}^2 - c_{n-2} = 0$ so $x_{n-2}$ is the largest real root of the $(n-2)$-th derivative of $x^n - \sum_{k=0}^{n-2} c_k \binom{n}{k} x^k$. Similarly, since $x_{n-3} \geq x_{n-2}$, $x_{n-3}^3 - 3c_{n-2} x_{n-3} - c_{n-3} =0$, and by Proposition~\hyperlink{P:2.1}{2.1}, $x_{n-3}$ is the largest real root of the $(n-3)$-th derivative of $x^n - \sum_{k=0}^{n-2} c_k \binom{n}{k} x^k$. It should be straightforward to show that $(c_{n-2}, \cdots, c_1, c_0) \in \tilde{\mathscr{C}}_n$ and $\varphi(c) = (x_{n-2}, \cdots, x_1, x_0)$ by mathematical induction. This implies that $\varphi$ is surjective. \smallskip

Then, we check that $\varphi \circ \psi = \Id_{\tilde{\mathcal{X}}_n}$ and $ \psi \circ \varphi = \Id_{\tilde{\mathscr{C}}_n}$. For $(x_{n-2}, \cdots, x_1, x_0) \in \tilde{\mathcal{X}}_n$, we have $\psi (x_{n-2}, \cdots, x_1, x_0) = (c_{n-2}, \cdots, c_1, c_0)$, where $c_l$ for $l \in \{0, 1, \cdots, n-2\}$ is inductively defined by (\ref{eq:2.8.l}) from $n-2$ back to $0$. Consider $\varphi(c) = (x_{n-2}(c), \cdots, x_1(c), x_0(c))$, we immediately get $x_{n-2}^2(c) = c_{n-2} = x^2_{n-2}$, which implies that $x_{n-2}(c) = x_{n-2}$. Similar to above, by Proposition~\hyperlink{P:2.1}{2.1}, we have $(x_{n-2}(c), \cdots, x_1(c), x_0(c)) = (x_{n-2}, \cdots, x_1, x_0)$. This implies that $\varphi \circ \psi = \Id_{\tilde{\mathcal{X}}_n}$. On the other hand, for $c \in \tilde{\mathcal{X}}_n$, we have $\varphi(c) = (x_{n-2}(c), \cdots, x_1(c), x_0(c))$ and it is straightforward to check that $\psi \circ \varphi = \Id_{\tilde{\mathscr{C}}_n}$. Hence, the map $\varphi \colon \tilde{\mathscr{C}}_n \rightarrow \tilde{\mathcal{X}}_n$ is a bijection and its inverse $\psi \colon  \tilde{\mathcal{X}}_n \rightarrow \tilde{\mathscr{C}}_n$ is a continuous map due to the fact that $\psi$ is determined by polynomials.

Last, we check that $\varphi \colon \tilde{\mathscr{C}}_n \rightarrow \tilde{\mathcal{X}}_n$ is also continuous. Fix $c \in  \tilde{\mathscr{C}}_n$, by collecting all the roots of the $l$-th derivative of $x^n - \sum_{k=0}^{n-2} c_k \binom{n}{k} x^k$ for all $l$ from $0$ to $n-1$, we may pick an $\epsilon > 0$ sufficiently small such that the closed balls with radius $\epsilon > 0$ centered at all these roots are disjoint when the roots are distinct. Moreover, by shrinking $\epsilon > 0$ if necessary, the closed balls with radius $\epsilon > 0$ centered at all the complex roots do not intersect the real line. By the argument principle, there exists a $\delta > 0$ such that for any $\tilde{c} \in \tilde{\mathscr{C}}_n$ with $\| c- \tilde{c}\| < \delta$ and for any $l \in \{0, \cdots, n-1\}$, if $\xi$ is a root of the $l$-th derivative of $x^n - \sum_{k=0}^{n-2} c_k \binom{n}{k} x^k$ with multiplicity $m$, then there are $m$ roots (counting multiplicities) in the open ball with radius $\epsilon > 0$ centered at $\xi$. We claim that $|x_l(c) - {x}_l(\tilde{c})| < \epsilon$ for any $l \in \{0, \cdots, n-2\}$, which implies that $\varphi$ is a continuous map. For $l = n-2$, we have $x_{n-2}(c) = \sqrt{c_{n-2}}$ and ${x}_{n-2}(\tilde{c}) = \sqrt{\tilde{c}_{n-2}}$. If $c_{n-2} = 0$, then by the argument principle, all two roots of $x^2 - \tilde{c}_{n-2}$ must be in the open ball with radius $\epsilon > 0$ centered at $x_{n-2}(c) = 0$. Hence, the largest real root ${x}_{n-2}(\tilde{c}) = \sqrt{\tilde{c}_{n-2}}$ must be in the open ball with radius $\epsilon > 0$ centered at $x_{n-2}(c) = 0$. If $c_{n-2} > 0$, then by the choice of $\epsilon$ and the argument principle, ${x}_{n-2}(\tilde{c})$ must be in the open ball with radius $\epsilon > 0$ centered at $x_{n-2}(c)$. Similarly, we use mathematical induction on the subscript from $n-2$ to $0$. Suppose the statement is true when $k = m > 0$. When $k = m-1$, if $x_{m-1}(c) > x_m(c)$, since $x^n - \sum_{k=0}^{n-2} \tilde{c}_k \binom{n}{k} x^k$ is a real polynomial, there exists a unique real root of the $(m-1)$-th derivative of $x^n - \sum_{k=0}^{n-2} \tilde{c}_k \binom{n}{k} x^k$ in the open ball with radius $\epsilon > 0$ centered at $x_{m-1}(c)$. In addition, since $x_{m-1}(c) > x_m(c) + \epsilon > {x}_m(\tilde{c})$ and by Proposition~\hyperlink{P:2.1}{2.1}, ${x}_{m-1}(\tilde{c})$ must be this real root in the open ball with radius $\epsilon > 0$ centered at $x_{m-1}(c)$. If $x_{m-1}(c) = x_m(c)$, then by the choice of $\epsilon$, ${x}_{m-1}(\tilde{c})$ must be in the open ball with radius $\epsilon > 0$ centered at a real root of the $(m-1)$-th derivative of $x^n - \sum_{k=0}^{n-2} c_k \binom{n}{k} x^k$. If this real root is less than $x_{m-1}(c)$, then we get ${x}_{m}(\tilde{c}) > x_{m}(c) - \epsilon > {x}_{m-1}(\tilde{c})$ which gives a contradiction. So this real root of the $(m-1)$-th derivative of $x^n - \sum_{k=0}^{n-2} c_k \binom{n}{k} x^k$ must be greater than or equal to $x_{m-1}(c)$. By Proposition~\hyperlink{P:2.1}{2.1}, $x_{m-1}(c)$ is the largest real root of the $(m-1)$-th derivative of $x^n - \sum_{k=0}^{n-2} c_k \binom{n}{k} x^k$. Hence, ${x}_m(\tilde{c})$ lies in the open ball with radius $\epsilon$ centered at $x_m(c)$. This finishes the proof.
\end{proof}

As a quick consequence, we have the following Corollary. 
\hypertarget{Cor:2.1}{\begin{fcor}
$\tilde{\mathscr{C}}_n$ can be decomposed into disjoint union of $2^{n-2}$ stratas, where for any $l \in \{1, \cdots, n-1\}$ there are $\binom{n-2}{l-1}$ stratas with dimension $l$.
\end{fcor}}
\begin{proof}
By Lemma~\hyperlink{L:2.2}{2.2}, the number of stratas in $\tilde{\mathscr{C}}_n$ corresponds to the number of stratas in $\tilde{\mathcal{X}}_n$. So, $\tilde{\mathscr{C}}_n$ can be decomposed into disjoint union of $2^{n-2}$ stratas of different dimensions. Moreover, any equality symbol in $\tilde{\mathcal{X}}_n$ reduces the dimension of the corresponding strata by one. Hence, for any $l \in \{1, \cdots, n-1\}$ there are $\binom{n-2}{l-1}$ stratas with dimension $l$.
\end{proof}

\hypertarget{R:2.2}{\begin{frmk}
The unique $1$-dimensional strata in $\tilde{\mathscr{C}}_n$ will be $\CY_n \coloneqq \bigl \{  (0, \cdots, 0, c_0) : c_0 > 0   \bigr\}$ defined in Definition~\hyperlink{D:2.4}{2.4}. Besides, the unique $(n-1)$-dimensional strata in $\tilde{\mathscr{C}}_n$ corresponds to the set $ \{  (x_{n-2}, \cdots, x_1, x_0) \colon    x_0 >   \cdots > x_{n-2} > 0    \}$ in $\tilde{\mathcal{X}}_n$ and we call it the generic strata of $\tilde{\mathscr{C}}_n$. The generic strata consists of many important strictly $\Upsilon$-stable general inverse $\sigma_k$ equations, for example, the $J$-equation and the deformed Hermitian--Yang--Mills equation.
\end{frmk}}

\hypertarget{R:2.3}{\begin{frmk}
In fact, $\varphi$ can be extended to a map from $\overline{\tilde{\mathscr{C}}_n}$, the closure of $\tilde{\mathscr{C}}_n$, to $\overline{\tilde{\mathcal{X}}_n}$, the closure of $\tilde{\mathcal{X}}_n$. Moreover, similar to the proof of Lemma~\hyperlink{L:2.2}{2.2}, the extension is still a homeomorphism.
\end{frmk}}

In Figure~\ref{fig:2.1}, we plot the homeomorphism between $\tilde{\mathcal{X}}_3$ and $\tilde{\mathscr{C}}_3$. Straight rays in $\tilde{\mathcal{X}}_3$ will be sent to bent curves in $\tilde{\mathscr{C}}_3$ as the graph shows.
From now on, for notational conventions, if only one $c \in \tilde{\mathscr{C}}_n$ (or $\overline{\tilde{\mathscr{C}}_n}$) shows up in a context, then $\varphi(c) = (x_{n-2}(c), \cdots, x_1(c), x_0(c))$ is abbreviated as $\varphi(c) = (x_{n-2}, \cdots, x_1, x_0)$. We hope that this will not cause any confusion.

\begin{figure}
\centering
\begin{tikzpicture} 
\node (a) at (0,0)
{
\begin{tikzpicture}[scale=0.4]
\begin{scope}
    \fill[color=red!50,opacity=0.3,thick,domain=0:9.6,samples=101]
    (0, 9.6) -- plot (\x, \x ) ;
  \end{scope}
\draw[->] (-0.2,0) -- (9.6,0) node[right] {$x_1$};
\draw[->] (0,-0.2) -- (0,9.6) node[above] {$x_0$};
\foreach \t in {4, ..., 12} %n = 16 
{
	\pgfmathsetmacro\n{ 16 }
	\pgfmathsetmacro\R{(186-  186*\t/\n)/256 }
	\pgfmathsetmacro\G{(12 + 88*\t/\n)/256 }	
	\pgfmathsetmacro\B{(47 + 117*\t/\n)/256}
  	\draw[color={rgb, 1:red, \R; green, \G; blue, \B},domain=0:9.6,samples=100]    plot ( {\t*\x/\n} ,\x)             ;
}
 \draw[color=OSUR,thick,domain=0.1:9.5]    plot ({0},\x)             ;
\draw[color=UCIB,dashed,thick,domain=0.1:9.6]    plot (\x, \x)             ;
 \draw [color=red,fill=white] (0, 0) circle[radius= 0.5 em]; 
\end{tikzpicture}
};
\node (b) at (a.east) [xshift=6.5cm]
{
\begin{tikzpicture}[scale=0.4]
\begin{scope}
    \fill[color=red!50,opacity=0.3,thick,domain=0:2,samples=101]
    (0, 4) -- plot (\x,{-2*(\x^(3/2))}) --(7.9,{-2*(2^(3/2))}) --(7.9,4);
  \end{scope}
\draw[->] (-0.2,0) -- (8,0) node[right] {$c_1$};
\draw[->] (0,-5.75) -- (0,4) node[above] {$c_0$};
\foreach \t in {4,5,6,7,8} %n = 16 
{
	\pgfmathsetmacro\n{ 16 }
	\pgfmathsetmacro\R{(186-  186*\t/\n)/256 }
	\pgfmathsetmacro\G{(12 + 88*\t/\n)/256 }	
	\pgfmathsetmacro\B{(47 + 117*\t/\n)/256}
	\clip (-0.2,{-2*(2^(3/2))}) rectangle (8,4);
  	\draw[color={rgb, 1:red, \R; green, \G; blue, \B},domain=0:2,samples=100]    plot ( {\x} ,{ ((\n/\t)^3 - (3*\n/\t) )*(\x^(3/2)) })             ;
}
\foreach \t in {9} %n = 16 
{
	\pgfmathsetmacro\n{ 16 }
	\pgfmathsetmacro\R{(186-  186*\t/\n)/256 }
	\pgfmathsetmacro\G{(12 + 88*\t/\n)/256 }	
	\pgfmathsetmacro\B{(47 + 117*\t/\n)/256}
	\clip (-0.2,{-2*(2^(3/2))}) rectangle (8,4);
  	\draw[color={rgb, 1:red, \R; green, \G; blue, \B},domain=0:6,samples=100]    plot ( {\x} ,{    ((\n/\t)^3 - (3*\n/\t) )*(\x^(3/2))  })             ;
}
\foreach \t in {10} %n = 16 
{
	\pgfmathsetmacro\n{ 16 }
	\pgfmathsetmacro\R{(186-  186*\t/\n)/256 }
	\pgfmathsetmacro\G{(12 + 88*\t/\n)/256 }	
	\pgfmathsetmacro\B{(47 + 117*\t/\n)/256}
	\clip (-0.2,{-2*(2^(3/2))}) rectangle (8,4);
  	\draw[color={rgb, 1:red, \R; green, \G; blue, \B},domain=0:4,samples=100]    plot ( {\x} ,{    ((\n/\t)^3 - (3*\n/\t) )*(\x^(3/2))  })             ;
}
\foreach \t in {11,12} %n = 16 
{
	\pgfmathsetmacro\n{ 16 }
	\pgfmathsetmacro\R{(186-  186*\t/\n)/256 }
	\pgfmathsetmacro\G{(12 + 88*\t/\n)/256 }	
	\pgfmathsetmacro\B{(47 + 117*\t/\n)/256}
	\clip (-0.2,{-2*(2^(3/2))}) rectangle (8,4);
  	\draw[color={rgb, 1:red, \R; green, \G; blue, \B},domain=0:3,samples=100]    plot ( {\x} ,{ ((\n/\t)^3 - (3*\n/\t) )*(\x^(3/2))  })             ;
}
\draw[color=UCIB,dashed,thick,domain=0.1:2]    plot (\x,{-2*(\x^(3/2))})             ;
\draw[color=OSUR,thick,domain=0.1:3.9]    plot ({0},\x)             ;
\draw [color=red,fill=white] (0, 0) circle[radius= 0.5 em]; 
\end{tikzpicture}
};
%\draw [<->] (a)--(b);
\draw [->] (a.10)        to [bend left=30] (b.170) node[align=left] at (4.75, 1.55) {$c_1 = x_1^2$; \\ $c_0 = x_0^3 - 3x_0 x_1^2$.};
\draw [->] (b.190)        to [bend left=30] (a.350) node[align=left] at (4.7, -1.7) {$x_1 = \sqrt{c_1}$; \\ $x_0 = $ largest real root of \\ \quad \quad \    $x^3 -3c_1x - c_0$.};
\end{tikzpicture}
\caption{Homeomorphism between $\tilde{\mathcal{X}}_3$ and $\tilde{\mathscr{C}}_3$.}
\label{fig:2.1}
\end{figure}
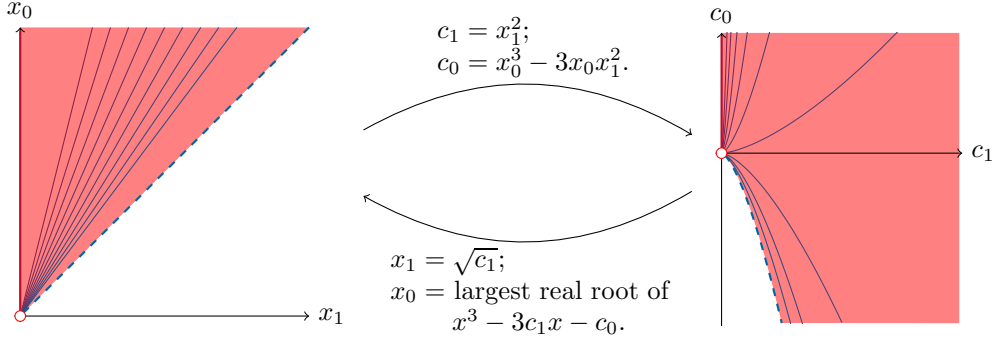

\hypertarget{L:2.3}{\begin{flemma}
Let $c \in \overline{\tilde{\mathscr{C}}_n}$ and $\varphi(c) = (x_{n-2}, \cdots, x_1, x_0)$, then for $l \in \{0, \cdots, n-2\}$, 
\begin{align*}
\sum_{k=l}^{n-2} c_k  \mybinom[0.8]{n-l-1}{k - l } x^{k-l}
\end{align*}is either a zero function when $x_{l} = x_{l+1}  = \cdots   = 0$ or a positive constant function $c_{l}$ when $x_{l} > x_{l+1} = \cdots    = 0$, or strictly increasing function on $[x_{l+1}, \infty)$ with $\sum_{k=l}^{n-2} c_k  \binom{n-l-1}{k - l} x_{l+1}^{k-l} \geq 0$ and equals $0$ if and only if $x_l = x_{l+1}$. In particular, if $c \in {\tilde{\mathscr{C}}_n}$, then $\sum_{k=0}^{n-2} c_k  \binom{n-1}{k } x_1^{k} > 0$.
\end{flemma}}

\begin{proof}
We prove this statement by mathematical induction on the degree $n$. When $n = 2$, it is straightforward. Suppose it is true when $n = m-1$. When $n = m$, by mathematical induction, 
\begin{align*}
\sum_{k=l}^{m-2} c_k  \mybinom[0.8]{m-l-1}{k - l} x^{k-l}
\end{align*}is either a zero function when $x_{l} = \cdots  = x_{m-2} = 0$ or a positive constant function $c_{l}$ when $x_{l} > x_{l+1} = \cdots  = x_{m-2} = 0$, or strictly increasing function on $[x_l, \infty)$ for $l \in \{1, \cdots, m-2\}$. For the case $l =0$, by plugging in $x_1$ to $x^m - \sum_{k=0}^{m-2}  c_k \binom{m}{k} x^k$ and the fact that $\varphi(c) \in \overline{\tilde{\mathcal{X}}_m}$, we obtain 
\begin{align*}
0 &\geq x_1^m  - \sum_{k=0}^{m-2} c_k \mybinom[0.8]{m}{k} x_1^k 
%= x_1  \Bigl(   x_1^{m-1}  - \sum_{k=1}^{m-2}  c_k \mybinom[0.8]{m-1}{k-1} x_1^{k-1}    \Bigr) - \sum_{k=0}^{m-2}  c_k \mybinom[0.8]{m-1}{k} x_1^{k} 
= - \sum_{k=0}^{m-2} c_k \mybinom[0.8]{m-1}{k} x_1^{k}.
\end{align*}Besides, for $y_1 \geq y_2 \geq x_1 \geq x_2$, by mathematical induction, we get
\begin{align*}
\label{eq:2.9}
&\kern-1em \sum_{k=0}^{m-2} c_k \mybinom[0.8]{m-1}{k} y_1^{k} - \sum_{k=0}^{m-2} c_k \mybinom[0.8]{m-1}{k} y_2^{k} \tag{2.9} \\
&= \int_{y_2}^{y_1} \frac{d}{dy} \Bigl(    \sum_{k=0}^{m-2} c_k \mybinom[0.8]{m-1}{k} y^{k}    \Bigr) dy = (m-1) \int_{y_2}^{y_1}    \sum_{k=1}^{m-2} c_k \mybinom[0.8]{m-2}{k-1} y^{k-1}  dy \geq 0.
\end{align*}If $x_0 = x_1 =  \cdots = x_{m-2} = 0$, then $c_0 = \cdots = c_{m-2} = 0$. If $x_0 > x_1 =  \cdots = x_{m-2} = 0$, then $c_0 > c_1 = \cdots = c_{m-2} = 0$ and $\sum_{k=0}^{m-2} c_k  \binom{m-1}{k } x^{k} = c_0$. Otherwise, by (\ref{eq:2.9}),  
\begin{align*}
\sum_{k=0}^{m-2} c_k  \mybinom[0.8]{m-1}{k } x^{k} = \sum_{k=0}^{m-2} c_k  \mybinom[0.8]{m-1}{k } x_1^{k} + (m-1) \int_{x_1}^{x}    \sum_{k=1}^{m-2} c_k \mybinom[0.8]{m-2}{k-1} y^{k-1}  dy
\end{align*}is strictly increasing on $[x_1, \infty)$ with $\sum_{k=0}^{m-2} c_k  \binom{m-1}{k } x_1^{k} \geq 0$. In particular, if $c \in {\tilde{\mathscr{C}}_n}$, then we have $\sum_{k=0}^{n-2} c_k  \binom{n-1}{k } x_1^{k} > 0$. This finishes the proof.
\end{proof}

\hypertarget{L:2.4}{\begin{flemma}
Let $c \in \overline{\tilde{\mathscr{C}}_n}$ and $\varphi(c) = (x_{n-2}, \cdots, x_1, x_0)$, then for $l \in \{0, \cdots, n-2\}$, 
\begin{align*}
\sum_{k=l}^{n-2} c_k  \mybinom[0.8]{n-l}{k - l} x^{k-l}
\end{align*}is either a zero function when $x_{l} = x_{l+1} =  \cdots    = 0$ or a positive constant function $c_{l}$ when $x_{l} > x_{l+1} = \cdots    = 0$, or strictly increasing function on $[x_{l+1}, \infty)$ with $\sum_{k=l}^{n-2} c_k  \binom{n-l}{k - l } x_{l+1}^{k-l} \geq x_{l+1}^{n-l}$ and equality happens if and only if $x_{l} = x_{l+1}$. In particular, if $c \in {\tilde{\mathscr{C}}_n}$, then $\sum_{k=0}^{n-2} c_k  \binom{n}{k } x_1^{k} > 0$.
\end{flemma}}

\begin{proof}
We prove this statement by mathematical induction on the degree $n$. When $n = 2$, it is straightforward. Suppose it is true when $n = m-1$. When $n = m$, by mathematical induction, 
\begin{align*}
\sum_{k=l}^{m-2} c_k  \mybinom[0.8]{m-l}{k - l} x^{k-l}
\end{align*}is either a zero function when $x_{l} = \cdots  = x_{m-2} = 0$ or a positive constant function $c_{l}$ when $x_{l} > x_{l+1} = \cdots  = x_{m-2} = 0$, or strictly increasing function on $[x_l, \infty)$ for $l \in \{1, \cdots, m-2\}$. For the case $l =0$, by plugging in $x_1$ to $x^m - \sum_{k=0}^{m-2}  c_k \binom{m}{k} x^k$ and the fact that $\varphi(c) \in \overline{\tilde{\mathcal{X}}_m}$, we obtain
\begin{align*}
0 \geq x_1^m  - \sum_{k=0}^{m-2} c_k \mybinom[0.8]{m}{k} x_1^k,
\end{align*}which implies that $\sum_{k=0}^{m-2} c_k \mybinom[0.8]{m}{k} x_1^k \geq x_1^m \geq 0$.
Besides, for $y_1 \geq y_2 \geq x_1 \geq x_2$, by mathematical induction, we get
\begin{align*}
\label{eq:2.10}
&\kern-1em \sum_{k=0}^{m-2} c_k \mybinom[0.8]{m}{k} y_1^{k} - \sum_{k=0}^{m-2} c_k \mybinom[0.8]{m}{k} y_2^{k} \tag{2.10} \\
&= \int_{y_2}^{y_1} \frac{d}{dy} \Bigl(    \sum_{k=0}^{m-2} c_k \mybinom[0.8]{m}{k} y^{k}    \Bigr) dy = m \int_{y_2}^{y_1}    \sum_{k=1}^{m-2} c_k \mybinom[0.8]{m-1}{k-1} y^{k-1}  dy \geq 0.
\end{align*}If $x_0 = x_1 =  \cdots = x_{m-2} = 0$, then $c_0 = \cdots = c_{m-2} = 0$. If $x_0 > x_1 =  \cdots = x_{m-2} = 0$, then $c_0 > c_1 = \cdots = c_{m-2} = 0$ and $\sum_{k=0}^{m-2} c_k  \binom{m}{k } x^{k} = c_0$. Otherwise, by (\ref{eq:2.10}),  
\begin{align*}
\sum_{k=0}^{m-2} c_k  \mybinom[0.8]{m}{k } x^{k} = \sum_{k=0}^{m-2} c_k  \mybinom[0.8]{m}{k } x_1^{k} + m \int_{x_1}^{x}    \sum_{k=1}^{m-2} c_k \mybinom[0.8]{m-1}{k-1} y^{k-1}  dy.
\end{align*}is strictly increasing on $[x_1, \infty)$ with $\sum_{k=0}^{m-2} c_k  \binom{m}{k } x_1^{k} \geq 0$. In particular, if $c \in {\tilde{\mathscr{C}}_n}$, then we have $\sum_{k=0}^{n-2} c_k  \binom{n}{k } x_1^{k} > 0$. This finishes the proof.
\end{proof}

\hypertarget{L:2.5}{\begin{flemma}
If $c \in \overline{\tilde{\mathscr{C}}_n}$, then for $t > 1$, $tc = 0 \in \overline{\tilde{\mathscr{C}}_n}$ when $c = 0$ and $tc  \in  {\tilde{\mathscr{C}}_n}$ otherwise. Moreover, for convenience, if we denote $(x_{n-2}(t), \cdots, x_1(t), x_0(t))$ the $(n-1)$-tuple $\varphi(tc) = (x_{n-2}(tc), \cdots, x_1(tc), x_0(tc))$, then for $t \geq 1$ and $l \in \{ 0, \cdots, n-2\}$, we have $\frac{d}{dt} x_l (t) = 0$ when $x_l = \cdots = x_{n-2} = 0$ and $\frac{d}{dt} x_k (t) > 0$ otherwise. In particular, if $c \in  {\tilde{\mathscr{C}}_n}$, then $\frac{d}{dt} x_0(t) > 0$.
\end{flemma}}

\begin{proof}
We use mathematical induction on the degree $n$ to prove it. When $n = 2$, it is straightforward. Suppose this statement is true when $n = m-1$. When $n = m$, by mathematical induction, for $l \in \{1, \cdots, m-2\}$, we have $\frac{d}{dt} x_l (t) = 0$ when $x_l = \cdots = x_{m-2} = 0$ and $\frac{d}{dt} x_l (t) > 0$ otherwise.
First, for the case $l = 0$, by plugging in $x_1(t)$ to $x^m - \sum_{k=0}^{m-2} t c_k \binom{m}{k} x^k$, we get
\begin{align*}
\label{eq:2.11}
x_1^m(t) - t \sum_{k=0}^{m-2}  c_k \mybinom[0.8]{m}{k} x_1^k(t) 
%&= x_1(t) \Bigl(   x_1^{n-1}(t)  - \sum_{k=1}^{n-2}  t c_k \mybinom[0.8]{n-1}{k-1} x_1^{k-1}(t)    \Bigr) - \sum_{k=0}^{n-2} t c_k \mybinom[0.8]{n-1}{k} x_1^{k}(t) \\
&= - t \sum_{k=0}^{m-2} c_k \mybinom[0.8]{m-1}{k} x_1^{k}(t).  \tag{2.11}
\end{align*}
Second, by Lemma~\hyperlink{L:2.3}{2.3} and by mathematical induction, we obtain
\begin{align*}
\label{eq:2.12}
&\kern-1em \sum_{k=0}^{m-2} c_k \mybinom[0.8]{m-1}{k} x_1^{k} (t) - \sum_{k=0}^{m-2} c_k \mybinom[0.8]{m-1}{k} x_1^{k}  \tag{2.12} \\
&= \int_{x_1}^{x_1(t)} \frac{d}{dx} \Bigl(    \sum_{k=0}^{m-2} c_k \mybinom[0.8]{m-1}{k} x^{k}    \Bigr) dx = (m-1) \int_{x_1}^{x_1(t)}    \sum_{k=1}^{m-2} c_k \mybinom[0.8]{m-2}{k-1} x^{k-1}  dx \geq 0.
\end{align*}By (\ref{eq:2.12}), we get
$\sum_{k=0}^{m-2} c_k \binom{m-1}{k} x_1^{k}(t) \geq  \sum_{k=0}^{m-2} c_k \binom{m-1}{k} x_1^{k} \geq 0$. Hence, by (\ref{eq:2.11}), $x_0(t)$ exists and $x_0(t) \geq x_1(t)$. To show that $x_0(t) \geq  x_0$ for $t \geq 1$, we have
\begin{align*}
x_0^m(t) - \sum_{k=0}^{m-2} c_k \mybinom[0.8]{m}{k} x_0^k(t) = (1 - 1/t) x_0^m(t) \geq    0.
\end{align*}This implies that $x_0(t) \geq x_0$ for $t \geq 1$. 

Last, for $\frac{d}{dt} x_0(t)$, if $c = 0$, then $x_0(t) = x_0 = 0$, which implies that $\frac{d}{dt} x_0(t) = 0$. Otherwise, there are two cases to consider: $x_0 > x_1 = \cdots = x_{m-2} = 0$ and the rest. For the case $x_0 > x_1 = \cdots = x_{m-2} = 0$, by mathematical induction we have $x_l(t) = x_l = 0$ for $l \in \{1, \cdots, m-2\}$ and $x_0(t) = \sqrt[m]{t} \cdot x_0$. For the rest, we have $\frac{d}{dt} x_1(t) > 0$ and by Lemma~\hyperlink{L:2.3}{2.3}, we get
\begin{align*}
\sum_{k=0}^{m-2} c_k \mybinom[0.8]{m-1}{k} x_1^{k}(t) >  \sum_{k=0}^{m-2} c_k \mybinom[0.8]{m-1}{k} x_1^{k} \geq 0
\end{align*}for $t >1$. Thus, by equality (\ref{eq:2.11}), we get $x_0(t) > x_1(t)$. In addition, by Lemma~\hyperlink{L:2.4}{2.4}, we get
\begin{align*}
\sum_{k=0}^{m-2} c_k \mybinom[0.8]{m}{k} x_0^{k}(t)  >\sum_{k=0}^{m-2} c_k \mybinom[0.8]{m}{k} x_1^{k}(t) >  \sum_{k=0}^{m-2} c_k \mybinom[0.8]{m}{k} x_1^{k} \geq 0
\end{align*}for $t >1$.
Hence, by implicit function theorem, we have
\begin{align*}
\label{eq:2.13}
m \frac{d}{dt} x_0(t) \cdot \Bigl[  x_0^{m-1}(t) - t \sum_{k=1}^{m-2}  c_k \mybinom[0.8]{m-1}{k-1} x_0^{k-1}(t)     \Bigr] = \sum_{k=0}^{m-2} c_k \mybinom[0.8]{m}{k} x_0^k(t), \tag{2.13}
\end{align*}
which implies that $\frac{d}{dt} x_0(t) >0$. This finishes the proof.
\end{proof}

\hypertarget{P:2.2}{\begin{fprop}
If $c \in \overline{\tilde{\mathscr{C}}_n}$ with $c \neq 0$ and we denote $(x_{n-2}(t), \cdots, x_1(t), x_0(t))$ the $(n-1)$-tuple $\varphi(tc) = (x_{n-2}(tc), \cdots, x_1(tc), x_0(tc))$, then $x_0(t) \geq \sqrt[n]{t} x_0$. In particular, $\lim_{t \rightarrow \infty} x_0(t) = \infty$.
\end{fprop}}

\begin{proof}
If $c \neq 0$, then by Lemma~\hyperlink{L:2.5}{2.5}, $t c \in \tilde{\mathscr{C}}_n$ for $t > 1$. In addition, by (\ref{eq:2.13}), we have
\begin{align*}
n \frac{d}{dt} x_0(t)  = \frac{ \sum_{k=0}^{n-2} c_k \binom{n}{k} x_0^k(t)}{x_0^{n-1}(t) - t \sum_{k=1}^{n-2}  c_k \binom{n-1}{k-1} x_0^{k-1}(t) } =  \frac{ x_0 \sum_{k=0}^{n-2} c_k \binom{n}{k} x_0^k(t)}{ t \sum_{k=0}^{n-2}  c_k \binom{n-1}{k} x_0^{k}(t) } > 0.
\end{align*}
On the other hand, by Lemma~\hyperlink{L:2.4}{2.4}, we have
\begin{align*}
\sum_{k=0}^{n-2} c_k \mybinom[0.8]{n}{k} x_0^k(t) - \sum_{k=0}^{n-2} c_k \mybinom[0.8]{n-1}{k} x_0^k(t) = \sum_{k=1}^{n-2} c_k \mybinom[0.8]{n-1}{k-1} x_0^k(t) >  \sum_{k=1}^{n-2} c_k \mybinom[0.8]{n-1}{k-1} x_1^k(t) \geq 0.
\end{align*}
This implies that for $t >1$, we have
\begin{align*}
n \frac{d}{dt} x_0(t)  =   \frac{ x_0 \sum_{k=0}^{n-2} c_k \binom{n}{k} x_0^k(t)}{ t \sum_{k=0}^{n-2}  c_k \binom{n-1}{k} x_0^{k}(t) }> \frac{x_0}{t}.
\end{align*}Hence, by standard calculus argument we have $x_0(t) \geq \sqrt[n]{t} \cdot x_0$. This finishes the proof.
\end{proof}

From the above results, we define the following function which is interesting in its own right. 
\hypertarget{D:2.8}{\begin{fdefi}
We define a function $\mathcal{T}$ from $\tilde{\mathscr{C}}_n$ to the extended real line $\overline{\mathbb{R}} = \mathbb{R} \cup \{ \infty, -\infty\}$ by
\begin{align*}
\mathcal{T}(c) &\coloneqq \sup \bigl\{  t \geq 1  : x_0 > x_1(t)    \bigr\}, \text{ for } c \in \tilde{\mathscr{C}}_n.
\end{align*}Here, $x_0$ is the largest real root of $x^n - \sum_{k=0}^{n-2} c_k \binom{n}{k} x^k$ and $x_1(t)$ is the largest real root of $x^{n-1} - t \sum_{k=1}^{n-2} c_k \binom{n-1}{k-1} x^{k-1}$.  
\end{fdefi}}

\hypertarget{P:2.3}{\begin{fprop}
If $c \in \CY_n$, then $\mathcal{T}(c) = \infty$. If $c \in   \tilde{\mathscr{C}}_n \backslash \CY_n$, then $\infty > \mathcal{T}(c) > 1$. 
\end{fprop}
\begin{proof}
If $c \in \CY_n \subset \tilde{\mathscr{C}}_n$, then we have $c_1 = \cdots = c_{n-2} = 0$ and $c_0 > 0$. Hence, $x_0 = \sqrt[n]{c_0}$ and $x_1(t) = 0$ for $t \geq 1$. Thus, $\mathcal{T}(c) = \infty$. If $c \in   \tilde{\mathscr{C}}_n \backslash \CY_n$, since $x_0 > x_1 = x_1(1)$, we get $\mathcal{T}(c) > 1$. In addition, since $c \notin \CY_n$, by Proposition~\hyperlink{P:2.2}{2.2}, we have $\lim_{t \rightarrow \infty} x_1(t) = \infty$. Hence, by the intermediate value theorem, we have $\mathcal{T}(c) < \infty$. This finishes the proof.
\end{proof}}

\hypertarget{P:2.4}{\begin{fprop}
Let $c \in \tilde{\mathscr{C}}_n \backslash \CY_n$, $x_0$ be the largest real root of $x^n - \sum_{k=0}^{n-2} c_k \binom{n}{k} x^k$, and $x_1(t)$ be the largest real root of $x^{n-1} - \sum_{k=1}^{n-2} tc_k \binom{n-1}{k-1} x^{k-1}$ for $t \geq 1$. Then, we have $x_0 = x_1 \bigl( \mathcal{T}(c) \bigr)$.
\end{fprop}}

\begin{proof}
By Proposition~\hyperlink{P:2.3}{2.3}, we have $\mathcal{T}(c) = \sup \bigl\{  t \geq 1  : x_0 > x_1(t)    \bigr\} < \infty$. Also, by Lemma~\hyperlink{L:2.5}{2.5}, $x_1(t)$ is strictly increasing for $t \geq 1$. Hence, by the continuity of $x_1(t)$, we have $x_0 = x_1 \bigl( \mathcal{T}(c) \bigr)$. 
\end{proof}

\hypertarget{P:2.5}{\begin{fprop}
The function $\mathcal{T}  \colon \tilde{\mathscr{C}}_n \rightarrow \overline{\mathbb{R}}$ is a continuous function. 
\end{fprop}}

\begin{proof}
There are two cases to consider: $c \in \CY_n$ or $c \in  \tilde{\mathscr{C}}_n \backslash \CY_n$. For the case $c \in  \tilde{\mathscr{C}}_n \backslash \CY_n$, by Proposition~\hyperlink{P:2.4}{2.4}, we have
$x_1\bigl ( \mathcal{T}(c)  \bigr) = x_0$.
This implies that
\begin{align*}
0 = x_1^{n-1} \bigl(\mathcal{T}(c) \bigr) - \sum_{k=1}^{n-2} \mathcal{T}(c) c_k \mybinom[0.8]{n-1}{k-1} x_1^{k-1}\bigl(\mathcal{T}(c) \bigr) =  x_0^{n-1} - \sum_{k=1}^{n-2} \mathcal{T}(c) c_k \mybinom[0.8]{n-1}{k-1} x_0^{k-1}. 
\end{align*}Hence, for $c \in  \tilde{\mathscr{C}}_n \backslash \CY_n$, we have
\begin{align*}
\mathcal{T}(c) = \frac{x_0^{n-1}}{\sum_{k=1}^{n-2}   c_k \binom{n-1}{k-1} x_0^{k-1} }.
\end{align*}By Lemma~\hyperlink{L:2.2}{2.2}, $\tilde{\mathcal{X}}_n$ is homeomorphic to $\tilde{\mathscr{C}}_n$, so $\sum_{k=1}^{n-2}   c_k \binom{n-1}{k-1} x_0^{k-1}$ is a continuous function on $\tilde{\mathscr{C}}_n$. Also, by Lemma~\hyperlink{L:2.4}{2.4}, we have $\sum_{k=1}^{n-2}   c_k \binom{n-1}{k-1} x_0^{k-1} >0$. Combine these, $\mathcal{T}$ is continuous at $c$. \smallskip

Now, if $c \in \CY_n$, then for any $\delta > 0$, we can always find $\tilde{c} \in \CY_n$ with $0 < \| c - \tilde{c} \| < \delta$, we immediately get $\limsup_{x \rightarrow c} \mathcal{T}(x) = \infty$. On the other hand, for any $\tilde{c} \in \tilde{\mathscr{C}}_n \backslash \CY_n$ with $0 < \| c - \tilde{c} \| < \delta$, we have the following rough estimate:  
\begin{align*}
0 \leq  {x}_1(\tilde{c})  \leq \sum_{k=1}^{n-2} \sqrt[n-k]{|\tilde{c}_k| \mybinom[0.8]{n-1}{k-1}}.
\end{align*}This implies that for $t \geq 1$, we get
\begin{align*}
{x}_1(t \tilde{c}) \leq \sum_{k=1}^{n-2} \sqrt[n-k]{|\tilde{c}_k| \mybinom[0.8]{n-1}{k-1}} \cdot t^{\frac{1}{n-k}} \leq \sum_{k=1}^{n-2} \sqrt[n-k]{|\tilde{c}_k| \mybinom[0.8]{n-1}{k-1}} \cdot t^{\frac{1}{2}}.
\end{align*}Hence, we get the following rough estimate
\begin{align*}
\mathcal{T}(\tilde{c}) \geq \Biggl( \frac{{x}_0(\tilde{c})}{\sum_{k=1}^{n-2} \sqrt[n-k]{|\tilde{c}_k| \binom{n-1}{k-1}}} \Biggr)^2,
\end{align*}which implies that $\liminf_{x \rightarrow c} \mathcal{T}(x) = \infty$. Hence, $\lim_{x \rightarrow c} \mathcal{T}(c) = \infty = \mathcal{T}(c)$ giving that $\mathcal{T}$ is a continuous function from $\tilde{\mathscr{C}}_n$ to $\overline{\mathbb{R}}$. This finishes the proof.  
\end{proof}

\subsection{Properties of Right-Noetherian Polynomials}
\label{sec:2.3}
In Section~\ref{sec:2.1}, we have introduced the background and notations of $\Upsilon$-stable general inverse $\sigma_k$ multilinear polynomials. In Section~\ref{sec:2.2}, we have introduced right-Noetherian polynomials and shown the correspondence between right-Noetherian polynomials and $\Upsilon$-stable general inverse $\sigma_k$ multilinear polynomials. Now, we are able to use this correspondence to find more properties of strictly $\Upsilon$-stable general inverse $\sigma_k$ equations.\smallskip

We state the following definition and a Positivstellensatz type result in \cite{lin2023c}, which not only show us how to compare $\Upsilon$-stable general inverse $\sigma_k$ multilinear polynomials but also help us understand which conditions the $C$-subsolution condition will be maintained. The capability to compare $\Upsilon$-stable general inverse $\sigma_k$ multilinear polynomials is very important, this leads us to create various polynomials which might not be obvious at first glimpse.

\hypertarget{D:2.9}{\begin{fdefi}[\texorpdfstring{$\Upsilon$}{}-dominance, L.~\cite{lin2023c}]
Let $f(\lambda) \coloneqq   \lambda_1 \cdots \lambda_n - \sum_{k = 0}^{n-1} c_k \sigma_k(\lambda)$ and $g(\lambda) \coloneqq   \lambda_1 \cdots \lambda_n - \sum_{k = 0}^{n-1} d_k \sigma_k(\lambda)$ be two $\Upsilon$-stable general inverse $\sigma_k$ multilinear polynomials. For $l \in \{0, \cdots, n-1\}$, we write $x_l$ the largest real root of the $l$-th derivative $r_f^{(k)}$ of the diagonal restriction $r_f$ of $f$ and $y_l$ the largest real root of the $l$-th derivative $r_g^{(l)}$ of the diagonal restriction $r_g$ of $g$. If $y_l \geq x_l$ for all $l \in \{0, \cdots, n-1\}$, then we say $g \gtrdot f$. 
\end{fdefi}}

\hypertarget{T:2.3}{\begin{fthm}[L.~\cite{lin2023c}]
Let $f(\lambda) \coloneqq   \lambda_1 \cdots \lambda_n - \sum_{k = 0}^{n-1} c_k \sigma_k(\lambda)$ and $g(\lambda) \coloneqq   \lambda_1 \cdots \lambda_n - \sum_{k = 0}^{n-1} d_k \sigma_k(\lambda)$ be two $\Upsilon$-stable general inverse $\sigma_k$ multilinear polynomials. Then $g \gtrdot f$ if and only if $\Gamma^n_g \subset \Gamma^n_f$.
\end{fthm}}

Sometimes we will just focus on comparing these roots, so we have the following definition.

\hypertarget{D:2.10}{\begin{fdefi}
For $(x_{n-2}, \cdots, x_1, x_{0}), (y_{n-2}, \cdots, y_1, y_{0}) \in \tilde{\mathcal{X}}_n$, we say $ (y_{n-2}, \cdots, y_1, y_{0}) \gtrdot (x_{n-2}, \cdots, x_1, x_{0})$ if $y_l \geq x_l$ for all $l \in \{0, \cdots, n-2\}$.
\end{fdefi}}

Based on previous results, the existence of a $C$-subsolution to any classical strictly $\Upsilon$-stable general inverse $\sigma_k$ equation should provide us with a priori estimates for the solution to this strictly $\Upsilon$-stable general inverse $\sigma_k$ equation. In Section~\ref{sec:3}, we will confirm that for any map $c \colon M \rightarrow \tilde{\mathscr{C}}_n$ with $c_1, \cdots, c_{n-2}$ constants, we have a priori estimates provided the existence of a $C$-subsolution to $c \colon M \rightarrow \tilde{\mathscr{C}}_n$.
Therefore, if a $C$-subsolution to a strictly $\Upsilon$-stable general inverse $\sigma_k$ equation is given, to apply the method of continuity to obtain the solvability, we are interested in those strictly $\Upsilon$-stable general inverse $\sigma_k$ equations such that the original $C$-subsolution will still be a $C$-subsolution to these equations.
That is to say, by Definition~\hyperlink{D:2.5}{2.5}, we want to find the space collecting all strictly $\Upsilon$-stable general inverse $\sigma_k$ equations such that the original $C$-subsolution will still be in the $\Upsilon_1$-cone of them.
In the author's thesis \cite{lin2023thesis}, the author suggests the following space should be the space that we should be looking for.

\hypertarget{T:2.4}{\begin{fthm}
Given $c, d \in \tilde{\mathscr{C}}_{n}$, we write $x_l(c)$ and $x_l(d)$ the largest real root of the $l$-th derivative of $f(x) =   x^n - \sum_{k=0}^{n-2} c_k \binom{n}{k} x^k$ and $g(x) =   x^n - \sum_{k=0}^{n-2} d_k \binom{n}{k} x^k$, respectively. Then the $\Upsilon_1$-cone of $g$ is contained in the $\Upsilon_1$-cone of $f$ if and only if for all $l \in \{1, \cdots, n-2\}$, we have
$f^{(l)}(x_l(d)) \geq 0$.  
\end{fthm}}
\begin{proof}
To show that the $\Upsilon_1$-cone of $g$ is contained in the $\Upsilon_1$-cone of $f$. By Theorem~\hyperlink{T:2.3}{2.3}, we are checking whether $x_l(d) \geq x_l(c)$ for all $l \in \{1, \cdots, n-1\}$. When $l = n-1$, $f^{(n-1)}(x) =  {n!} \cdot x $, $g^{(n-1)}(x) = n! \cdot x$, and $x_{n-1}(c) = 0 = x_{n-1}(d)$. When $l = n-2$, $f^{(n-2)}(x) =  \frac{n!}{2!} ( x^2  - c_{n-2})  $, $g^{(n-2)}(x) = \frac{n!}{2!} ( x^2  - d_{n-2}) $, $x_{n-1}(c) = \sqrt{c_{n-2}}$, and $x_{n-1}(d) = \sqrt{d_{n-2}}$.
So, $x_{n-2}(d) \geq x_{n-2}(c)$ if and only if $f^{(n-2)}(x_{n-2}(d)) = \frac{n!}{2!}  (x_{n-2}^2(d) - c_{n-2} ) = \frac{n!}{2!}  (x_{n-2}^2(d) - x_{n-2}^2(c) ) \geq 0$. \smallskip

Now, we use mathematical induction on $l$ from $n-2$ to $1$. We have already shown the case $l = n-2$. Suppose the statement is true when $l = s \geq 2$. When $l  = s-1$, if $x_{s-1}(d) \geq x_{s-1}(c)$, then since $x_{s-1}(c)$ is the largest real root of $f^{(s-1)}$, we get $f^{(s-1)}(x_{s-1}(d)) \geq 0$. On the other hand, if $f^{(s-1)}(x_{s-1}(d)) \geq 0$, then by mathematical induction, we have $x_{s-1}(d) \geq x_{s}(d) \geq x_{s}(c)$. By Proposition~\hyperlink{P:2.1}{2.1}, we obtain $x_{s-1}(d) \geq x_{s-1}(c)$. This finishes the proof.
\end{proof}

As a quick consequence, we get the following Corollary.

\hypertarget{Cor:2.2}{\begin{fcor}
Let $d \in \tilde{\mathscr{C}}_{n}$. Then for any ${c} \in \tilde{\mathscr{C}}_{n}$ satisfies the following $n-2 \times n-2$ linear system:
\begin{align*}
\label{eq:2.14}
\begin{cases}
&  x_{n-2}^2(d)   - {c}_{n-2}  \geq 0; \\
&\kern4em \vdots \\
&  x_2^{n-2}(d) - \sum_{k = 2}^{n-2} {c}_k \binom{n-2}{k-2} x_2^{k-2}(d) \geq 0; \\
& x_1^{n-1}(d) - \sum_{k = 1}^{n-2} {c}_k \binom{n-1}{k-1} x_1^{k-1}(d) \geq  0,
\end{cases} \tag{2.14}
\end{align*}the $\Upsilon_1$-cone of $x^n - \sum_{k=0}^{n-2} d_k \binom{n}{k} x^k$ will be contained in the $\Upsilon_1$-cone of $x^n - \sum_{k=0}^{n-2} c_k \binom{n}{k} x^k$. In other words, for any ${c} \in \tilde{\mathscr{C}}_{n}$ satisfies (\ref{eq:2.14}), $({c}_{n-2}, \cdots, {c}_1)$ lies in one of the polyhedrons, containing $(-R, -R^2, \cdots, -R^{n-2})$ for $R >0$ sufficiently large, defined by $n-2$ hypersurfaces passing through $({d}_{n-2}, \cdots, {d}_1)$ with the following $n-2$ linearly independent vectors as normal vectors:
\begin{align*}
\begin{pmatrix}
1 \\
0 \\
0 \\
\vdots \\
0 \\
0
\end{pmatrix}, \quad 
\begin{pmatrix}
\binom{3}{1} x_{n-3}(d) \\
1 \\
0 \\
\vdots \\
0 \\
0
\end{pmatrix},  \cdots,\quad 
\begin{pmatrix}
\binom{n-2}{n-4} x_{2}^{n-4}(d) \\
\binom{n-2}{n-5} x_{2}^{n-5}(d)  \\
\binom{n-2}{n-6} x_{2}^{n-6}(d) \\
\vdots \\
1 \\
0
\end{pmatrix}, \quad
\begin{pmatrix}
\binom{n-1}{n-3} x_{1}^{n-3}(d) \\
\binom{n-1}{n-4} x_{1}^{n-4}(d)  \\
\binom{n-1}{n-5} x_{1}^{n-5}(d) \\
\vdots \\
\binom{n-1}{1} x_{1}(d) \\
1
\end{pmatrix}.
\end{align*}
\end{fcor}}

By Theorem~\hyperlink{T:2.4}{2.4} and Corollary~\hyperlink{Cor:2.2}{2.2}, we consider the following polyhedron determined by $n-2 \times n-2$ linear system (\ref{eq:2.14}) for a $\Upsilon$-stable general inverse $\sigma_k$ multilinear polynomial.

\hypertarget{D:2.11}{\begin{fdefi}
Let $d  \in \tilde{\mathscr{C}}_{n}$, then we define the following polyhedron in $\tilde{\mathscr{C}}_{n}$:
    \begin{align*}
        \tilde{\mathcal{P}}^d \coloneqq \Bigl \{  {c} \in   \tilde{\mathscr{C}}_{n} \colon     x_l^{n-l}(d)  - \sum_{k=l}^{n-2} {c}_k\mybinom[0.8]{n-l}{k-l} x_l^{k-l}(d)  \geq 0, \quad \forall l \in \{1, \cdots, n-2\}   \Bigr \}.  
    \end{align*} 
\end{fdefi}}

\begin{figure} 
\centering
\begin{tikzpicture}[scale=0.4]
\begin{scope}
    %\clip(0, 0) -- (0, 3) -- (4.2325, 3) --  (4.2325, 0);
    \fill[color=red!50,opacity=0.3,thick,domain=0:2,samples=101]
    (0, 4) -- plot (\x,{-2*(\x^(3/2))}) --(5.9,{-2*(2^(3/2))}) --(5.9,4);
  \end{scope}
\begin{scope}
    %\clip(0, 0) -- (0, 3) -- (4.2325, 3) --  (4.2325, 0);
    \fill[color=blue!50,opacity=0.3,thick,domain=0:2,samples=101]
    (0, 4) -- plot (\x,{-2*(\x^(3/2))}) --(4,{-2*(2^(3/2))}) --(4,4);
  \end{scope}
  \draw[color=UCIB,thick,domain={-2*(2^(3/2))}:4]    plot ({4},\x)             ;
 \filldraw [color=UCIB] (4, -4) circle[radius= 0.15 em] node[right] {$(d_1, d_0)$}; 
  \draw[->] (-0.2,0) -- (6,0) node[right] {$c_1$};
  \draw[->] (0,-5.75) -- (0,4) node[above] {$c_0$};
  \draw[color=red,thick,domain=0.1:3.8]    plot ({0},\x)             ;
  \draw[color=red,dashed,thick,domain=0.1:2]    plot (\x,{-2*(\x^(3/2))})             ;
 \draw [color=red,fill=white] (0, 0) circle[radius= 0.4 em]; 
% \draw[color=UCIB,dashed,domain=0:1]    plot ( {1 + 2*\x} ,{-5.1*\x})            ;
% \draw[color=UCIB,dashed,domain=0:1]    plot ( {3} ,{-5.1*\x})            ;
\end{tikzpicture}
\caption{Polyhedron $\tilde{\mathcal{P}}^d$ of $\lambda_1 \lambda_2 \lambda_3 - d_1(\lambda_1 + \lambda_2 + \lambda_3) - d_0 = 0$.}
\label{fig:2.2}
\end{figure}
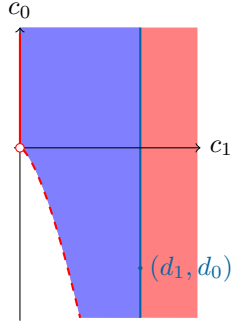

In Figure~\ref{fig:2.2}, we plot the polyhedron $\tilde{\mathcal{P}}^d$ for a $d = (d_1, d_0) \in \tilde{\mathscr{C}}_3$, which is the purple shaded region. Notice that the pink shaded region is the space $\tilde{\mathscr{C}}_3$.
Similar to before, we can extend this definition to any strictly $\Upsilon$-stable general inverse $\sigma_k$ equation $c \colon M \rightarrow \tilde{\mathscr{C}}_n$, we have the following. 
\hypertarget{D:2.12}{\begin{fdefi}
Let $d  \colon M \rightarrow \tilde{\mathscr{C}}_{n}$, then we say $c  \colon M \rightarrow \tilde{\mathscr{C}}_{n}$ satisfies $c \in \tilde{\mathcal{P}}^d$ if for any $z \in M$,  
$c(z) \in \tilde{\mathcal{P}}^{d(z)}$. 
\end{fdefi}}

\hypertarget{R:2.4}{\begin{frmk}
Given a strictly $\Upsilon$-stable general inverse $\sigma_k$ equation $d \colon M \rightarrow \tilde{\mathscr{C}}_{n}$, if there exists a $C$-subsolution $X_{\underline{u}}$ to $d$, then $X_{\underline{u}}$ is also a $C$-subsolution to $c \colon M \rightarrow \tilde{\mathscr{C}}_{n}$ for any $c \in \tilde{\mathcal{P}}^d$.
\end{frmk}}

We state the following Lemma in \cite{lin2023c}, which involves both $\Upsilon$-stable general inverse $\sigma_k$ equation and right-Noetherian polynomial, and the Lemma is very useful in many estimates.

\hypertarget{L:2.6}{\begin{flemma}[L.~\cite{lin2023c}]}
Let $c \in \overline{\tilde{\mathscr{C}}_n}$. If the $\Upsilon_l$-cone of $\lambda_1 \cdots \lambda_n - \sum_{k=0}^{n-2} c_k \sigma_k(\lambda)$ not equals $\Gamma_n$ for some $l \in \{0, \cdots, n-2\}$, then either
\begin{align*}
\partial \Upsilon_l  \cap  \partial \Upsilon_{l+1} = \emptyset \, \text{ or } \, \{(x_{l}, \cdots, x_{l})\}. 
\end{align*}If $\partial \Upsilon_l  \cap  \partial \Upsilon_{l+1} \neq  \emptyset$, then $x_{l} = x_{l+1}$.
In particular, if $\Gamma^n_f$ is strictly $\Upsilon$-stable, then
\begin{align*}
\partial \Upsilon_0 \cap  \partial \Upsilon_{1} = \partial \Gamma^n_f \cap \partial \Upsilon_{1} = \emptyset.
\end{align*}
\end{flemma}

From now on, for notational conventions, for any degree $n$ $\Upsilon$-stable general inverse $\sigma_k$ multilinear polynomial $\lambda_1 \cdots \lambda_n - \sum_{k=0}^{n-2} c_k \sigma_k(\lambda)$ and for any $l \in \{0, \cdots, n-2\}$, we will write 
\[ \bigcup_{s \in S_n} \{ \lambda_{s(l+1)} \cdots \lambda_{s(n)} - \sum_{k=l}^{n-2} c_k \sigma_k(\lambda; s(1), \cdots, s(l)) = 0\} \] 
instead of $\partial \Upsilon_l$ for convenience. But keep in mind that they are not the same. For example, $\lambda_1 \cdots \lambda_n$ is a $\Upsilon$-stable general inverse $\sigma_k$ multilinear polynomial. In this case, $\partial \Upsilon_0$ is the boundary of the positive orthant, but $\lambda_1 \cdots \lambda_n = 0$ will be the union of hyperplanes $\{\lambda_i = 0\}$. We hope that this will not cause any confusion.\smallskip

We can now state some results, we have the following.

\hypertarget{L:2.7}{\begin{flemma}
Let $c \in \overline{\tilde{\mathscr{C}}_n}$. Then for $(\lambda_1, \cdots, \lambda_n) \in \{ \lambda_1 \cdots \lambda_n - \sum_{k=0}^{n-2} c_k \sigma_k(\lambda) = 0\}$, for $l \in \{0, \cdots, n-2\}$, and for any $s \in S_n$, we have
\begin{align*}
\sum_{k = l}^{n-2} c_k \sigma_{k-l}(\lambda_{; s(1), \cdots,   s(l+1)}) &\geq \sum_{k = l}^{n-2} c_k \mybinom[0.8]{n-l-1}{k-l} x^{k-l}_{l+1}; \\
\sum_{k = l}^{n-2} c_k \sigma_{k-l}(\lambda_{; s(1), \cdots, s(l)}) &\geq  \sum_{k = l}^{n-2} c_k \mybinom[0.8]{n-l}{k-l}x_l^{k-l}  = x_l^{n-l}.
\end{align*}Moreover, for any $l \in \{0, \cdots, n-3\}$, if $(c_{n-2}, \cdots, c_{l+1}) \neq 0$, then $\sum_{k = l}^{n-2} c_k \sigma_{k-l}(\lambda_{; s(1), \cdots,   s(l+1)}) = \sum_{k = l}^{n-2} c_k \binom{n-l-1}{k-l} x^{k-l}_{l+1}$ if and only if $x_0 = \cdots = x_{l+1}$ and $\lambda_{s(l+2)} = \cdots = \lambda_{s(n)} = x_{l+1}$.
Similarly, for any $l \in \{0, \cdots, n-3\}$, if $(c_{n-2}, \cdots, c_{l+1}) \neq 0$, then $\sum_{k = l}^{n-2} c_k \sigma_{k-l}(\lambda_{; s(1), \cdots,   s(l)}) = \sum_{k = l}^{n-2} c_k \binom{n-l}{k-l} x^{k-l}_{l}$ if and only if $x_0 = \cdots = x_l$ and $\lambda_{s(l+1)} = \cdots = \lambda_{s(n)} = x_{l}$.
In addition, for any $l \in \{0, \cdots, n-2\}$, $\sum_{k = l}^{n-2} c_k \binom{n-l}{k-l}x_l^{k-l}  = x_l^{n-l} = 0$ if and only if $(c_{n-2}, \cdots, c_l) = 0$. Similarly, for any $l \in \{0, \cdots, n-2\}$, $\sum_{k = l}^{n-2} c_k \binom{n-l-1}{k-l} x^{k-l}_{l+1} = 0$ if and only if $x_l = x_{l+1}$. In particular, if $c \in  \tilde{\mathscr{C}}_n$, then $\sum_{k = 0}^{n-2} c_k \sigma_{k}(\lambda)  \geq   x_0^n    > 0$ and $\sum_{k = 0}^{n-2} c_k \sigma_{k}(\lambda_{; i}) \geq   \sum_{k = 0}^{n-2} c_k \binom{n-1}{k}x_1^{k}   > 0$.
\end{flemma}}
\begin{proof}
We prove both statements by mathematical induction on the index $l$ from $n-2$ to $0$ simultaneously. The proofs are alternated, to be more precise, suppose we prove the statement is true for the quantity $\sum_{k = l}^{n-2} c_k \sigma_{k-l}(\lambda_{; s(1), \cdots,   s(l+1)})$, then we are able to prove the statement is hence true for the quantity $\sum_{k = l}^{n-2} c_k \sigma_{k-l}(\lambda_{; s(1), \cdots, s(l)})$, and last we can prove the statement is also true for the quantity $\sum_{k = l-1}^{n-2} c_k \sigma_{k-l+1}(\lambda_{; s(1), \cdots, s(l)})$. By following this alternating way, we are able to prove both statements eventually.\smallskip

First, when $l = n-2$, then we have $\sum_{k = l}^{n-2} c_k \sigma_{k-l}(\lambda_{; s(1), \cdots,   s(l+1)}) = c_{n-2} \geq 0$ and equals $0$ if and only if $x_{n-2}  = x_{n-1} = 0$. Similarly, we have $\sum_{k = l}^{n-2} c_k \sigma_{k-l}(\lambda_{; s(1), \cdots, s(l)}) = c_{n-2} \geq 0$ and equals $0$ if and only if $c_{n-2} = 0$. When $l = n-3$, we have $\sum_{k = l}^{n-2} c_k \sigma_{k-l}(\lambda_{; s(1), \cdots,   s(l+1)}) = c_{n-3} + c_{n-2} ( \lambda_{s(n-1)} + \lambda_{s(n)} )$. If $c_{n-2} = 0$, then $\sum_{k = l}^{n-2} c_k \sigma_{k-l}(\lambda_{; s(1), \cdots,   s(l+1)}) = c_{n-3} \geq 0$ and equals $0$ if and only if $x_{n-3} = x_{n-2} = 0$. If $c_{n-2} > 0$, then $\sum_{k = l}^{n-2} c_k \sigma_{k-l}(\lambda_{; s(1), \cdots,   s(l+1)}) = c_{n-3} + c_{n-2} ( \lambda_{s(n-1)} + \lambda_{s(n)} ) \geq c_{n-3} + 2c_{n-2} \sqrt{\lambda_{s(n-1)} \lambda_{s(n)}} \geq c_{n-3} + 2c_{n-2}^{3/2} = c_{n-3} + 2c_{n-2}x_{n-2}$. By Lemma~\hyperlink{L:2.6}{2.6}, equality holds if and only if $x_0 = \cdots = x_{n-2}$ and $\lambda_{s(n-1)} = \lambda_{s(n)} = x_{n-2}$. Similarly, we have $\sum_{k = l}^{n-2} c_k \sigma_{k-l}(\lambda_{; s(1), \cdots,   s(l)}) = c_{n-3} + c_{n-2} ( \lambda_{s(n-2)} +  \lambda_{s(n-1)} + \lambda_{s(n)} )$. If $c_{n-2} = 0$, then $\sum_{k = l}^{n-2} c_k \sigma_{k-l}(\lambda_{; s(1), \cdots,   s(l)}) = c_{n-3} \geq 0$ and equals $0$ if and only if $(c_{n-2}, c_{n-3}) = 0$. If $c_{n-2}> 0$, then by Lemma~\hyperlink{L:2.6}{2.6}, we get $\sum_{k = l}^{n-2} c_k \sigma_{k-l}(\lambda_{; s(1), \cdots,   s(l)}) = c_{n-3} + c_{n-2} ( \lambda_{s(n-2)} +  \lambda_{s(n-1)} + \lambda_{s(n)} ) \geq c_{n-3} + 3c_{n-2} \sqrt[3]{\lambda_{(s(n-2))} \lambda_{(s(n-1))} \lambda_{(s(n))}    } \geq c_{n-3} + 3c_{n-2} x_{0} \geq c_{n-3} + 3c_{n-2} x_{n-3}$. Equality holds if and only if $x_0 = \cdots = x_{n-3}$ and $\lambda_{s(n-2)} = \lambda_{s(n-1)} = \lambda_{s(n)} = x_{n-3}$.
\smallskip
 
Second, suppose both statements are true when $l = m > 0$. When $l = m-1$, without loss of generality, we only need to consider the quantity $\sum_{k = m-1}^{n-2} c_k \sigma_{k-m+1}(\lambda_{; 1, \cdots, m})$. There are two cases to consider. For the first case, if the quantity $\sum_{k = m}^{n-2} c_k \sigma_{k-m}(\lambda_{; 1, \cdots, m}) = 0$, then $(c_{n-2}, \cdots, c_m)  = 0$. This implies that $\sum_{k = m-1}^{n-2} c_k \sigma_{k-m+1}(\lambda_{; 1, \cdots, m}) = c_{m-1} = x_{m-1}^{n-m+1} \geq 0$ and equals $0$ if and only if $c_{m-1} = 0$. In this case, if $c_{m-1} = 0$, then $x_{m-1} = x_{m} = \cdots = x_{n-1} = 0$. For the second case, if the quantity $\sum_{k = m}^{n-2} c_k \sigma_{k-m}(\lambda_{; 1, \cdots, m}) > 0$, then
by the fact that $\lambda_{m+1} \cdots \lambda_n \geq  \sum_{k=m}^{n-2} c_k \sigma_{k-m}(\lambda_{; 1, \cdots, m})$, there exists a $t \geq 1$ such that 
\begin{align*}
\lambda_{m+1} \cdots \lambda_n - t \sum_{k=m}^{n-2} c_k \sigma_{k-m}(\lambda_{; 1, \cdots, m}) = 0. 
\end{align*}
We view $(\lambda_{m+1}, \cdots, \lambda_n)$ as a point on the new $\Upsilon$-stable general inverse $\sigma_k$ multilinear polynomial $\lambda_{m+1} \cdots \lambda_n - t \sum_{k=m}^{n-2} c_k \sigma_{k-m}(\lambda_{; 1, \cdots, m}) = 0$ and $\sum_{k=m-1}^{n-2} c_k \sigma_{k-m+1}(\lambda_{; 1, \cdots, m})$ as a function on $\lambda_{m+1} \cdots \lambda_n - t \sum_{k=m}^{n-2} c_k \sigma_{k-m}(\lambda_{; 1, \cdots, m}) = 0$. For convenience, we say $\lambda_{m+1} \geq \cdots \geq \lambda_n$ and treat $\lambda_n$ as a function of variables $\{ \lambda_{m+1}, \cdots, \lambda_{n-1} \}$. For $j \in \{m+1, \cdots, n-1\}$, we have
\begingroup
\allowdisplaybreaks
\begin{align*}
\label{eq:2.15}
&\kern-2em \frac{\partial}{\partial \lambda_j}  \Bigl(  \sum_{k=m-1}^{n-2} c_k \sigma_{k-m+1}(\lambda_{; 1, \cdots, m}) \Bigr)  \tag{2.15} \\
&=  \sum_{k =m}^{n-2} c_k \sigma_{k-m}(\lambda_{; 1, \cdots, m, j}) + \frac{\partial \lambda_n}{\partial \lambda_j}  \sum_{k =m}^{n-2} c_k \sigma_{k-m}(\lambda_{; 1, \cdots, m, n}) \\
&=  \lambda_j  \Bigl (       \frac{\lambda_{m+1} \cdots \lambda_n}{t \lambda_j} -  \sum_{k=m+1}^{n-2} c_k \sigma_{k-m-1}(\lambda_{; 1, \cdots, m, j})       \Bigr)     \\
&\kern2em + \lambda_n \frac{\partial \lambda_n}{\partial \lambda_j}  \Bigl (       \frac{\lambda_{m+1} \cdots \lambda_{n-1}}{t} -  \sum_{k=m+1}^{n-2} c_k \sigma_{k-m-1}(\lambda_{; 1, \cdots, m, n})       \Bigr).              
\end{align*}
\endgroup
On the other hand, on $\lambda_{m+1} \cdots \lambda_n - t \sum_{k=m}^{n-2} c_k \sigma_{k-m}(\lambda_{; 1, \cdots, m}) = 0$, by taking the partial derivative with respect to $\lambda_j$ for $j \in \{m+1, \cdots, n-1\}$, we obtain
\begin{align*}
\label{eq:2.16}
&\kern-2em \frac{\lambda_{m+1} \cdots \lambda_n}{t \lambda_j} -  \sum_{k=m+1}^{n-2} c_k \sigma_{k-m-1}(\lambda_{; 1, \cdots, m, j})  \tag{2.16} \\
&= - \frac{\partial \lambda_n}{\partial \lambda_j} \Bigl (  \frac{\lambda_{m+1} \cdots \lambda_{n-1}}{t}  -  \sum_{k=m+1}^{n-2} c_k \sigma_{k-m-1}(\lambda_{; 1, \cdots, m, n})       \Bigr). 
\end{align*}
By combining equations (\ref{eq:2.15}) and (\ref{eq:2.16}), we have
\begin{align*}
\frac{\partial}{\partial \lambda_j}  \Bigl(  \sum_{k=m-1}^{n-2} c_k \sigma_{k-m+1}(\lambda_{; 1, \cdots, m}) \Bigr) 
&= (\lambda_j - \lambda_n) \Bigl (       \frac{\lambda_{m+1} \cdots \lambda_n}{ t \lambda_j} -  \sum_{k=m+1}^{n-2} c_k \sigma_{k-m-1}(\lambda_{; 1, \cdots, m, j})       \Bigr) \geq 0.
\end{align*}Since we assume that $(c_{n-2}, \cdots, c_m) \neq 0$, by Lemma~\hyperlink{L:2.6}{2.6}, we have \[ \frac{1}{t}{\lambda_{m+1} \cdots \lambda_n}/\lambda_j -  \sum_{k=m+1}^{n-2} c_k \sigma_{k-m-1}(\lambda_{; 1, \cdots, m, j})  = 0\] if and only if $\lambda_{m+1} = \cdots = \lambda_n = x_m(t) = x_{m+1}(t)$. If $(\lambda_{m+1}, \cdots \lambda_n) \neq (x_m(t), \cdots, x_m(t))$, we have $   \frac{1}{t}{\lambda_{m+1} \cdots \lambda_n}/\lambda_j -  \sum_{k=m+1}^{n-2} c_k \sigma_{k-m-1}(\lambda_{; 1, \cdots, m, j})  > 0$ and this implies that the function $\sum_{k=m-1}^{n-2} c_k \sigma_{k-m+1}(\lambda_{; 1, \cdots, m-1})$ attains its global minimum at $\lambda_{m+1} = \cdots = \lambda_n$. That is to say, 
\begin{align*}
 \sum_{k=m-1}^{n-2} c_k \sigma_{k-m+1}(\lambda_{; 1, \cdots, m-1}) \geq \sum_{k = m-1}^{n-2} c_k \mybinom[0.8]{n-m}{k-m+1} x^{k-m+1}_{m}(t)  \geq \sum_{k = m-1}^{n-2} c_k \mybinom[0.8]{n-m}{k-m+1} x^{k-m+1}_{m}.
\end{align*}
Moreover, if $(c_{n-2}, \cdots, c_{m}) \neq 0$, then by Lemma~\hyperlink{L:2.3}{2.3} and Lemma~\hyperlink{L:2.5}{2.5}, we have
\[ \sum_{k = m-1}^{n-2} c_k \sigma_{k-m+1}(\lambda_{; 1, \cdots,   m}) = \sum_{k = m-1}^{n-2} c_k \mybinom[0.8]{n-m}{k-m+1} x^{k-m+1}_{m}\] 
if and only if $\lambda_{m+1} = \cdots = \lambda_n = x_m(t)$ and $t = 1$. By Lemma~\hyperlink{L:2.6}{2.6}, we have
$x_0 = \cdots = x_{m}$ and $\lambda_{m+1} = \cdots = \lambda_{n} = x_{m}$.
\smallskip

Then, we consider the quantity $\sum_{k = m-1}^{n-2} c_k \sigma_{k-m+1}(\lambda_{; 1, \cdots, m-1})$. By mathematical induction,  
\begin{align*}
\sum_{k = m-1}^{n-2} c_k \sigma_{k-m+1}(\lambda_{; 1, \cdots, m-1}) = \sum_{k = m-1}^{n-2} c_k \sigma_{k-m+1}(\lambda_{; 1, \cdots, m}) + \lambda_m \sum_{k = m}^{n-2} c_k \sigma_{k-m}(\lambda_{; 1, \cdots, m})  \geq 0.
\end{align*}There are two cases to consider. For the first case, if $\sum_{k = m-1}^{n-2} c_k \sigma_{k-m+1}(\lambda_{; 1, \cdots, m-1})= 0$, then we immediately get $ \sum_{k = m-1}^{n-2} c_k \sigma_{k-m+1}(\lambda_{; 1, \cdots, m}) = 0$, which implies that $x_{m-1} = x_m$. In addition, if $\lambda_m > 0$, then $\sum_{k = m}^{n-2} c_k \sigma_{k-m}(\lambda_{; 1, \cdots, m}) = 0$, which implies that $x_m = \cdots = x_{n-2} = 0$. Otherwise, if $\lambda_m = 0$, then by Lemma~\hyperlink{L:2.6}{2.6}, we have $(c_{n-2}, \cdots, c_m) = 0$, which also implies that $x_m = \cdots = x_{n-2} = 0$. In conclusion, if $\sum_{k = m-1}^{n-2} c_k \sigma_{k-m+1}(\lambda_{; 1, \cdots, m-1})= 0$, then $x_{m-1} = \cdots = x_{n-2} = 0$.\smallskip

For the second case, if $\sum_{k = m-1}^{n-2} c_k \sigma_{k-m+1}(\lambda_{; 1, \cdots, m-1}) > 0$, then there exists a $t \geq 1$ such that
\begin{align*}
\lambda_m \cdots \lambda_n - t \sum_{k = m-1}^{n-2} c_k \sigma_{k-m+1}(\lambda_{; 1, \cdots, m-1}) = 0.
\end{align*} 
By Lemma~\hyperlink{L:2.6}{2.6}, this implies that $t\sum_{k = m-1}^{n-2} c_k \sigma_{k-m+1}(\lambda_{; 1, \cdots, m-1}) \geq x_{m-1}^{n-m+1}(t)$.
Furthermore, by Proposition~\hyperlink{P:2.2}{2.2}, we have $x_{m-1}(t) \geq t^{\frac{1}{n-m+1}} x_{m-1}$. Combine these, we get 
\[ \sum_{k = m-1}^{n-2} c_k \sigma_{k-m+1}(\lambda_{; 1, \cdots, m-1}) \geq x_{m-1}^{n-m+1}. \]
Moreover, if $(c_{n-2}, \cdots, c_{m}) \neq 0$, then by Lemma~\hyperlink{L:2.3}{2.3}, Lemma~\hyperlink{L:2.5}{2.5}, and Lemma~\hyperlink{L:2.6}{2.6}, we have $\sum_{k = m-1}^{n-2} c_k \sigma_{k-m+1}(\lambda_{; 1, \cdots,   m-1}) = \sum_{k = m-1}^{n-2} c_k \binom{n-m+1}{k-m+1} x^{k-m+1}_{m-1}$ if and only if $\lambda_{m} = \cdots = \lambda_n = x_{m-1}(t)$ and $t = 1$. By Lemma~\hyperlink{L:2.6}{2.6}, we have
$x_0 = \cdots = x_{m-1}$ and $\lambda_{m} = \cdots = \lambda_{n} = x_{m-1}$.
This finishes the proof.
\end{proof}

\subsection{Basic Formulas of Symmetric Functions}
\label{sec:2.4}
In this subsection, we state some standard lemmas for symmetric functions. Let $(M, \omega)$ be a compact connected Kähler manifold of complex dimension $n$ with a Kähler form $\omega$ and $[X] \in H^{1,1}(M; \mathbb{R})$. We consider the Hermitian endomorphism $\Lambda = \omega^{-1} X_u = \omega^{-1} (X + \sqrt{-1} \partial \bar{\partial} u)$ at a fixed point $p \in M$.

\hypertarget{L:2.8}{\begin{flemma}
If $F(\Lambda) = f(\lambda_1, \dots, \lambda_n)$ is a smooth function in the eigenvalues $\{ \lambda_1, \cdots, \lambda_n\}$ of a Hermitian matrix $\Lambda$, then at a diagonal matrix $\Lambda$ with distinct eigenvalues $\lambda_i$, we get
\begin{align*}
\frac{\partial F}{\partial \Lambda_{i }^j}(\Lambda) = \delta_{ij} f_i(\lambda);\quad  \frac{\partial^2 F}{\partial \Lambda_{i }^j \Lambda_{r  }^s} (\Lambda) = f_{ir}(\lambda) \delta_{ij} \delta_{rs} + \frac{f_i - f_j}{ \lambda_i - \lambda_j}(\lambda) (1 - \delta_{ij})\delta_{is} \delta_{jr}, 
\end{align*}where $f_i(\lambda) =  \frac{\partial f}{\partial \lambda_i}(\lambda)$ and $f_{ir} = \frac{\partial^2 f}{\partial \lambda_i \partial \lambda_r} (\lambda)$.
\end{flemma}}

We denote $\lambda$ the eigenvalues of the Hermitian endomorphism $\omega^{i \bar{k}} \bigl (  X + \sqrt{-1} \partial \bar{\partial}  {u}   \bigr )_{j \bar{k}}$. Since we are on a Kähler manifold, we can pick the following local coordinates to simplify our computation.

\hypertarget{L:2.9}{\begin{flemma}At any point $p \in M$, there exist local holomorphic coordinates near $p$ such that 
\begin{align*}
\omega_{i \bar{j}} (p) = \delta_{ij}; \quad     ( X_u     )_{i \bar{j}} (p) = \lambda_i \delta_{ij};  \quad \omega_{i \bar{j}, k} (p) = 0,   
\end{align*}for all $i, j, k \in \{ 1, \dots, n\}$.
\end{flemma}}

From now on, without further notice, we always use the above coordinates. We denote $\Lambda$ the Hermitian endomorphism $\omega^{i \bar{k}} \bigl (  X + \sqrt{-1} \partial \bar{\partial}  {u}   \bigr )_{j \bar{k}}$. Then the first and second derivatives of $\Lambda$ will be the following.

\hypertarget{L:2.10}{\begin{flemma}The first and second derivatives of $\Lambda$ are
\begin{align*}
\frac{\partial \Lambda_i^j }{\partial \bar{z}_k} &=  \tensor[]{\omega}{^{j \bar{p}}_{, \bar{k}}}     (  X_u      )_{i \bar{p}} +  \tensor[]{\omega}{^{j \bar{p}}}     (  X_u      )_{i \bar{p},\bar{k}} =  - \omega^{j \bar{b}} \omega_{a \bar{b}, \bar{k}} \omega^{a \bar{p}}     (  X_u      )_{i \bar{p}} +  \tensor[]{\omega}{^{j \bar{p}}}    (  X_u      )_{i \bar{p},\bar{k}},  \\ 
\frac{\partial^2 \Lambda_i^j }{\partial z_l \partial \bar{z}_k} &= \tensor[]{\omega}{^{j \bar{p}}_{, \bar{k}l}}    (  X_u      )_{i \bar{p}} + \tensor[]{\omega}{^{j \bar{p}}_{, \bar{k}}}     (  X_u      )_{i \bar{p},l} +  \tensor[]{\omega}{^{j \bar{p}}_{, l}}     (  X_u      )_{i \bar{p},\bar{k}} + \tensor[]{\omega}{^{j \bar{p}}}    (  X_u      )_{i \bar{p},\bar{k}l}   \\
&=  \omega^{j \bar{d}} \omega_{c \bar{d}, l} \omega^{c \bar{b}} \omega_{a \bar{b}, \bar{k}} \omega^{a \bar{p}}    (  X_u      )_{i \bar{p}} - \omega^{j \bar{b}} \omega_{a \bar{b}, \bar{k} l} \omega^{a \bar{p}}     (  X_u      )_{i \bar{p}} + \omega^{j \bar{b}} \omega_{a \bar{b}, \bar{k}} \omega^{a \bar{d}} \omega_{c \bar{d}, l} \omega^{c \bar{p}}    (  X_u      )_{i \bar{p}} \\
&\kern2em - \omega^{j \bar{b}} \omega_{a \bar{b}, \bar{k}} \omega^{a \bar{p}}    (  X_u      )_{i \bar{p},l} - \omega^{j \bar{b}} \omega_{a \bar{b}, l} \omega^{a \bar{p}}     (  X_u      )_{i \bar{p},\bar{k}} + \omega^{j \bar{p}}    (  X_u      )_{i \bar{p},\bar{k}l} , 
\end{align*}where we denote $   ( X_u      )_{i \bar{j}} = X_{i \bar{j}} + u_{i \bar{j}}$ and $\Lambda$ is the Hermitian endomorphism $\omega^{-1}   ( X_u     )$. 
\end{flemma}}

If we evaluate at any fixed point $p \in M$ and we use the coordinates in Lemma~\hyperlink{L:2.9}{2.9}, we can simplify the first and second derivatives of $\Lambda$.

\hypertarget{L:2.11}{\begin{flemma}At any fixed point $p$, by picking the coordinates in Lemma~\hyperlink{L:2.9}{2.9}, we get
\begin{align*}
\frac{\partial \Lambda_i^j}{\partial \bar{z}_k} (p) =      (  X_u      )_{i \bar{j},\bar{k}}; \quad \frac{\partial^2 \Lambda_i^j}{\partial z_l \partial \bar{z}_k} (p) =   - \lambda_i \omega_{i \bar{j}, \bar{k} l}    +     (  X_u       )_{i \bar{j},\bar{k}l}.
%&= \tilde{\lambda}_i \omega_{i \bar{a}, \bar{k}} \omega_{a \bar{j}, l} + \tilde{\lambda}_i \omega_{a \bar{j}, \bar{k}} \omega_{i \bar{a}, l} - \tilde{\lambda}_i \omega_{i \bar{j}, \bar{k} l}    - \omega_{a \bar{j}, \bar{k}}  \tensor[]{\hat{\chi}}{_{i \bar{a}, l}} - \omega_{a \bar{j}, l}  \tensor[]{\hat{\chi}}{_{i \bar{a}, \bar{k}}}  +  \tensor[]{\hat{\chi}}{_{i \bar{j}, \bar{k}l}}.   
\end{align*}
\end{flemma}}

\section{A Priori Estimates}
\label{sec:3}
In this section, we study a priori estimates for the following map
\begin{align*}
c \colon M^n \longrightarrow \tilde{\mathscr{C}}_{n} \text{ given by } c(z) =  (c_{n-2}, \cdots, c_{1}, c_0(z)  )
\end{align*}
provided that a $C$-subsolution to $c \colon M^n \longrightarrow \tilde{\mathscr{C}}_{n}$ is given. Some of the ideas in this section come from the author's previous works \cite{lin2023d, lin2023c, lin2023thesis}. The author was unable to find a unified approach due to the lack of a better understanding of $\Upsilon$-stable general inverse $\sigma_k$ multilinear polynomials and right-Noetherian polynomials. The author had to study the generic stratas of $\tilde{\mathscr{C}}_3$ and $\tilde{\mathscr{C}}_4$ separately and the arguments were complicated.
In this paper, we find a unified approach to handle all dimensions simultaneously. Another important thing the author would like to point out is that the following a priori estimates also hold for other stratas of different dimensions not just the Calabi--Yau set $\CY_n$ and the generic strata of $\tilde{\mathscr{C}}_n$. Though it is still unknown what roles these stratas of different dimensions play on Kähler manifolds or projective manifolds.\smallskip

First, let us summarize the proof of our a priori estimates. Under the assumption of $C$-subsolution, we apply the Alexandroff--Bakelman--Pucci estimate to get a $C^0$ estimate. This $C^0$ estimate can be obtained following the proof in Székelyhidi \cite{szekelyhidi2018fully}, which is based on the method that Błocki \cite{blocki2005uniform, blocki2011uniform} used in the case of the complex Monge--Ampère equation. We will skip the proof of this $C^0$ estimate because it follows verbatim.

Second, we use the maximum principle to obtain that the $C^2$ norm can be bounded by the $C^1$ norm. The method is inspired by Hou--Ma--Wu \cite{hou2010second} for the complex Hessian equations and used by Székelyhidi \cite{szekelyhidi2018fully}. The interested reader is referred to \cite{collins20201, szekelyhidi2018fully} and the references therein. Once we have the above type inequality, by a blow-up argument due to Dinew--Kołodziej \cite{dinew2017liouville}, we can get an indirect $C^1$ estimate. 

Last, to get $C^{2, \alpha}$ estimate, we follow the proof of the complex version of the Evans--Krylov theory in Siu \cite{siu2012lectures}, we can exploit the convexity of the solution sets to obtain a $C^{2, \alpha}$ estimates by a blow-up argument. Furthermore, for higher regularity, we apply the standard Schauder estimates and bootstrapping. \smallskip

To start with, we consider the following equation on $M^n$:
\begin{align*}
	\label{eq:3.1}
    	X^n =  d_{n-2} \mybinom[0.8]{n}{n-2} X^{n-2} \wedge \omega^2 + \cdots +  \mybinom[0.8]{n}{1} d_1 X \wedge \omega^{n-1} +  \mybinom[0.8]{n}{0} d_0(z)   \omega^{n} , \tag{3.1}
\end{align*}where $d_1, \cdots, d_{n-2}$ are constants and $d_0$ is a function on $M$ satisfying $d(z) =  (d_{n-2}, \cdots, d_{1}, d_0(z)  )  \in \tilde{\mathscr{C}}_{n}$ for all $z \in M$. So equation (\ref{eq:3.1}) can be rewritten as a map $d \colon M^n \rightarrow \tilde{\mathscr{C}}_{n}$. \smallskip

Throughout this section, first, we always assume that there exists a $C$-subsolution $\underline{u} \colon M   \rightarrow \mathbb{R}$ to the map $d \colon M^n \rightarrow \tilde{\mathscr{C}}_n$. Then, we also call $X_{\underline{u}} \coloneqq X + \sqrt{-1} \partial \bar{\partial} \underline{u}$ this $C$-subsolution and by changing representative, we may assume $X$ is this $C$-subsolution. In Section~\ref{sec:4}, we want to apply the method of continuity to obtain the solvability of equation (\ref{eq:3.1}), so we are interested in maps in the set $\tilde{\mathcal{P}}^d$. For any $c \in \tilde{\mathcal{P}}^{d}$, by Remark~\hyperlink{R:2.4}{2.4}, $X$ is also a $C$-subsolution to $c$.
In this paper, we only consider maps $c \in \tilde{\mathcal{P}}^{d}$ with $c_1, \cdots, c_{n-2}$ constants. With such $c$, we consider the following equation:
\begin{align*}
\label{eq:3.2}
h_c(z, \lambda) \coloneqq \frac{ \sum_{k=1}^{n-2} c_k  \sigma_k(\lambda) + c_0(z)}{\lambda_1 \cdots \lambda_n} = 1, \tag{3.2}
\end{align*}where $z \in M$ and $\lambda_i$ are the eigenvalues of $\omega^{-1}X_u$ at $z$. Note that we abbreviate $\lambda = \{\lambda_1, \cdots, \lambda_n\}$ and we always assume $\lambda_1 \geq \cdots \geq \lambda_n$ unless further notice. Most of the time, to save spaces, we will abbreviate $h = h_c(\lambda)$, $h_i = \partial h_c /\partial \lambda_i, h_{ij} = \partial^2 h_c/\partial \lambda_i \partial \lambda_j$ for $i, j \in \{1, \cdots, n\}$ for notational convention. 
Unless otherwise specified, in this paper, we always assume that the maps $c \in \tilde{\mathcal{P}}^d$ are those with $c_{1}, \cdots, c_{n-2}$ constants and range in a compact subset $S$ of $\tilde{\mathscr{C}}_n$.

\subsection{The \texorpdfstring{$C^2$}{} Estimate}
\label{sec:4.1.1}
Define a Hermitian endomorphism $\Lambda \coloneqq \omega^{-1}X_u$, where $X_u = X + \sqrt{-1} \partial \bar{\partial} u$, and let $\lambda = \{ \lambda_1, \cdots, \lambda_n \}$ be the eigenvalues of $\Lambda$. We consider the following function $G(\Lambda) = \log(1 + \lambda_1) =g(\lambda_1, \cdots, \lambda_n)$ and the following test function
\begin{align*}
\label{eq:3.3}
U \coloneqq - Au + G(\Lambda),  \tag{3.3}
\end{align*}where $A \gg 0$ will be determined later. We want to apply the maximum principle to $U$, but since the eigenvalues of $\Lambda$ might not be distinct at the maximum point $q \in M$ of $U$, we do a perturbation here. Assume $\lambda_1$ is large, otherwise, we are done, then 
\begin{itemize}
\hypertarget{pert for all dim}{\item} we pick the constant matrix $B$ to be a diagonal matrix with real entries 
\begin{align*}
B_{11} = \lambda_1;\quad B_{22} = \lambda_2/2;\quad \cdots;\quad  B_{n-1 \ n-1} = \lambda_{n-1}/(n-1);\quad  B_{nn} = 0.
\end{align*}
\end{itemize}
By defining $\tilde{\Lambda} = \Lambda + B$, then $\tilde{\Lambda}$ has distinct eigenvalues near $q \in M$. The eigenvalues of $\tilde{\Lambda}$ define smooth functions near the maximum point $q$. And we can check that the maximum point $q$ of $U$ in equation (\ref{eq:3.3}) is still the maximum point of the following locally defined test function
\begin{align*}
\label{eq:3.4}
\tilde{U} \coloneqq -Au + G(\tilde{\Lambda}). \tag{3.4}
\end{align*}
Near the maximum point $q$ of $\tilde{U}$, we always use the coordinates in Lemma~\hyperlink{L:2.9}{2.9} unless otherwise noted. We instantly get the following lemma. 
\hypertarget{L:3.1}{\begin{flemma}
At the maximum point $q$ of $\tilde{U}$, by taking the first derivative of $\tilde{U}$ at $q$, we get
\begin{align*}
\label{eq:3.5}
0 &= -A u_k(q) + \frac{1}{1+ \tilde{\lambda}_1} ( X_u )_{1\bar{1}, k}, \tag{3.5}
%0 &\geq   -A u_{k\bar{k}} (p) - \frac{1}{(1+ \lambda_1)^2} ( X_u )_{1\bar{1}, k},
\end{align*}where we denote $u_k = \partial u/\partial z_k$ and $( X_u )_{1\bar{1}, k} = \partial (X_u)_{1 \bar{1}} / \partial z_k$.
\end{flemma}}

\begin{proof}
First, since $\tilde{U} = -Au + G(\tilde{\Lambda})$, if we take the first derivative, we obtain
\begin{align*}
\frac{\partial}{\partial z_k} \tilde{U} &= - A u_k + \frac{\partial G}{\partial \Lambda^j_i} (\tilde{\Lambda}) \frac{\partial \tilde{\Lambda}^j_i}{\partial z_k}.
\end{align*}At the maximum point $q$, we have $0 = - A u_k(q) + \frac{1}{1+ \tilde{\lambda}_1} ( X_u )_{1\bar{1}, k}$, which finishes the proof.
\end{proof}

For any $c \in \tilde{\mathcal{P}}^d$ with $c_{1}, \cdots, c_{n-2}$ constants, we may define the following operator \hypertarget{dim n operator}{$\mathcal{L}_c$} by
\begin{align*}
\label{eq:3.6}
\mathcal{L}_c \coloneqq - \sum_{i, j, k} \frac{\partial H_c}{\partial \Lambda^k_i} (z, {\Lambda}) \omega^{k \bar{j}}(z)  \frac{\partial^2}{\partial z_i \partial \bar{z}_j}, \tag{3.6}
\end{align*}where $\Lambda$ is the Hermitian endomorphism $\omega^{-1} X_u$ at $z$ and $H_c(z, \Lambda)   = h_c(z, {\lambda}_1, \cdots,  {\lambda}_n)$ is defined by $h_c(z, {\lambda}) =   \sum_{k=0}^{n-2} c_k  \sigma_k( {\lambda})   /  {\lambda}_1  \cdots   {\lambda}_n$. We immediately have the following.

\hypertarget{L:3.2}{\begin{flemma}
By taking $h_c( {\lambda}) =   \sum_{k=0}^{n-2} c_k  \sigma_k( {\lambda})   /  \sigma_n(\lambda)$ and $g( {\lambda})=\log(  1 +  {\lambda}_1)$, we have 
\begin{align*}
h_i &= - \frac{    \sum_{k=0}^{n-2} c_k  \sigma_k (  {\lambda}_{;i})    }{   {\lambda}_1  \cdots   {\lambda}_n   {\lambda}_i }; &h_{ij}& = \frac{  \sum_{k=0}^{n-2} c_k \sigma_k(\lambda_{; i, j})   }{  {\lambda}_1  \cdots   {\lambda}_n   {\lambda}_i   {\lambda}_j} (1+ \delta_{ij}); \\
g_i &= \delta_{1i} \frac{1}{ 1 +   {\lambda}_1}; &g_{ij}&=-\delta_{1i}\delta_{1j}\frac{1}{( 1 +  {\lambda}_1)^2}. 
\end{align*}Here, we denote $h_i \coloneqq  {\partial  h_c}/{\partial \lambda_i}, \, g_i \coloneqq  {\partial g}/{\partial \lambda_i}$, $h_{ij} \coloneqq  {\partial^2 h_c}/{\partial \lambda_i \partial \lambda_j}$, and $g_{ij} \coloneqq  {\partial^2 g}/{\partial \lambda_i \partial \lambda_j}$. 
\end{flemma}}

\hypertarget{L:3.3}{\begin{flemma}
Let $c \in \tilde{\mathscr{C}}_n$, then for any point on the set $\{   h = 1\}$, we have
\begin{align*}
-h_i = \frac{    \sum_{k=0}^{n-2} c_k  \sigma_k (  {\lambda}_{;i})  }{   {\lambda}_1  \cdots   {\lambda}_n   {\lambda}_i } = \frac{1}{\lambda_i} -  \frac{\sum_{k=1}^{n-2} c_k \sigma_{k-1} (\lambda_{; i}) }{\lambda_1 \cdots \lambda_n}  > 0.
\end{align*}
\end{flemma}}
\begin{proof}
On $h  = 1$, this implies that $\lambda_i \bigl(  \lambda_1 \cdots \lambda_n /\lambda_i  -  \sum_{k=1}^{n-2} c_k \sigma_{k-1}(\lambda_{; i}) \bigr) = \sum_{k=0}^{n-2} c_k \sigma_k(\lambda_{;i})$. Since $c \in \tilde{\mathscr{C}}_n$, we have $ \lambda_1 \cdots \lambda_n/\lambda_i  -  \sum_{k=1}^{n-2} c_k \sigma_{k-1}(\lambda_{; i}) > 0$ on $\{ \lambda_1 \cdots \lambda_n  -  \sum_{k=0}^{n-2} c_k \sigma_k(\lambda) = 0\}$. This implies that $ \sum_{k=0}^{n-2} c_k \sigma_k(\lambda_{;i}) > 0$ and $-h_i > 0$. This finishes the proof.
\end{proof}

Assuming that $\lambda_1 \geq \cdots \lambda_n$, then we have the following upper bound estimate of $h_n$ which plays an important role in getting a unified approach to all dimensions.
\hypertarget{L:3.4}{\begin{flemma}
    Let $c \in \tilde{\mathscr{C}}_n$, $h(\lambda)= \sum_{k=0}^{n-2} c_k \sigma_k ( {\lambda})  /  {\lambda}_1   \cdots   {\lambda}_n$, and assume that $\lambda_1 \geq \cdots \geq \lambda_n$, then we have the following estimates of $h_n$ on $\{  h = 1\}$:
    \begingroup
    \allowdisplaybreaks
    \begin{align*}
        0 > -\frac{1}{\lambda_n}\Bigl( 1 -     {\sum_{k=1}^{n-2} c_k \mybinom[0.8]{n-1}{k-1} {x}_{0}^{k-n}}    \Bigr) \geq h_n \geq  -\frac{1}{\lambda_n},
    \end{align*}
    \endgroup
    where ${x}_0$ is the largest real root of the strictly right-Noetherian polynomial $x^{n} -  \sum_{k=0}^{n-2} c_k \binom{n}{k} x^{k}$.
\end{flemma}}
\begin{proof}
We can rewrite $h_n$ as
\begin{align*}
h_n = - \frac{\sum_{k=0}^{n-2} c_k \sigma_k(\lambda_{; n}) }{\lambda_1 \cdots \lambda_n \lambda_n} = -\frac{1}{\lambda_n} + \frac{\sum_{k=1}^{n-2} c_k \sigma_{k-1} (\lambda_{;n})}{\lambda_1 \cdots \lambda_n}.
\end{align*}
There are two cases to consider. For the case $\sum_{k=1}^{n-2} c_k \sigma_{k-1}(\lambda_{; n}) = 0$, by Lemma~\hyperlink{L:2.7}{2.7}, we have $(c_{n-2}, \cdots, c_1) = 0$. So, we immediately get $\sum_{k=1}^{n-2} c_k \binom{n-1}{k-1} {x}_{0}^{k-n} = 0$ and $h_n = -  {1}/{\lambda_n}$. For the case $\sum_{k=1}^{n-2} c_k \sigma_{k-1}(\lambda_{; n}) > 0$, by fixing $\lambda_n > 0$, we view $(\lambda_1, \cdots, \lambda_{n-1})$ as a point on 
\begin{align*}
\label{eq:3.7}
\lambda_1 \cdots \lambda_{n-1}  - \sum_{k=0}^{n-2} \bigl(  c_{k+1} + \frac{c_k}{\lambda_n} \bigr) \sigma_k(\lambda_{; n}) = 0. \tag{3.7}
\end{align*}Here, $\lambda_1 \cdots \lambda_{n-1}  - \sum_{k=0}^{n-2} \bigl(  c_{k+1} + \frac{c_k}{\lambda_n} \bigr) \sigma_k(\lambda_{; n})$ is a degree $n-1$ strictly $\Upsilon$-stable general inverse $\sigma_k$ multilinear polynomial. Now, we consider the quantity $\frac{\sum_{k=1}^{n-2} c_k \sigma_{k-1}(\lambda_{; n})}{\lambda_1 \cdots \lambda_{n-1}}$ on (\ref{eq:3.7}), we get
\begin{align*}
\label{eq:3.8}
\frac{\partial}{\partial \lambda_i} \Biggl( \frac{\sum_{k=1}^{n-2} c_k \sigma_{k-1}(\lambda_{; n})}{\lambda_1 \cdots \lambda_{n-1}} \Biggr) = - \frac{ \lambda_j \sum_{k=1}^{n-2}c_k\sigma_{k-1}(\lambda_{; i, n}) + \lambda_i \frac{\partial \lambda_j}{\partial \lambda_i} \sum_{k=1}^{n-2}c_k\sigma_{k-1}(\lambda_{; j, n} )}{\lambda_1 \cdots \lambda_{n-1} \lambda_i \lambda_j}. \tag{3.8}
\end{align*}
On the other hand, on $\lambda_1 \cdots \lambda_n - \sum_{k=0}^{n-2} c_k \sigma_k(\lambda) = 0$, we have
\begin{align*}
\label{eq:3.9}
0 &= \frac{1}{\lambda_i} \bigl( \lambda_1 \cdots \lambda_n - \lambda_i \sum_{k=1}^{n-2} c_k \sigma_{k-1}(\lambda_{; i}) \bigr) + \frac{1}{\lambda_j}  \frac{\partial \lambda_j}{\partial \lambda_i} \bigl( \lambda_1 \cdots \lambda_n - \lambda_j \sum_{k=1}^{n-2} c_k \sigma_{k-1}(\lambda_{; i}) \bigr) \tag{3.9} \\
&= \frac{1}{\lambda_i} \sum_{k=0}^{n-2} c_k \sigma_k(\lambda_{; i}) + \frac{1}{\lambda_j} \frac{\partial \lambda_j}{\partial \lambda_i}  \sum_{k=0}^{n-2} c_k \sigma_k(\lambda_{; j}). 
 \end{align*}By combining equations (\ref{eq:3.8}) and (\ref{eq:3.9}), we obtain
 \begin{align*}
 \label{eq:3.10}
&\kern-2em  \frac{\partial}{\partial \lambda_i} \Biggl( \frac{\sum_{k=1}^{n-2} c_k \sigma_{k-1}(\lambda_{; n})}{\lambda_1 \cdots \lambda_{n-1}} \Biggr) \tag{3.10} \\
&= - \frac{\sum_{k=0}^{n-2} c_k \sigma_k (\lambda_{;j})  \sum_{k=1}^{n-2}c_k\sigma_{k-1}(\lambda_{; i, n}) - \sum_{k=0}^{n-2} c_k \sigma_k (\lambda_{;i}) \sum_{k=1}^{n-2}c_k\sigma_{k-1}(\lambda_{; j, n} )}{\lambda_1 \cdots \lambda_{n-1} \lambda_i  \sum_{k=0}^{n-2} c_k \sigma_k (\lambda_{;j}) }.
 \end{align*}For the numerator of equation (\ref{eq:3.10}), if we treat $\{\lambda_1, \cdots, \lambda_n\}$ as variable, then by taking the partial derivative with respect to $\lambda_n$, we get
 \begin{align*}
 &\kern-2em \frac{\partial}{\partial \lambda_n} \Bigl( \sum_{k=0}^{n-2} c_k \sigma_k (\lambda_{;j})  \sum_{k=1}^{n-2}c_k\sigma_{k-1}(\lambda_{; i, n}) - \sum_{k=0}^{n-2} c_k \sigma_k (\lambda_{;i}) \sum_{k=1}^{n-2}c_k\sigma_{k-1}(\lambda_{; j, n} ) \Bigr) \\
 &= \sum_{k=1}^{n-2} c_k \sigma_{k-1} (\lambda_{; j, n})  \sum_{k=1}^{n-2}c_k\sigma_{k-1}(\lambda_{; i, n}) - \sum_{k=1}^{n-2} c_k \sigma_k (\lambda_{; i, n}) \sum_{k=1}^{n-2}c_k\sigma_{k-1}(\lambda_{; j, n} ) = 0.
 \end{align*}So the value of $\sum_{k=0}^{n-2} c_k \sigma_k (\lambda_{;j})  \sum_{k=1}^{n-2}c_k\sigma_{k-1}(\lambda_{; i, n}) - \sum_{k=0}^{n-2} c_k \sigma_{k-1} (\lambda_{;i}) \sum_{k=1}^{n-2}c_k\sigma_{k-1}(\lambda_{; j, n} )$ is independent of the value of $\lambda_n$. Besides, by Lemma~\hyperlink{L:2.7}{2.7} and $c \in \tilde{\mathscr{C}}_n$, for $\lambda_i \geq \lambda_j \geq \lambda_n$, we have
\begin{align*}
\label{eq:3.11}
\sum_{k=1}^{n-2} c_k \sigma_{k-1}(\lambda_{j, n}) \geq \sum_{k=1}^{n-2} c_k \sigma_{k-1}(\lambda_{i, n}) \geq \sum_{k=1}^{n-2} c_k \sigma_{k-1}(\lambda_{i, j})  > \sum_{k=1}^{n-2} c_k \mybinom[0.8]{n-2}{k-1} x_2^{k-1} \geq 0. \tag{3.11}
\end{align*}
Now, the idea is to decrease the value of $\lambda_n$ and see whether the quantity $\sum_{k=0}^{n-2}c_k \sigma_{k}(\lambda_{; i})$ or the quantity $\sum_{k=0}^{n-2}c_k \sigma_{k}(\lambda_{; j})$ attains $0$ first. Note that these two quantities both share the following variables $\{\lambda_1, \cdots, \lambda_n\}\backslash \{\lambda_i, \lambda_j\}$, so we can treat them as a single function $\sum_{k=0}^{n-2} c_k \sigma_k(\lambda_i)$ evaluated at two points, which only one entry differs.
We consider how $\lambda_n$ changes on the level set of $\sum_{k=0}^{n-2}c_k \sigma_{k}(\lambda_{; i}) = 0$. By taking the partial derivative with respect to $\lambda_n$, we get
\begin{align*}
0 = \frac{\partial}{\partial \lambda_j } \Bigl(   \sum_{k=0}^{n-2}c_k \sigma_{k}(\lambda_{; i})     \Bigr) = \sum_{k = 1}^{n-2} c_k \sigma_{k-1}(\lambda_{; i, j})  + \frac{\partial \lambda_n}{\partial \lambda_j} \sum_{k = 1}^{n-2} c_k \sigma_{k-1}(\lambda_{; i, n}), 
\end{align*} which implies that $\frac{\partial \lambda_n}{\partial \lambda_j} < 0$ by (\ref{eq:3.11}). For $\lambda_i > \lambda_j$, by decreasing the value of $\lambda_n$, we see that $\sum_{k=0}^{n-2} c_k \sigma_k(\lambda_{; i})$ attains $0$ before $\sum_{k=0}^{n-2} c_k \sigma_k(\lambda_{; j})$ attains $0$. Hence, for this new $\lambda_n$, we have $\sum_{k=0}^{n-2} c_k \sigma_k(\lambda_{; j}) > \sum_{k=0}^{n-2} c_k \sigma_k(\lambda_{; i}) = 0$ and (\ref{eq:3.10}) becomes
\begin{align*}
\label{eq:3.12}
\frac{\partial}{\partial \lambda_i} \Biggl( \frac{\sum_{k=1}^{n-2} c_k \sigma_{k-1}(\lambda_{; n})}{\lambda_1 \cdots \lambda_{n-1}} \Biggr) = - \frac{\sum_{k=0}^{n-2} c_k \sigma_k (\lambda_{;j})  \sum_{k=1}^{n-2}c_k\sigma_{k-1}(\lambda_{; i, n})  }{\lambda_1 \cdots \lambda_{n-1} \lambda_i  \sum_{k=0}^{n-2} c_k \sigma_k (\lambda_{;j}) } < 0. \tag{3.12}
\end{align*}The last inequality is due to (\ref{eq:3.11}).
Back to the original fix $\lambda_n$, by inequality (\ref{eq:3.12}), we can decrease the largest $\lambda$'s and increase the smallest $\lambda$'s, the quantity $\frac{\sum_{k=1}^{n-2} c_k \sigma_{k-1}(\lambda_{; n})}{\lambda_1 \cdots \lambda_{n-1}}$ becomes larger. By continuing this process and arguing cautiously, eventually, we see that the quantity $\frac{\sum_{k=1}^{n-2} c_k \sigma_{k-1}(\lambda_{; n})}{\lambda_1 \cdots \lambda_{n-1}}$ attains its maximum when $\lambda_1 = \cdots = \lambda_{n-1}$. That is,
\begin{align*}
\frac{\sum_{k=1}^{n-2} c_k \sigma_{k-1}(\lambda_{; n})}{\lambda_1 \cdots \lambda_{n-1}} \leq \frac{\sum_{k=1}^{n-2} c_k \binom{n-1}{k-1} y_0^{k-1}  }{y_0^{n-1}} =   \sum_{k=1}^{n-2} c_k \mybinom[0.8]{n-1}{k-1} y_0^{k-n},
\end{align*}where $y_0$ is the largest real root of $\lambda_n y^{n-1}  - \lambda_n \sum_{k=1}^{n-2} c_k \sigma_{k-1} \binom{n-1}{k-1} y^{k-1} - \sum_{k=0}^{n-2} c_k \binom{n-1}{k} y^{k}$. Due to the assumption that $x_0 \geq \lambda_n$, we have $y_0 \geq x_0 > x_1 \geq 0$. By taking the partial derivative of $\sum_{k=1}^{n-2} c_k \binom{n-1}{k-1}y_0^{k-n}$ with respect to $y_0$, we get
\begin{align*}
\frac{\partial}{\partial y_0} \Bigl (      \sum_{k=1}^{n-2} c_k \mybinom[0.8]{n-1}{k-1} y_0^{k-n}       \Bigr)  = -(n-1) y_0^{-n} \sum_{k=1}^{n-2} c_k \mybinom{n-2}{k-1} y_0^{k-1} < 0.
\end{align*}This implies that we have the following upper bound:
\begin{align*}
\frac{\sum_{k=1}^{n-2} c_k \sigma_{k-1}(\lambda_{; n})}{\lambda_1 \cdots \lambda_{n-1}} \leq \sum_{k=1}^{n-2} c_k \mybinom[0.8]{n-1}{k-1} y_0^{k-n} \leq \sum_{k=1}^{n-2} c_k \mybinom[0.8]{n-1}{k-1} x_0^{k-n} < 1,
\end{align*}where the last inequality is due to the fact that $x_0 > x_1$.
This finishes the proof.
\end{proof}

With Lemma~\hyperlink{L:3.4}{3.4}, we can get the following uniform estimates.

\hypertarget{L:3.5}{\begin{flemma}
Let $c \in \tilde{\mathscr{C}}_n$, $h(\lambda)= \sum_{k=0}^{n-2} c_k \sigma_k ( {\lambda})  /  {\lambda}_1   \cdots   {\lambda}_n$, and assume that $\lambda_1 \geq \cdots \geq \lambda_n$, then we have the following estimates on $\{  h = 1\}$:    
\begingroup
\allowdisplaybreaks
\begin{align*}
\label{eq:3.13}
& \frac{n}{\lambda_n} \geq -\sum_i h_i  >  \frac{1}{\lambda_n}\Bigl( 1 -    {\sum_{k=1}^{n-2} c_k \mybinom[0.8]{n-1}{k-1} {x}_{0}^{k-n}} \Bigr) ;  \tag{3.13}  \\
\label{eq:3.14}
& \frac{n}{\lambda_n^2} \geq \sum_i h_i^2  >  \frac{1}{\lambda_n^2} \Bigl( 1 -      {\sum_{k=1}^{n-2} c_k \mybinom[0.8]{n-1}{k-1} x_0^{k-n} }  \Bigr)^2; \tag{3.14} \\
\label{eq:3.15}
&\frac{n-1}{\lambda_i \lambda_n^2}> \sum_j h_{i j} h_j \geq  -\frac{n+1}{\lambda_i \lambda_n^2}; \tag{3.15} \\
\label{eq:3.16}
&\frac{n(n+1)}{\lambda_n^4} \geq \sum_{i, j} h_i h_{ij} h_j  >  -\frac{n(n-1)}{\lambda_n^4}. \tag{3.16}
\end{align*}
\endgroup
\end{flemma}}

\begin{proof}
First, we have $h_i = -  {\sum_{k=0}^{n-2} c_k \sigma_k(\lambda_{; i}) }/{\lambda_1 \cdots \lambda_n \lambda_i}$ and this implies that
\begin{align*}
h_i - h_j &= - \frac{\sum_{k=0}^{n-2} c_k \sigma_k(\lambda_{; i}) }{\lambda_1 \cdots \lambda_n \lambda_i} +  \frac{\sum_{k=0}^{n-2} c_k \sigma_k(\lambda_{; j}) }{\lambda_1 \cdots \lambda_n \lambda_j} = \frac{ - \lambda_j \sum_{k=0}^{n-2} c_k \sigma_k(\lambda_{; i}) + \lambda_i \sum_{k=0}^{n-2} c_k \sigma_k(\lambda_{; j}) }{\lambda_1 \cdots \lambda_n \lambda_i \lambda_j} \\
%&= \frac{- \lambda_j \bigl( \sum_{k=0}^{n-2} c_k \sigma_k(\lambda) - \lambda_i \sum_{k=1}^{n-2} c_k \sigma_{k-1}(\lambda_{; i})    \bigr) + \lambda_i \bigl( \sum_{k=0}^{n-2} c_k \sigma_k(\lambda) - \lambda_j \sum_{k=1}^{n-2} c_k \sigma_{k-1}(\lambda_{; j})    \bigr) }{\lambda_1 \cdots \lambda_n \lambda_i \lambda_j} \\
%&= \frac{\lambda_i - \lambda_j}{\lambda_i \lambda_j} + \frac{  \sum_{k=1}^{n-2} c_k \sigma_{k-1}(\lambda_{; i})  -  \sum_{k=1}^{n-2} c_k \sigma_{k-1}(\lambda_{; j})   }{\lambda_1 \cdots \lambda_n  } \\
&= \frac{\lambda_i - \lambda_j}{\lambda_1 \cdots \lambda_n \lambda_i \lambda_j} \Bigl(  \lambda_1 \cdots \lambda_n -  \lambda_i \lambda_j \sum_{k = 2}^{n-2} c_k \sigma_{k-2}(\lambda_{; i, j})      \Bigr).
\end{align*}Hence, we have $0 > h_1 \geq \cdots \geq h_n$. Thus, we have the following upper and lower bound 
\begin{align*}
\label{eq:3.17}
h_n > \sum_i h_i  \geq n h_n. \tag{3.17}
\end{align*}
By Lemma~\hyperlink{L:3.4}{3.4} and inequality (\ref{eq:3.17}), we immediately get 
\begin{align*}
-\frac{1}{\lambda_n}\Bigl( 1 -    {\sum_{k=1}^{n-2} c_k \mybinom[0.8]{n-1}{k-1} {x}_{0}^{k-n}} \Bigr) > \sum_i h_i \geq  -\frac{n}{\lambda_n}.
\end{align*}
Second, similarly, we also have
\begin{align*}
\frac{n}{\lambda_n^2} \geq \sum_i h_i^2  >  \frac{1}{\lambda_n^2}\Bigl( 1 -      {\sum_{k=1}^{n-2} c_k \mybinom[0.8]{n-1}{k-1} x_0^{k-n} }  \Bigr)^2.
\end{align*}
Third, for $i \neq j$, we have
\begin{align*}
h_{ij} &= \frac{  \sum_{k=0}^{n-2} c_k \sigma_k(\lambda_{; i, j})   }{  {\lambda}_1  \cdots   {\lambda}_n   {\lambda}_i   {\lambda}_j} = \frac{\sum_{k=0}^{n-2} c_k \sigma_k(\lambda_{; j})}{\lambda_1 \cdots \lambda_n \lambda_i \lambda_j} - \frac{\sum_{k=1}^{n-2} c_k \sigma_{k-1}(\lambda_{; i, j})} {\lambda_1 \cdots \lambda_n   \lambda_j}    \\
&= \frac{1}{\lambda_i \lambda_j} - \frac{\sum_{k=1}^{n-2} c_k \sigma_{k-1}(\lambda_{; j})}{\lambda_1 \cdots \lambda_n \lambda_i} - \frac{\sum_{k=1}^{n-2} c_k \sigma_{k-1}(\lambda_{; i, j})} {\lambda_1 \cdots \lambda_n   \lambda_j}.
\end{align*}In addition, we have $\lambda_1 \cdots \lambda_n/\lambda_j - \sum_{k=1}^{n-2}c_k \sigma_{k-1}(\lambda_{; j}) > 0$ and by Lemma~\hyperlink{L:2.7}{2.7}, we get
\begin{align*}
\frac{\lambda_1 \cdots \lambda_n}{ \lambda_i  } - \sum_{k=1}^{n-2} c_k \sigma_{k-1}(\lambda_{; i, j}) &> \sum_{k=1}^{n-2} c_k \sigma_{k-1} (\lambda_{; i}) - \sum_{k=1}^{n-2} c_k \sigma_{k-1}(\lambda_{; i, j})  
= \lambda_j \sum_{k=2}^{n-2} c_k \sigma_{k-2}(\lambda_{; i, j}) \geq 0.
\end{align*}By combining these, for $i \neq j$, we obtain
\begin{align*}
\label{eq:3.18}
\frac{1}{\lambda_i \lambda_j}  \geq h_{ij} = \frac{1}{\lambda_i \lambda_j} - \frac{\sum_{k=1}^{n-2} c_k \sigma_{k-1}(\lambda_{; j})}{\lambda_1 \cdots \lambda_n \lambda_i} - \frac{\sum_{k=1}^{n-2} c_k \sigma_{k-1}(\lambda_{; i, j})} {\lambda_1 \cdots \lambda_n   \lambda_j} \geq - \frac{1}{\lambda_i \lambda_j}. \tag{3.18}
\end{align*}Also, we always have the following estimate for $h_j$:
\begin{align*}
\label{eq:3.19}
0 > h_j = - \frac{\sum_{k=0}^{n-2} c_k \sigma_k(\lambda_{; j}) }{\lambda_1 \cdots \lambda_n \lambda_j}  = - \frac{1}{\lambda_j} +  \frac{\sum_{k=1}^{n-2} c_k \sigma_{k-1} (\lambda_{; j})}{\lambda_1 \cdots \lambda_n} \geq - \frac{1}{\lambda_j}. \tag{3.19}
\end{align*}
Hence, for $i \in \{1, \cdots, n\}$, by inequalities (\ref{eq:3.18}) and (\ref{eq:3.19}), we have the following lower bound:
\begin{align*}
\sum_j h_{ij} h_j = h_{ii} h_i + \sum_{j \neq i} h_{ij} h_j  = -2 \frac{h_i^2}{\lambda_i} + \sum_{j \neq i} h_{ij} h_j  \geq - \frac{2}{\lambda_i^3} - \sum_{j \neq i} \frac{1}{\lambda_i \lambda_j^2} \geq - \frac{n+1}{\lambda_i \lambda_n^2}.
\end{align*}On the other hand, we get the following upper bound:
\begin{align*}
\sum_j h_{ij} h_j = h_{ii} h_i + \sum_{j \neq i} h_{ij} h_j  = -2 \frac{h_i^2}{\lambda_i} + \sum_{j \neq i} h_{ij} h_j  \leq  \sum_{j \neq i} \frac{1}{\lambda_i \lambda_j^2} \leq  \frac{n-1}{\lambda_i \lambda_n^2}.
\end{align*}By combining these, we obtain the following upper and lower bound:
\begin{align*}
\label{eq:3.20}
\frac{n-1}{\lambda_i \lambda_n^2}> \sum_l h_{i l} h_l \geq  -\frac{n+1}{\lambda_i \lambda_n^2}. \tag{3.20}
\end{align*}
Last, by inequality (\ref{eq:3.20}), we have
\begin{align*}
\frac{n(n+1)}{\lambda_n^4} \geq \sum_i \frac{n+1}{\lambda_i^2 \lambda_n^2} \geq  \sum_{i, j} h_i h_{ij} h_j > - \sum_i \frac{n-1}{\lambda_i^2 \lambda_n^2} \geq - \frac{n(n-1)}{\lambda_n^4}.
\end{align*}This finishes the proof.
\end{proof}

Now, by taking the first and second derivatives of the equation $H_c(z, \Lambda) = 1$, we have the following Lemma. The proof should be straightforward; we apply Lemma~\hyperlink{L:2.8}{2.8}, Lemma~\hyperlink{L:2.9}{2.9}, Lemma~\hyperlink{L:2.11}{2.11}, and Lemma~\hyperlink{L:3.2}{3.2}. 
Or one can check the following reference \cite{lin2020, lin2023d} for more details.
\hypertarget{L:3.6}{\begin{flemma}
Let $H_c( z, {\Lambda} )= 1$, then we have
\begin{align*}
0  &= \sum_{i ,j}  \frac{\partial H (   {\Lambda}     ) }{\partial  {\Lambda}^j_i}    \frac{\partial  {\Lambda}^j_i }{\partial  {z}_k} + \frac{  \frac{\partial c_0}{\partial z_k} }{\sigma_n(\Lambda)}; \\
0 &=  \sum_{i ,j} \Bigl( \frac{\partial^2 H (   {\Lambda}     ) }{\partial \Lambda^j_i \partial \Lambda_r^s}    \frac{\partial  {\Lambda}_i^j }{\partial \bar{z}_k} \frac{\partial  {\Lambda}_r^s}{\partial z_k} +    \frac{\partial H (   {\Lambda}     )}{\partial \Lambda_i^j }    \frac{\partial^2  {\Lambda}_i^j }{\partial z_k \partial \bar{z}_k}  +     \frac{\partial}{\partial \Lambda^j_i} \Bigl(  \frac{  2 }{\sigma_n(\Lambda)} \Bigr) \Re \Bigl( \frac{\partial c_0}{\partial \bar{z}_k}  \frac{\partial \Lambda^j_i}{\partial {z}_k} \Bigr)   \Bigr) + \frac{   \frac{\partial^2 c_0}{\partial z_k \partial \bar{z}_k } }{\sigma_n (\Lambda)}.  
% &\kern2em +  \frac{\partial c_1}{\partial z_k} \frac{\partial}{\partial \Lambda^j_i} \frac{  \sigma_1(\Lambda) }{\sigma_3(\Lambda)} \frac{\partial \Lambda^j_i}{\partial \bar{z}_k}  +  \frac{\partial c_1}{\partial \bar{z}_k} \frac{\partial}{\partial \Lambda^j_i} \frac{  \sigma_1(\Lambda) }{\sigma_3(\Lambda)} \frac{\partial \Lambda^j_i}{\partial {z}_k} +   \frac{\partial c_0}{\partial z_k} \frac{\partial}{\partial \Lambda^j_i} \frac{  1 }{\sigma_3(\Lambda)} \frac{\partial \Lambda^j_i}{\partial \bar{z}_k} +   \frac{\partial c_0}{\partial \bar{z}_k} \frac{\partial}{\partial \Lambda^j_i} \frac{  1 }{\sigma_3(\Lambda)} \frac{\partial \Lambda^j_i}{\partial {z}_k}
\end{align*}In particular, at the maximum point $q \in M$ of $\tilde{U}$, we have
\begin{align*}
\label{eq:3.21}
0 &=  \sum_i h_i   (  X_u     )_{i \bar{i},  {k}} + \frac{   \frac{\partial c_0}{\partial z_k} }{\lambda_1 \cdots \lambda_n};   \tag{3.21}      \\
\label{eq:3.22}
0 &= \sum_{i, j}  h_{ij}   (  X_u     )_{i \bar{i}, \bar{k}}    (  X_u     )_{j \bar{j}, {k}} + \sum_{i \neq j}  \Bigl(   \frac{ 1 }{  \lambda_i \lambda_j}  -  \frac{\sum_{k=2}^{n-2} c_k \sigma_{k-2}(\lambda_{; i, j})}{\lambda_1 \cdots \lambda_n}   \Bigr)  | (X_u)_{j \bar{i},k}    |^2  \tag{3.22} \\
&\kern2em +  \sum_i h_i \Bigl(   ( X_u   )_{i \bar{i},k \bar{k}} -   {\lambda}_i \omega_{i \bar{i},k \bar{k}} \Bigr) 
 + \frac{   \frac{\partial^2 c_0}{\partial z_k \partial \bar{z}_k } }{ \lambda_1 \cdots \lambda_n } - \sum_i \frac{2 \Re \Bigl(  \frac{\partial c_0}{\partial \bar{z}_k}  (X_u)_{i \bar{i}, {k}} \Bigr)  }{\lambda_1 \cdots \lambda_n \lambda_i}.  
\end{align*}
\end{flemma}}

\begin{proof}
The first and second derivatives should be straightforward. At the maximum point, suppose the eigenvalues are pairwise distinct satisfying $\lambda_1 > \cdots > \lambda_n$. Since $\Lambda$ is a diagonal matrix, then
\begin{align*}
    0 =  \sum_i h_i   (  X_u     )_{i \bar{i},  {k}}+ \frac{  \frac{\partial c_0}{\partial z_k}  }{\lambda_1 \cdots \lambda_n }.
\end{align*}
 This is also true when the eigenvalues are not pairwise distinct. For the second derivative, if the eigenvalues at the maximum point $q$ are pairwise distinct, then 
\begin{align*}
&\kern-2em 0
= \sum_{i, j}  h_{ij}   (  X_u     )_{i \bar{i}, \bar{k}}    (  X_u     )_{j \bar{j}, {k}} + \sum_{i \neq j}  \Bigl(   \frac{ 1 }{  \lambda_i \lambda_j}  -  \frac{\sum_{k=2}^{n-2} c_k \sigma_{k-2}(\lambda_{; i, j})}{\lambda_1 \cdots \lambda_n}   \Bigr)  | (X_u)_{j \bar{i},k}    |^2  \\  
&\kern2em+  \sum_i h_i \Bigl(   ( X_u   )_{i \bar{i},k \bar{k}} -   {\lambda}_i \omega_{i \bar{i},k \bar{k}} \Bigr)  + \frac{  \frac{\partial^2 c_0}{\partial z_k \partial \bar{z}_k } }{ \lambda_1 \cdots \lambda_n }   - \frac{2}{\lambda_1 \cdots \lambda_n \lambda_i}   \Re \Bigl( \frac{\partial c_0}{\partial \bar{z}_k} (X_u)_{i \bar{i}, {k}} \Bigr).
\end{align*}This is also true when the eigenvalues are not pairwise distinct.
\end{proof}

As a quick consequence of Theorem~\hyperlink{T:2.1}{2.1}, we have the following general Hadamard's inequality.
\hypertarget{Cor:3.1}{\begin{fcor}[General Hadamard's inequality]
Let $c  \in \tilde{\mathscr{C}}_n$ and $X$ be a positive definite Hermitian matrix with eigenvalues $\{\lambda_1, \cdots, \lambda_n\}$ satisfying $\lambda_1 \cdots \lambda_n - \sum_{k=0}^{n-2} c_k \sigma_k(\lambda) = 0$. Then the diagonal entries $\{X_1, \cdots, X_n \}$ of $X$ satisfies 
\begin{align*}
X_1 \cdots X_n - \sum_{k=0}^{n-2} c_k \sigma_k(X) \geq 0.
\end{align*}
\end{fcor}}

\begin{proof}
By the Schur--Horn theorem, the diagonal entries $\{X_1, \cdots, X_n \}$ of $X$ is in the permutation polytope generated by the eigenvalues $\{ \lambda_1, \cdots, \lambda_n\}$ of $X$.
Since $c = (c_{n-2}, \cdots, c_1, c_0) \in \tilde{\mathscr{C}}_n$, by Theorem~\hyperlink{T:2.1}{2.1}, the level set $\lambda_1 \cdots \lambda_n - \sum_{k=0}^{n-2} c_k \sigma_k = 0$ is convex. By combining these, we get
$X_1 \cdots X_n - \sum_{k=0}^{n-2} c_k \sigma_k(X) \geq 0,$
which finishes the proof.
\end{proof}

To simplify estimates and compute asymptotic behavior, for the remainder of this section, we let $O_i$ be the Big $O$ notation that describes the limiting behavior when $\lambda_i$ approaches infinity. So $O_i(1)$ means the quantity will be bounded by a uniform constant if $\lambda_i$ is sufficiently large.

\hypertarget{L:3.7}{\begin{flemma}
Let $d \colon M^n \rightarrow \tilde{\mathscr{C}}_n$ be given by $d(z) = ( d_{n-2}, \cdots, d_1, d_0(z)   )$, $X$ be a $C$-subsolution to $d$, ${c}(z) = (c_{n-2}, \cdots, c_1, c_0(z)) \in \tilde{\mathcal{P}}^d$ with $c(z) \in S$ for any $z \in M$, and $u$ is a solution to $c \colon M^n \rightarrow \tilde{\mathscr{C}}_n$. Then there exists uniform constants $N  = N(d, X, S)> 0$ and $\kappa = \kappa(d, X, S) > 0$, which are independent of ${c}$, such that if $\lambda_1 > N$,  we have
$\sum_i h_i u_{i \bar{i}} \geq -\kappa \sum_i h_i$.
%\textcolor{red}{Maybe can be improved}.
\end{flemma}}
\begin{proof}
If $X$ is a $C$-subsolution to $d \colon M^n \rightarrow \tilde{\mathscr{C}}_{n}$ given by $d(z) =  ( d_{n-2}, \cdots, d_1, d_0(z)   )$, then since the $\Upsilon_1$-cone of $d$ is an open set, there exists a $\kappa > 0$ sufficiently small such that $X - 2 \kappa \omega$ is still in the $\Upsilon_1$-cone of $d$. In addition, by choosing $\epsilon > 0$ and $\delta > 0$ sufficiently small, for any $l \in \{1, \cdots, n-1\}$, we have
\begin{align*}
\label{eq:3.23.l}
(1-\delta) (X - 2 \kappa \omega)^{n - l} - \sum_{k = l}^{n-2}   d_k \mybinom[0.8]{n-l}{k-l} (X - 2 \kappa \omega)^{k-l} \wedge \omega^{n-k} > \epsilon \omega^{n- l }. \tag{3.23.$l$}
\end{align*}
In addition, for any $c \in \tilde{\mathcal{P}}^d$ with $c(z) \in S$ for any $z \in M$, the $\Upsilon_1$-cone of $d$ is contained in the $\Upsilon_1$-cone of $c$. Hence, for any $z \in M$ and $l \in \{1, \cdots, n-1\}$, by inequality (\ref{eq:3.23.l}), we get 
\begin{align*}
\label{eq:3.24.l}
(1-\delta)  (X - 2 \kappa \omega)^{n - l} - \sum_{k = l}^{n-2}  c_k \mybinom[0.8]{n-l}{k-l} (X - 2 \kappa \omega)^{k-l} \wedge \omega^{n-k} > \epsilon \omega^{n- l }. \tag{3.24.$l$}
\end{align*}
Note that $u_{ i \bar{i}} = \lambda_i - X_{i \bar{i}}$, so we can write
\begingroup
\allowdisplaybreaks
\begin{align*}
\label{eq:3.25}
\sum_i  h_i   ( u_{i \bar{i}}  + \kappa     ) &= \sum_i h_i  \bigl (  {\lambda}_i - X_{i \bar{i}}    + \kappa   \bigr ) = - \sum_i  \frac{\sum_{k=0}^{n-2} c_k \sigma_k (\lambda_{; i})  }{   {\lambda}_1   \cdots   {\lambda}_n     {\lambda}_i }     (   {\lambda}_i - X_{i \bar{i}}     + \kappa     ), \tag{3.25}
\end{align*}
\endgroup
We always assume that $\lambda_1 \geq \cdots \geq \lambda_n$ for convenience. There are $n$ cases to consider: 

$\bullet$ If $\lambda_{n}$ is less than $(X_{n \bar{n}} - \kappa)/n$: Then equation (\ref{eq:3.25}) becomes 
\begin{align*}
\sum_i  h_i   ( u_{i \bar{i}}  + \kappa     ) &\geq (X_{n \bar{n}} - \kappa)  \frac{\sum_{k=0}^{n-2} c_k \sigma_k (\lambda_{; n})  }{   {\lambda}_1   \cdots   {\lambda}_n     {\lambda}_n }   -  \sum_i  \frac{\sum_{k=0}^{n-2} c_k \sigma_k (\lambda_{; i})  }{   {\lambda}_1   \cdots   {\lambda}_n    }   \\
&\geq   (X_{n \bar{n}} - \kappa)  \frac{\sum_{k=0}^{n-2} c_k \sigma_k (\lambda_{; n})  }{   {\lambda}_1   \cdots   {\lambda}_n     {\lambda}_n }   -  n  \frac{\sum_{k=0}^{n-2} c_k \sigma_k (\lambda_{; n})  }{   {\lambda}_1   \cdots   {\lambda}_n    } \geq 0,
\end{align*}which finishes the proof.\smallskip

$\bullet$ If $\lambda_{n-1}$ is sufficiently large and $\lambda_n$ is uniformly bounded from below: Since $S$ is a compact subset, there exists a $N > 0$ such that for any $c \in \tilde{\mathcal{P}}^d$ with $c(z) \in S$ and $\lambda_{n-1} \geq N$, we have
\begin{align*}
\lambda_n = \frac{\sum_{k=0}^{n-2} c_k \sigma_k(\lambda_{; n}) }{ \lambda_1 \cdots \lambda_{n-1} - \sum_{k=1}^{n-2} c_k \sigma_{k-1}(\lambda_{; n})  } \leq \frac{C}{N}.
\end{align*} 
In this case, by choosing $N = N_1 > 0$ sufficiently large, we get a contradiction and finish the proof.\smallskip

Now, we argue the remaining cases inductively, for $l$ from $1$ to $n-2$:\smallskip

$\bullet$ If $\lambda_{n-l-1}$ is sufficiently large, $\lambda_{n-l}$ is uniformly bounded from above by $N_l$, and $\lambda_{n}$ is uniformly bounded from below: By enlarging $N > N_l$ if necessary, for any $c \in \tilde{\mathcal{P}}^d$ with $c(z) \in S$ and $\lambda_{n-l-1} \geq N$, we can rewrite the original equation $\lambda_1 \cdots \lambda_n - \sum_{k=0}^{n-2} c_k \sigma_k(\lambda) = 0$ as
\begin{align*}
\lambda_1 \cdots \lambda_n &= \sum_{k=0}^{n-2} c_k \sigma_k(\lambda_{; n-l, \cdots, n}) + (\lambda_{n-l} + \cdots + \lambda_n) \sum_{k=1}^{n-2} c_k \sigma_{k-1}(\lambda_{; n-l, \cdots, n}) \\
&\kern2em + \cdots  + \lambda_{n-l} \cdots \lambda_n \sum_{k=l}^{n-2} c_k \sigma_{k-l} (\lambda_{; n-l, \cdots, n}).
\end{align*}This implies that ${C}/{N}  >  \lambda_{n-l} \cdots \lambda_n - c_{n-l-1} - c_{n-l} (\lambda_{n-l} + \cdots + \lambda_n) - \cdots - c_{n-2}\sigma_{l-1}(\lambda_{; 1, \cdots, n-l-1})$. By choosing $N = N_{l +1} > 0$ sufficiently large, we get 
\begin{align*}
\label{eq:3.26}
\epsilon > \frac{\delta}{2}\lambda_{n-l} \cdots \lambda_n >  \lambda_{n-l} \cdots \lambda_n - \sum_{k=n-l-1}^{n-2} c_k \sigma_{k-n+l+1}(\lambda_{; 1, \cdots, n-l-1})  > 0. \tag{3.26}
\end{align*}
In this case, equation (\ref{eq:3.25}) becomes 
\begingroup
\allowdisplaybreaks
\begin{align*}
\label{eq:3.27}
\sum_i  h_i   ( u_{i \bar{i}}  + \kappa     ) &= - \sum_i  \frac{\sum_{k=0}^{n-2} c_k \sigma_k (\lambda_{; i})  }{   {\lambda}_1   \cdots   {\lambda}_n     {\lambda}_i }     (   {\lambda}_i - X_{i \bar{i}}     + \kappa     ) \tag{3.27} \\
&=  \sum_{j = n-l}^n (X_{ j \overline{j}} - 2 \kappa)   \frac{\sum_{k=0}^{n-2} c_k \sigma_k (\lambda_{; j})  }{   {\lambda}_1   \cdots   {\lambda}_n     {\lambda}_{j} }   - \sum_i   \frac{\sum_{k=0}^{n-2} c_k \sigma_k (\lambda_{; i})  }{   {\lambda}_1   \cdots   {\lambda}_n  }   - \kappa \sum_{i} h_i  \\
&\geq \sum_{j = n-l}^n (X_{ j \overline{j}} - 2 \kappa - \lambda_j )   \frac{     \sum_{k=n-l-1}^{n-2} c_k \sigma_{k-n+l+1} (\lambda_{; 1, \cdots, n-l-1, j})  }{   {\lambda}_{n-l}   \cdots   {\lambda}_n     {\lambda}_{j} }   \\  
&\kern2em - \kappa \sum_{i} h_i  +  \lambda_{n-l-1}^{-1} O_{n-l-1}(1).   
\end{align*}
\endgroup
By Corollary~\hyperlink{Cor:3.1}{3.1} and inequality (\ref{eq:3.24.l}), we have 
\begin{align*}
(1- \delta)  (X_{n-l \  \overline{n-l}} - 2\kappa)  \cdots (X_{n \bar{n}} - 2 \kappa)    -  \sum_{k=n-l-1}^{n-2}   c_k \sigma_{k-n+l+1}\bigl( (X-2\kappa)_{; 1, \cdots, n-l-1} \bigr)  > \epsilon.
\end{align*}
In addition, for any $(\lambda_{n-l}, \cdots, \lambda_n)$ in this case, there exists a corresponding $\tilde{\delta} \in (0, \delta/2)$ such that
\begin{align*}
\tilde{\delta} \lambda_{n-l} \cdots \lambda_n = \lambda_{n-l} \cdots \lambda_n - \sum_{k=n-l-1}^{n-2} c_k \sigma_{k-n+l+1}(\lambda_{; 1, \cdots, n-l-1}),
\end{align*}which implies that $(\lambda_{n-l}, \cdots, \lambda_n)$ is on the following set
\begin{align*}
\Bigl\{ (1 - \tilde{\delta} ) \lambda_{n-l} \cdots \lambda_n  - \sum_{k=n-l-1}^{n-2} c_k \sigma_{k-n+l+1}(\lambda_{; 1, \cdots, n-l-1}) = 0\Bigr\}.
\end{align*}For this case, we can see that $(c_{n-l-1}, \cdots, c_{n-2}) \neq 0$, otherwise we get a contradiction. By Lemma~\hyperlink{L:2.5}{2.5} and Theorem~\hyperlink{T:2.3}{2.3}, we have
\begin{align*}
  \bigl(X_{n-l \  \overline{n-l}} - 2 \kappa,  \cdots, X_{n \bar{n}} - 2 \kappa \bigr)  \in \Bigl \{ (1 - \tilde{\delta} ) \lambda_{n-l} \cdots \lambda_n  - \sum_{k=n-l-1}^{n-2} c_k \sigma_{k-n+l+1}(\lambda_{; 1, \cdots, n-l-1}) > 0 \Bigr\}.
\end{align*}By Theorem~\hyperlink{T:2.1}{2.1}, $\{ (1 - \tilde{\delta} ) \lambda_{n-l} \cdots \lambda_n  - \sum_{k=n-l-1}^{n-2} c_k \sigma_{k-n+l+1}(\lambda_{; 1, \cdots, n-l-1}) = 0\}$ is a convex set. Hence, by the supporting hyperplane theorem, we have
\begin{align*}
\label{eq:3.28}
\sum_{j = n-l}^n     ( X_{j \bar{j}} - 2 \kappa - \lambda_j  )  \frac{     \sum_{k=n-l-1}^{n-2} c_k \sigma_{k-n+l+1} (\lambda_{; 1, \cdots, n-l-1, j})  }{   {\lambda}_{n-l}   \cdots   {\lambda}_n     {\lambda}_{j} }  \geq 0. \tag{3.28}
\end{align*}
By combining inequalities (\ref{eq:3.27}) and (\ref{eq:3.28}) and by Lemma~\hyperlink{L:3.5}{3.5}, we get
\begin{align*}
&\kern-2em \sum_i  h_i   ( u_{i \bar{i}}  + \kappa     ) \\
&\geq   -\kappa \sum_i h_i   +  \lambda_{n-l-1}^{-1} O_{n-l-1}(1)  > \frac{\kappa}{\lambda_n} \Bigl( 1 -    {\sum_{k=1}^{n-2} c_k \mybinom[0.8]{n-1}{k-1} {x}_{0}^{k-n}} \Bigr)   +  \lambda_{n-l-1}^{-1} O_{n-l-1}(1) \\
&\geq \frac{\kappa}{x_0} \Bigl( 1 -  \sup_{c \in S}  {\sum_{k=1}^{n-2} c_k \mybinom[0.8]{n-1}{k-1} {x}_{0}^{k-n}}(c) \Bigr) +  \lambda_{n-l-1}^{-1} O_{n-l-1}(1) \\
&\geq  \frac{\kappa}{\inf_{c \in S} x_0(c)} \Bigl( 1 -  \sup_{c \in S}  {\sum_{k=1}^{n-2} c_k \mybinom[0.8]{n-1}{k-1} {x}_{0}^{k-n}(c) } \Bigr) +  \lambda_{n-l-1}^{-1} O_{n-l-1}(1) \geq 0.
\end{align*}provided that $N_{l+1}$ is sufficiently large. Notice that by Lemma~\hyperlink{L:2.2}{2.2}, $x_0(c)$ is a positive continuous function on $\tilde{\mathscr{C}}_n$. Since $S$ is a compact subset of $\tilde{\mathscr{C}}_n$, $x_0(c)$ attains its minimum on $S$ with minimum value greater than $0$. Similarly, for any $c \in \tilde{\mathscr{C}}_n$, $c$ corresponds to a strictly right-Noetherian polynomial, which implies that $x_0(c) > x_1(c)$ and hence $\sum_{k=1}^{n-2} c_k \binom{n-1}{k-1} x_0^{k-n}(c) < 1$. Thus, $ {\sum_{k=1}^{n-2} c_k \binom{n-1}{k-1} {x}_{0}^{k-n}(c) }$ is a continuous function on $\tilde{\mathscr{C}}_n$ with value less than $1$. Since $S$ is a compact subset of $\tilde{\mathscr{C}}_n$, $ {\sum_{k=1}^{n-2} c_k \binom{n-1}{k-1} {x}_{0}^{k-n}(c) }$ attains its maximum on $S$ with maximum value less than $1$. This finishes the proof. 
\end{proof}

Now, we let $C$ be a positive constant depending only on the stated data, but which may change from line to line. With all the previous work, we can finish the proof of the following $C^2$ estimate.

\hypertarget{T:3.1}{\begin{fthm}
Suppose $d \colon M^n \rightarrow  \tilde{\mathscr{C}}_n$ is a strictly $\Upsilon$-stable general inverse $\sigma_k$ equation with $d_1, \cdots, d_{n-2}$ constants and range in $S$ and $X$ is a $C$-subsolution to $d$. For any $c \in \tilde{\mathcal{P}}^d$ with $c_1, \cdots, c_{n-2}$ constants and range in $S$, if $u \colon M   \rightarrow \mathbb{R}$ is a smooth function solving $c \colon M^n \rightarrow \tilde{\mathscr{C}}_n$, then there exists a constant $C$ such that 
\begin{align*}
|\partial \bar{\partial} u | \leq C \bigl ( 1 + \sup_M \bigl|\nabla u\bigr|^2  \bigr).
\end{align*}Here, $C = C(M, X, S, d, \omega, \osc_M u, \|c_0\|_{C^2})$ is a constant  and $\nabla$ is the Levi-Civita connection with respect to $\omega$.
\end{fthm}}

\begin{proof}
We use the maximum principle to prove this statement. For any $c \in \tilde{\mathcal{P}}^d$ with $c_1, \cdots, c_{n-2}$ constants and range in $S$, we can define the elliptic operator $\mathcal{L}_c$ in equation (\ref{eq:3.6}). First, by applying the operator $\mathcal{L}_c$ to $G(\tilde{\Lambda})$, at the maximum point $q$, we obtain
\begingroup
\allowdisplaybreaks
\begin{align*}
\label{eq:3.29}
\kern-2em\mathcal{L}_c  \bigl (  G( \tilde{\Lambda} )   \bigr )  \tag{3.29}
&= - \sum_{i, j, k} h_k   g_{ij}  \frac{\partial \tilde{\Lambda}_i^i}{\partial z_k} \frac{\partial \tilde{\Lambda}_j^j}{\partial \bar{z}_k} - \sum_k h_k   \sum_{i \neq j} \frac{g_i - g_j}{ \tilde{\lambda}_i - \tilde{\lambda}_j}   \frac{\partial \tilde{\Lambda}_j^i}{\partial z_k} \frac{\partial \tilde{\Lambda}_i^j}{\partial \bar{z}_k} - \sum_{i,k} h_k   g_i   \frac{\partial^2 \tilde{\Lambda}_i^i}{\partial z_k \partial \bar{z}_k}    \\
&=  \sum_k h_k \frac{ \bigl|    (  X_u     )_{1 \bar{1},k }    \bigr |^2}{(1 + \tilde{\lambda}_1)^2}   + \sum_k h_k  \frac{ {\lambda}_1}{1 + \tilde{\lambda}_1} \omega_{1\bar{1},k \bar{k}}  -  \sum_k h_k \frac{ ( X_u   )_{1 \bar{1},k \bar{k}}}{1+ \tilde{\lambda}_1}   \\
&\kern2em  - \sum_{k}  h_k \sum_{j \neq 1}  \frac{ \bigl |    (  X_u    )_{j \bar{1},k}     \bigr |^2   +      \bigl |    (  X_u    )_{1 \bar{j},k}    \bigr  |^2}{(1 + \tilde{\lambda}_1)( \tilde{\lambda}_1 - \tilde{\lambda}_j)}    \\
&\geq   \sum_i h_i \frac{ \bigl |    (  X_u     )_{1 \bar{1}, i}    \bigr |^2}{(1 + \tilde{\lambda}_1)^2}    - \sum_i h_i \frac{ (  X_u   )_{1 \bar{1}, i \bar{i}}}{1+ \tilde{\lambda}_1}     -     \sum_{j \neq 1}  \frac{   h_j  \bigl |    (  X_u    )_{1 \bar{j}, j}    \bigr  |^2}{(1 + \tilde{\lambda}_1)( \tilde{\lambda}_1 - \tilde{\lambda}_j)} +  C \sum_i h_i.   
% &\kern2em - \sum_{i}  h_i \sum_{j \neq 1}  \frac{1}{(1 + \tilde{\lambda}_1)( \tilde{\lambda}_1 - \tilde{\lambda}_j)}  \Bigr(  \bigl |    (  X_u    )_{j \bar{1}, i}    \bigr  |^2 +      \bigl |    (  X_u    )_{1 \bar{j}, i}    \bigr  |^2      \Bigr ) . 
\end{align*}
\endgroup
Second, by equation (\ref{eq:3.22}) and inequality (\ref{eq:3.29}), we have
\begingroup
\allowdisplaybreaks
\begin{align*}
\label{eq:3.30}
  0  &= \sum_{i, j}  h_{ij}   (  X_u     )_{i \bar{i}, \bar{k}}    (  X_u     )_{j \bar{j}, {k}} + \sum_{i \neq j} \Bigl(  \frac{ 1 }{ \lambda_i \lambda_j} - \frac{\sum_{k=2}^{n-2} c_k \sigma_{k-2} (\lambda_{; i, j}) }{\lambda_1 \cdots  \lambda_n  } \Bigr) | (X_u)_{j \bar{i},k}    |^2  \tag{3.30} \\
&\kern2em +  \sum_i h_i \bigl(   ( X_u   )_{i \bar{i},k \bar{k}} -   {\lambda}_i \omega_{i \bar{i},k \bar{k}} \bigr)   + \frac{   \frac{\partial^2 c_0}{\partial z_k \partial \bar{z}_k } }{ \lambda_1 \cdots \lambda_n } - \sum_i \frac{2 \Re \bigl(  \frac{\partial c_0}{\partial \bar{z}_k}  (X_u)_{i \bar{i}, {k}} \bigr)  }{\lambda_1  \cdots  \lambda_n \lambda_i} \\
&= \sum_{i, j}  h_{ij}  \Bigl(  (  X_u     )_{i \bar{i}, \bar{k}}  - \frac{\sum_l h_l (X_u)_{l \bar{l}, \bar{k} } }{\sum_{l} h_l^2 } h_i     \Bigr) \Bigl( (  X_u     )_{j \bar{j}, {k}} - \frac{\sum_l h_l (X_u)_{l \bar{l}, {k} } }{\sum_{l} h_l^2 } h_j     \Bigr)  \\
&\kern2em + \frac{2 \sum_{i, j, l} h_j h_{ij} h_l \Re \bigl(  (X_u)_{l \bar{l}, {k}} (X_u)_{i \bar{i}, \bar{k} }     \bigr) }{\sum_i h_i^2}  - \sum_{i, j} h_i h_j h_{ij} \frac{ | \sum_l h_l (X_u)_{l \bar{l}, k }    |^2 }{ (\sum_i h_i^2)^2} \\ 
&\kern2em + \sum_{i \neq j} \Bigl(  \frac{ 1 }{ \lambda_i \lambda_j} - \frac{\sum_{k=2}^{n-2} c_k \sigma_{k-2}(\lambda_{; i, j})}{\lambda_1  \cdots  \lambda_n} \Bigr) | (X_u)_{j \bar{i},k}    |^2    +  \sum_i h_i \Bigl(   ( X_u   )_{i \bar{i},k \bar{k}} -   {\lambda}_i \omega_{i \bar{i},k \bar{k}} \Bigr)   \\ 
&\kern2em + \frac{   \frac{\partial^2 c_0}{\partial z_k \partial \bar{z}_k } }{ \lambda_1 \cdots \lambda_n } - \sum_i \frac{2 \Re \bigl(  \frac{\partial c_0}{\partial \bar{z}_k}  (X_u)_{i \bar{i}, {k}} \bigr)  }{\lambda_1 \cdots \lambda_n \lambda_i} \\
&= \sum_{i, j}  h_{ij}  \Bigl(  (  X_u     )_{i \bar{i}, \bar{k}}  - \frac{\sum_l h_l (X_u)_{l \bar{l}, \bar{k} } }{\sum_{l} h_l^2 } h_i     \Bigr) \Bigl( (  X_u     )_{j \bar{j}, {k}} - \frac{\sum_l h_l (X_u)_{l \bar{l}, {k} } }{\sum_{l} h_l^2 } h_j     \Bigr)  \\
&\kern2em - \frac{2 \sum_{i, j} h_j h_{ij}  \Re \bigl(   \frac{\partial c_0}{\partial \bar{z}_k}     (X_u)_{i \bar{i}, {k} }     \bigr) }{ \lambda_1 \cdots \lambda_n \sum_i h_i^2}   -  \frac{ \sum_{i, j} h_i h_j h_{ij} |  \frac{\partial c_0}{\partial z_k}|^2 }{  \lambda_1^2  \cdots \lambda_n^2 ( \sum_i h_i^2)^2}     \\  
&\kern2em + \sum_{i \neq j} \Bigl(  \frac{ 1 }{ \lambda_i \lambda_j} - \frac{ \sum_{k=2}^{n-2} c_k \sigma_{k-2}(\lambda_{; i, j})}{\lambda_1 \cdots \lambda_n} \Bigr) | (X_u)_{j \bar{i},k}    |^2   +  \sum_i h_i \Bigl(   ( X_u   )_{i \bar{i},k \bar{k}} -   {\lambda}_i \omega_{i \bar{i},k \bar{k}} \Bigr)   \\ 
&\kern2em + \frac{   \frac{\partial^2 c_0}{\partial z_k \partial \bar{z}_k } }{ \lambda_1 \cdots \lambda_n } - \sum_i \frac{2 \Re \bigl(  \frac{\partial c_0}{\partial \bar{z}_k}  (X_u)_{i \bar{i}, {k}} \bigr)  }{\lambda_1 \cdots \lambda_n \lambda_i} \\
&\geq  - \frac{2 \sum_{i, j} h_j h_{ij}  \Re \bigl(    \frac{\partial c_0}{\partial \bar{z}_k}    (X_u)_{i \bar{i}, {k} }     \bigr)  }{ \lambda_1  \cdots \lambda_n \sum_i h_i^2} - \sum_i \frac{2 \Re \bigl(   \frac{\partial c_0}{\partial \bar{z}_k}  (X_u)_{i \bar{i}, {k}} \bigr)  }{\lambda_1  \cdots  \lambda_n \lambda_i}  - \frac{ \sum_{i, j} h_i h_j h_{ij}  |  \frac{\partial c_0}{\partial z_k}|^2 }{  \lambda_1^2 \cdots  \lambda_n^2 ( \sum_i h_i^2)^2}      \\ 
&\kern2em +\sum_{i \neq j}  \Bigl(  \frac{ 1 }{ \lambda_i \lambda_j} - \frac{\sum_{k=2}^{n-2} c_k \sigma_{k-2}(\lambda_{; i, j}) }{\lambda_1 \cdots \lambda_n} \Bigr) | (X_u)_{j \bar{i},k}    |^2  +  \sum_i  h_i  ( X_u   )_{i \bar{i},k \bar{k}} + C \sum_i h_i \lambda_i  \\
&\kern2em + \frac{   \frac{\partial^2 c_0}{\partial z_k \partial \bar{z}_k } }{ \lambda_1 \cdots \lambda_n }, 
\end{align*}
\endgroup
where the inequality on the last line is due to the convexity of $h =1$ and $|\omega_{i \bar{i}, k \bar{k}}| \leq C$. Since by Theorem~\hyperlink{T:2.1}{2.1}, the solution set $\{ \lambda_1 \cdots \lambda_n = \sum_{k=0}^{n-2} c_k \sigma_k(\lambda)  \}$ is convex. In addition, we have  
\begin{align*}
(  X_u     )_{j \bar{j}, {k}}  - \frac{\sum_l h_l (X_u)_{l \bar{l}, \bar{k} } }{\sum_{l} h_l^2 } h_j = (  X_u     )_{j \bar{j}, {k}}  - \frac{   \frac{\partial c_0}{\partial {z}_k} h_j }{\lambda_1 \cdots \lambda_n \sum_{l} h_l^2 } 
\end{align*}is a tangent vector on the solution set. Hence, we obtain that 
\begin{align*}
    \sum_{i, j}  h_{ij}  \Bigl(  (  X_u     )_{i \bar{i}, \bar{k}}  - \frac{\sum_l h_l (X_u)_{l \bar{l}, \bar{k} } }{\sum_{l} h_l^2 } h_i     \Bigr) \Bigl( (  X_u     )_{j \bar{j}, {k}} - \frac{\sum_l h_l (X_u)_{l \bar{l}, {k} } }{\sum_{l} h_l^2 } h_j     \Bigr)  \geq 0.
\end{align*} Now, by setting $k = 1$, inequality (\ref{eq:3.30}) gives
\begingroup
\allowdisplaybreaks
\begin{align*}
\label{eq:3.31}
&\kern-2em -\sum_i h_i   (  X_u     )_{1 \bar{1},i \bar{i}} \tag{3.31} \\
&= - \sum_i h_i (X_u  )_{i  \bar{i},1 \bar{1}}  + \sum_i h_i \bigl (    (X_u  )_{i  \bar{i},1 \bar{1}} -   ( X_u   )_{1 \bar{1}, i \bar{i}}  \bigr )   \\
&= - \sum_i h_i (X_u  )_{i  \bar{i},1 \bar{1}}  + \sum_i h_i \bigl (    (X   )_{i  \bar{i},1 \bar{1}} -   ( X    )_{1 \bar{1}, i \bar{i}}  \bigr ) \\
&\geq - \frac{2 \sum_{i, j} h_j h_{ij}  \Re \bigl(   \frac{\partial c_0}{\partial \bar{z}_1}     (X_u)_{i \bar{i}, {1} }     \bigr)}{ \lambda_1 \cdots \lambda_n \sum_i h_i^2} - \sum_i \frac{2 \Re \bigl(  \frac{\partial c_0}{\partial \bar{z}_1}  (X_u)_{i \bar{i}, {1}} \bigr)  }{\lambda_1 \cdots \lambda_n \lambda_i} -  \frac{ \sum_{i, j} h_i h_j h_{ij} |  \frac{\partial c_0}{\partial z_1}|^2 }{  \lambda_1^2 \cdots \lambda_n^2 ( \sum_i h_i^2)^2}   \\
&\kern2em +\sum_{i \neq j}  \Bigl(  \frac{ 1 }{ \lambda_i \lambda_j} - \frac{\sum_{k=2}^{n-2} c_k \sigma_{k-2}(\lambda_{; i, j})}{\lambda_1 \cdots \lambda_n} \Bigr) | (X_u)_{j \bar{i},k}    |^2   + \frac{   \frac{\partial^2 c_0}{\partial z_1 \partial \bar{z}_1 } }{ \lambda_1 \cdots \lambda_n }  + C \sum_i h_i   (1 + \lambda_i).  
% &\geq  \sum_{1 \neq j}  \frac{ 1 }{  \lambda_i \lambda_j}   | (X_u)_{j \bar{1},1}    |^2 +  C \sum_i h_i (1 + \lambda_i).
% &\geq  \sum_{j \neq 1}   \frac{c_1 \sigma_1(\lambda_{;i, j}) +c_0 }{\lambda_1 \lambda_2 \lambda_3 \lambda_i \lambda_j}   | (X_u)_{j \bar{1},1}    |^2 +  C \sum_i h_i (1 + \lambda_i)   +    C.
\end{align*}
\endgroup
Third, by combining Lemma~\hyperlink{L:3.1}{3.1}, Lemma~\hyperlink{L:3.7}{3.7}, and inequalities (\ref{eq:3.29}), (\ref{eq:3.30}), and (\ref{eq:3.31}), at the maximum point $q$, if $\lambda_1$ is sufficiently large, then we have
\begingroup
\allowdisplaybreaks
\begin{align*}
\label{eq:3.32}
\mathcal{L}_c   ( \tilde{U}   )  
&\geq   A   \sum_i h_i u_{i \bar{i}}  +  C \sum_i h_i  + \sum_i h_i \frac{\bigl |    (  X_u     )_{1 \bar{1},i}    \bigr |^2}{(1 + \tilde{\lambda}_1)^2}  -     \sum_{j \neq 1}  \frac{   h_j  \bigl |    (  X_u    )_{1 \bar{j}, j}    \bigr  |^2}{(1 + \tilde{\lambda}_1)( \tilde{\lambda}_1 - \tilde{\lambda}_j)}   \tag{3.32}      \\
&\kern2em  - \sum_i h_i  \frac{(  X_u   )_{1 \bar{1}, i\bar{i}}}{1+ \tilde{\lambda}_1}  \\
&\geq   A \sum_i h_i u_{i \bar{i}} +  C \sum_i h_i  + \sum_i  h_i \frac{  \bigl |    (  X_u     )_{1 \bar{1},i}    \bigr |^2}{(1 + \tilde{\lambda}_1)^2} -     \sum_{j \neq 1}  \frac{   h_j  \bigl |    (  X_u    )_{1 \bar{j}, j}    \bigr  |^2}{(1 + \tilde{\lambda}_1)( \tilde{\lambda}_1 - \tilde{\lambda}_j)}    \\
&\kern2em  +    \sum_{j \neq i} \Bigl(  \frac{ 1 }{ \lambda_i \lambda_j}  - \frac{\sum_{k=2}^{n-2} c_k \sigma_{k-2}(\lambda_{; i, j})}{\lambda_1 \cdots \lambda_n} \Bigr) \frac{| (X_u)_{j \bar{i}, 1}    |^2}{1+ \tilde{\lambda}_1}     \\  
&\kern2em   - \frac{2 \sum_{i, j} h_j h_{ij}  \Re \bigl(   \frac{\partial c_0}{\partial \bar{z}_1}    (X_u)_{i \bar{i}, {1} }     \bigr)}{ (1 + \tilde{\lambda}_1)  \lambda_1 \cdots \lambda_n \sum_i h_i^2} - \sum_i \frac{2 \Re \bigl(   \frac{\partial c_0}{\partial \bar{z}_1}  (X_u)_{i \bar{i}, {1}} \bigr)  }{(1 + \tilde{\lambda}_1)\lambda_1 \cdots \lambda_n \lambda_i}  \\
&\kern2em - \frac{  \sum_{i, j} h_i h_j h_{ij}  |  \frac{\partial c_0}{\partial z_1}|^2 }{  (1 + \tilde{\lambda}_1) \lambda_1^2 \cdots \lambda_n^2 ( \sum_i h_i^2)^2}    + \frac{   \frac{\partial^2 c_0}{\partial z_1 \partial \bar{z}_1 } }{ (1 + \tilde{\lambda}_1) \lambda_1 \cdots \lambda_n }    \\
&\geq  \bigl( C-A \kappa \bigr)  \sum_i h_i  +  \sum_i h_i \frac{ \bigl |    (  X_u     )_{1 \bar{1},i}    \bigr |^2}{(1 + \tilde{\lambda}_1)^2}   -     \sum_{j \neq 1}  \frac{   h_j  \bigl |    (  X_u    )_{1 \bar{j}, j}    \bigr  |^2}{(1 + \tilde{\lambda}_1)( \tilde{\lambda}_1 - \tilde{\lambda}_j)}      \\  
&\kern2em   +  \sum_{j \neq 1}  \Bigl(  \frac{ 1 }{ \lambda_1 \lambda_j} - \frac{\sum_{k=2}^{n-2} c_k \sigma_{k-2} (\lambda_{; 1, j}) }{\lambda_1 \cdots  \lambda_n } \Bigr) \frac{| (X_u)_{j \bar{1},1}    |^2}{1+ \tilde{\lambda}_1}  \\  
&\kern2em  - \frac{2 \sum_{i, j} h_j h_{ij}  \Re \bigl(    \frac{\partial c_0}{\partial \bar{z}_1}     (X_u)_{i \bar{i}, {1} }     \bigr)}{ (1 + \tilde{\lambda}_1)  \lambda_1 \cdots \lambda_n \sum_i h_i^2} - \sum_i \frac{2 \Re \bigl(   \frac{\partial c_0}{\partial \bar{z}_1}  (X_u)_{i \bar{i}, {1}} \bigr)  }{(1 + \tilde{\lambda}_1)\lambda_1  \cdots   \lambda_n \lambda_i}   \\
&\kern2em - \frac{  \sum_{i, j} h_i h_j h_{ij}  |  \frac{\partial c_0}{\partial z_1}|^2 }{  (1 + \tilde{\lambda}_1) \lambda_1^2 \cdots  \lambda_n^2 ( \sum_i h_i^2)^2}    + \frac{   \frac{\partial^2 c_0}{\partial z_1 \partial \bar{z}_1 } }{ (1 + \tilde{\lambda}_1) \lambda_1 \cdots  \lambda_n }  \\
% &\kern2em - \sum_{i}  h_i \sum_{j \neq 1}  \frac{ \bigl |    (  X_u    )_{j \bar{1}, i}    \bigr  |^2 +      \bigl |    (  X_u    )_{1 \bar{j}, i}    \bigr  |^2 }{(1 + \tilde{\lambda}_1)( \tilde{\lambda}_1 - \tilde{\lambda}_j)}        \\
&=    \bigl( C-A \kappa \bigr)  \sum_i h_i  +  h_1 \frac{\bigl |    (  X_u     )_{1 \bar{1}, 1}    \bigr |^2}{(1 + \tilde{\lambda}_1)^2}  - \frac{2 \sum_{i} h_i h_{1i}  \Re \bigl(   \frac{\partial c_0}{\partial \bar{z}_1}    (X_u)_{1 \bar{1}, {1} }     \bigr)}{ (1 + \tilde{\lambda}_1)  \lambda_1 \cdots  \lambda_n \sum_i h_i^2}  \\
&\kern2em  -  \frac{2 \Re \bigl(   \frac{\partial c_0}{\partial \bar{z}_1}  (X_u)_{1 \bar{1}, {1}} \bigr)  }{(1 + \tilde{\lambda}_1) \lambda_1 \cdots \lambda_n \lambda_1 }   +  \sum_{j \neq 1} h_j \frac{  \bigl |    (  X_u     )_{1 \bar{1}, j}    \bigr |^2}{(1 + \tilde{\lambda}_1)^2}    \\  
&\kern2em+  \sum_{j \neq 1}     \frac{ \lambda_2 \cdots \lambda_n - \lambda_j \sum_{k=2}^{n-2} c_k \sigma_{k-2}(\lambda_{; 1, j})}{\lambda_1 \cdots  \lambda_n \lambda_j }  \frac{| (X_u)_{j \bar{1},1}    |^2}{1+ \tilde{\lambda}_1}  \\
&\kern2em -  \sum_{j \neq 1}   \frac{   h_j   \bigl |    (  X_u    )_{1 \bar{j}, j}    \bigr  |^2}{(1 + \tilde{\lambda}_1)( \tilde{\lambda}_1 - \tilde{\lambda}_j)}     - \frac{2 \sum_{j \neq 1} \sum_i h_i h_{ij}  \Re \bigl(   \frac{\partial c_0}{\partial \bar{z}_1}   (X_u)_{j \bar{j}, {1} }     \bigr)}{ (1 + \tilde{\lambda}_1)  \lambda_1 \cdots \lambda_n \sum_i h_i^2}  \\
&\kern2em   - \sum_{j \neq 1} \frac{2 \Re \bigl(   \frac{\partial c_0}{\partial \bar{z}_1}  (X_u)_{j \bar{j}, {1}} \bigr)  }{(1 + \tilde{\lambda}_1)\lambda_1 \cdots \lambda_n \lambda_j}  - \frac{  \sum_{i, j} h_i h_j h_{ij}  |  \frac{\partial c_0}{\partial z_1}|^2 }{  (1 + \tilde{\lambda}_1) \lambda_1^2 \cdots \lambda_n^2 ( \sum_i h_i^2)^2}    + \frac{  \frac{\partial^2 c_0}{\partial z_1 \partial \bar{z}_1 } }{ (1 + \tilde{\lambda}_1) \lambda_1 \cdots \lambda_n }.
\end{align*}
\endgroup
We can simplify some terms in inequality (\ref{eq:3.32}). We denote $T_{j \bar{a}} \coloneqq  (X_u)_{j \bar{a}, 1}- (X_u)_{1 \bar{a}, j} = X_{j \bar{a}, 1}- X_{1 \bar{a}, j}$, then for $j \neq 1$, by Lemma~\hyperlink{L:2.7}{2.7}, we have: 
\begingroup
\allowdisplaybreaks
\begin{align*}
\label{eq:3.33}
&\kern-1em   h_j  \frac{ \bigl |    (  X_u     )_{1 \bar{1},j}    \bigr |^2}{(1 + \tilde{\lambda}_1)^2}          +     \frac{ \lambda_2 \cdots \lambda_n - \lambda_j \sum_{k=2}^{n-2} c_k \sigma_{k-2} (\lambda_{; 1, j}) }{\lambda_1 \cdots \lambda_n \lambda_j }  \frac{| (X_u)_{j \bar{1},1}    |^2}{1+ \tilde{\lambda}_1}  \tag{3.33}  \\
&=    h_j \frac{ \bigl |    (  X_u     )_{j \bar{1},1}   - T_{j \bar{1}} \bigr |^2}{(1 + \tilde{\lambda}_1)^2}   +    \frac{ \lambda_2 \cdots \lambda_n  - \lambda_j \sum_{k=2}^{n-2} c_k \sigma_{k-2} (\lambda_{; 1, j}) }{\lambda_1 \cdots \lambda_n \lambda_j}  \frac{| (X_u)_{j \bar{1},1}    |^2}{1+ \tilde{\lambda}_1}       \\
&\geq  2  h_j \frac{ \bigl |    (  X_u     )_{j \bar{1},1}    \bigr |^2}{(1 + \tilde{\lambda}_1)^2} +  2  h_j \frac{ \bigl |    T_{j\bar{1}} \bigr |^2}{(1 + \tilde{\lambda}_1)^2} +      \frac{ \lambda_2 \cdots \lambda_n  - \lambda_j \sum_{k=2}^{n-2} c_k \sigma_{k-2} (\lambda_{; 1, j}) }{\lambda_1 \cdots \lambda_n \lambda_j } \frac{| (X_u)_{j \bar{1},1}    |^2}{1+ \tilde{\lambda}_1}    \\ 
&=    \frac{ \lambda_2 \cdots \lambda_n - \lambda_j \sum_{k=2}^{n-2} c_k \sigma_{k-2} (\lambda_{; 1, j}) +  2 \lambda_j \sum_{k=1}^{n-2} c_k \sigma_{k-1}(\lambda_{; 1, j})  }{(1 + \tilde{\lambda}_1)^2\lambda_1 \cdots \lambda_n \lambda_j }    \bigl |    (  X_u     )_{j \bar{1},1}    \bigr |^2   \\
&\kern2em +2 \Bigl( \frac{ \sum_{k=1}^{n-2} c_k \sigma_{k-1}(\lambda_{; j})}{\lambda_1 \cdots \lambda_n} - \frac{1}{\lambda_j} \Bigr) \frac{ \bigl |    T_{j\bar{1}} \bigr |^2}{(1 + \tilde{\lambda}_1)^2}    \\
&\geq 2 \Bigl( \frac{ \sum_{k=1}^{n-2} c_k \sigma_{k-1}(\lambda_{; j})}{\lambda_1 \cdots \lambda_n} - \frac{1}{\lambda_j} \Bigr) \frac{ \bigl |    T_{j \bar{1}} \bigr |^2}{(1 + \tilde{\lambda}_1)^2} \geq  \frac{ C}{(1 + \tilde{\lambda}_1)^2 }h_j \geq \frac{C}{\lambda_1^2} h_j.  
% &\geq   \sum_{j \neq 1}  h_j \frac{ \bigl |    (  X_u     )_{1 \bar{1},j}    \bigr |^2}{(1 + \tilde{\lambda}_1)^2}   +   \sum_{j \neq 1}  \frac{ 1 }{(1+ \tilde{\lambda}_1)   \lambda_1 \lambda_j}  \Bigl(  \frac{\lambda_1}{1+ \tilde{\lambda}_1} | (X_u)_{1 \bar{1},j} |^2 - \frac{\lambda_1 - B_{11}}{1+ B_{11}} | S_j   |^2 \Bigr) \\
% &\geq    \sum_{j \neq 1}  \frac{ h_j \lambda_j +  1 }{(1+ \tilde{\lambda}_1)^2    \lambda_j}      | (X_u)_{1 \bar{1},j} |^2    -  \sum_{j \neq 1}  \frac{  | S_j   |^2 }{(1+ \tilde{\lambda}_1)    \lambda_j}         \\
% &=   \sum_{j \neq 1}  \frac{ c_1 \bigl |    (  X_u     )_{1 \bar{1},j}    \bigr |^2}{(1 + \tilde{\lambda}_1)^2 \lambda_1 \lambda_2 \lambda_3}   -   \sum_{j \neq 1}  \frac{  | S_j   |^2 }{(1+ \tilde{\lambda}_1)    \lambda_j}         
% \geq    -   \sum_{j \neq 1}  \frac{   | S_j   |^2 }{(1+ \tilde{\lambda}_1)     \lambda_j}         \geq -C,
\end{align*}
\endgroup
In addition, for $j \neq 1$, we can also simplify the following terms in inequality (\ref{eq:3.32}), we have:
\begingroup
\allowdisplaybreaks
\begin{align*}
\label{eq:3.34}
&\kern-2em -   \frac{  h_j    \bigl |    (  X_u    )_{1 \bar{j}, j}    \bigr  |^2}{(1 + \tilde{\lambda}_1)( \tilde{\lambda}_1 - \tilde{\lambda}_j)}  - \frac{2  \sum_{i} h_i h_{ij}  \Re \bigl(   \frac{\partial c_0}{\partial \bar{z}_1}    (X_u)_{j \bar{j}, 1 }     \bigr) }{(1 + \tilde{\lambda}_1) \lambda_1 \cdots \lambda_n \sum_i h_i^2 }   - \frac{2 \Re \bigl(   \frac{\partial c_0}{\partial \bar{z}_1}   (X_u)_{j \bar{j}, {1}} \bigr)  }{(1 + \tilde{\lambda}_1)\lambda_1 \cdots \lambda_n \lambda_j} \tag{3.34} \\
&= -   \frac{  h_j    \bigl |    (  X_u    )_{1 \bar{j}, j}    \bigr  |^2}{(1 + \tilde{\lambda}_1)( \tilde{\lambda}_1 - \tilde{\lambda}_j)}  - \frac{2 \sum_{i} h_i h_{ij}  \Re \bigl(   \frac{\partial c_0}{\partial \bar{z}_1}     (X_u)_{1 \bar{j}, j }     \bigr) }{(1 + \tilde{\lambda}_1) \lambda_1 \cdots  \lambda_n \sum_i h_i^2 }    - \frac{2 \Re \bigl(   \frac{\partial c_0}{\partial \bar{z}_1}  (X_u)_{1 \bar{j}, {j}} \bigr)  }{(1 + \tilde{\lambda}_1)\lambda_1 \cdots  \lambda_n \lambda_j} \\
&\kern2em  - \frac{2\sum_{l} h_l h_{lj}  \Re \bigl(   \frac{\partial c_0}{\partial \bar{z}_1}  T_{j \bar{j}}     \bigr)}{(1 + \tilde{\lambda}_1) \lambda_1 \cdots  \lambda_n \sum_l h_l^2 }    - \frac{2 \Re \bigl(   \frac{\partial c_0}{\partial \bar{z}_1}  T_{j \bar{j}}  \bigr)  }{(1 + \tilde{\lambda}_1)\lambda_1 \cdots \lambda_n \lambda_j} \\
&= -   \frac{  h_j  }{(1 + \tilde{\lambda}_1)( \tilde{\lambda}_1 - \tilde{\lambda}_j)} \Biggl|        (  X_u    )_{1 \bar{j}, j}    + \frac{\sum_i h_i h_{ij}  \frac{\partial c_0}{\partial z_1}  (\tilde{\lambda}_1  - \tilde{\lambda}_j )  }{\lambda_1 \cdots \lambda_n h_j \sum_i h_i^2}   + \frac{  \frac{\partial c_0}{\partial {z}_1} (\tilde{\lambda}_1  - \tilde{\lambda}_j )  }{\lambda_1 \cdots \lambda_n \lambda_j h_j  }       \Biggr|^2  \\
&\kern2em+ \frac{   (\tilde{\lambda}_1 - \tilde{\lambda}_j )  }{(1 + \tilde{\lambda}_1) \lambda_1^2  \cdots  \lambda_n^2 h_j   }   \Biggl|         \frac{\sum_i h_i h_{ij}   \frac{\partial c_0}{\partial z_1}     }{  \sum_i h_i^2}    + \frac{   \frac{\partial c_0}{\partial {z}_1}   }{  \lambda_j    }       \Biggr|^2    - \frac{2  \sum_{i} h_i h_{ij}  \Re \bigl(   \frac{\partial c_0}{\partial \bar{z}_1}   T_{j \bar{j}}  \bigr) }{(1 + \tilde{\lambda}_1) \lambda_1 \cdots \lambda_n \sum_i h_i^2 }  \\
&\kern2em - \frac{2 \Re \bigl(   \frac{\partial c_0}{\partial \bar{z}_1} T_{j \bar{j}}   \bigr) }{(1 + \tilde{\lambda}_1)\lambda_1 \cdots \lambda_n \lambda_j} \\
&\geq  \frac{   (\tilde{\lambda}_1 - \tilde{\lambda}_j )  }{(1 + \tilde{\lambda}_1) \lambda_1^2 \cdots \lambda_n^2 h_j   }   \Biggl|         \frac{\sum_i h_i h_{ij}  \frac{\partial c_0}{\partial z_1}     }{  \sum_i h_i^2}    + \frac{   \frac{\partial c_0}{\partial {z}_1}   }{  \lambda_j    }       \Biggr|^2     - \frac{2  \sum_{i} h_i h_{ij}  \Re \bigl(   \frac{\partial c_0}{\partial \bar{z}_1}  T_{j \bar{j}}  \bigr) }{(1 + \tilde{\lambda}_1) \lambda_1 \cdots \lambda_n \sum_i h_i^2 }  \\
&\kern2em - \frac{2 \Re \bigl(   \frac{\partial c_0}{\partial \bar{z}_1}  T_{j \bar{j}} \bigr)  }{(1 + \tilde{\lambda}_1)\lambda_1 \cdots \lambda_n \lambda_j} \\
&\geq \frac{C}{\lambda_1^2 \cdots \lambda_n^2 h_j}    \Biggl|         \frac{\sum_i h_i h_{ij}       }{  \sum_i h_i^2}    + \frac{   1   }{  \lambda_j    }       \Biggr|^2 + \frac{C \sum_i h_i h_{ij}}{\lambda_1 \cdots \lambda_n \lambda_1 \sum_i h_i^2}  -      \frac{C}{\lambda_1 \cdots \lambda_n \lambda_1 \lambda_j}.
\end{align*}
\endgroup
By Lemma~\hyperlink{L:3.5}{3.5}, we have the following observation:
\begin{align*}
 \Biggl|    \frac{\sum_i h_i h_{ij}       }{  \sum_i h_i^2} \Biggr| \leq  \frac{ \frac{n+1}{\lambda_j \lambda_n^2}    }{\frac{1}{\lambda_n^2}   \bigl( 1 -      {\sum_{k=1}^{n-2} c_k \binom{n-1}{k-1} x_0^{k-n} }  \bigr)^2   }    \leq   \frac{  {n+1}    }{     \lambda_j  \bigl( 1 -   \sup_{c \in S}  {\sum_{k=1}^{n-2} c_k \binom{n-1}{k-1} x_0^{k-n}(c) }  \bigr)^2   } \leq \frac{C}{\lambda_j}.
\end{align*}
The last inequality is due to the fact that for any $c \in \tilde{\mathscr{C}}_n$, $c$ corresponds to a strictly right-Noetherian polynomial, which implies that $x_0 > x_1$ and hence $\sum_{k=1}^{n-2} c_k \binom{n-1}{k-1} x_0^{k-n}(c) < 1$. Thus, $\sum_{k=1}^{n-2} c_k \binom{n-1}{k-1} x_0^{k-n}(c)$ is a positive continuous function on $\tilde{\mathscr{C}}_n$ with value less than $1$. Since $S$ is a compact subset of $\tilde{\mathscr{C}}_n$, the continuous function $\sum_{k=1}^{n-2} c_k \binom{n-1}{k-1} x_0^{k-n}(c)$ attains its maximum on $S$ with maximum value less than $1$.\smallskip

With the above observation, we estimate all the terms in inequality (\ref{eq:3.34}). For convenience, we divide all the terms by $h_j$. By Lemma~\hyperlink{L:2.7}{2.7}, we obtain 
\begingroup
\allowdisplaybreaks
\begin{align*}
\label{eq:3.35}
\frac{   \Bigl|         \frac{\sum_i h_i h_{ij}       }{  \sum_i h_i^2}    + \frac{   1   }{  \lambda_j    }       \Bigr|^2 }{\lambda_1^2 \cdots \lambda_n^2 h_j^2}     &\leq  \frac{C}{\lambda_1^2 \cdots \lambda_n^2 \lambda_j^2 h_j^2} = \frac{C}{\bigl( \sum_{k=0}^{n-2} c_k \sigma_{k} (\lambda_{; j}) \bigr)^2 } \leq \frac{C}{ \bigl(\sum_{k=0}^{n-2}  c_k \binom{n-1}{k} x_1^k \bigr)^2  } \tag{3.35} \\
&\leq \frac{C}{   \bigl( \inf_{c \in S} \sum_{k=0}^{n-2}  c_k \binom{n-1}{k} x_1^k(c) \bigr)^2    } \leq C; \\
\label{eq:3.36}
\biggl | \frac{\sum_i h_i h_{ij}}{\lambda_1 \cdots \lambda_n \lambda_1 h_j \sum_i h_i^2} \biggr |  &\leq   \frac{-C}{\lambda_1 \cdots \lambda_n \lambda_1 \lambda_j h_j}  =  \frac{C}{\lambda_1 \sum_{k=0}^{n-2} c_k \sigma_k(\lambda_{; j})}    \leq \frac{C}{ \lambda_1  \sum_{k=0}^{n-2}  c_k \binom{n-1}{k} x_1^k     } \tag{3.36} \\
&\leq \frac{C}{   \lambda_1 \inf_{c \in S} \sum_{k=0}^{n-2}  c_k \binom{n-1}{k} x_1^k(c)     } \leq \frac{C}{\lambda_1};   \\
\label{eq:3.37}
\biggl |  \frac{1}{\lambda_1 \cdots \lambda_n \lambda_1 \lambda_j h_j} \biggr |  &\leq \frac{1}{\lambda_1 \bigl(  \inf_{c \in S} \sum_{k=0}^{n-2} c_k \binom{n-1}{k} x_1^k(c)  \bigr)}  \leq \frac{C}{\lambda_1}. \tag{3.37}
\end{align*}
\endgroup
Similar to before, for any $c \in \tilde{\mathscr{C}}_n$, we have $x_0(c) > x_1(c)$ and hence $\sum_{k=0}^{n-2} c_k \binom{n-1}{k} x_1^k(c) > 0$. $\sum_{k=0}^{n-2} c_k \binom{n-1}{k} x_1^k(c)$ is a continuous function on $\tilde{\mathscr{C}}_n$ with value greater than $0$. Since $S$ is a compact subset of $\tilde{\mathscr{C}}_n$, so the continuous function $\sum_{k=0}^{n-2} c_k \binom{n-1}{k} x_1^k(c)$ attains its minimum on $S$ with minimum value greater than $0$. \smallskip

By combining inequalities (\ref{eq:3.34}), (\ref{eq:3.35}), (\ref{eq:3.36}), and (\ref{eq:3.37}), we get
\begin{align*}
\label{eq:3.38}
&\kern-2em - \frac{  h_j    \bigl |    (  X_u    )_{1 \bar{j}, j}    \bigr  |^2}{(1 + \tilde{\lambda}_1)( \tilde{\lambda}_1 - \tilde{\lambda}_j)}  - \frac{2  \sum_{i} h_i h_{ij}  \Re \bigl(   \frac{\partial c_0}{\partial \bar{z}_1}    (X_u)_{j \bar{j}, 1 }     \bigr) }{(1 + \tilde{\lambda}_1) \lambda_1 \cdots \lambda_n \sum_i h_i^2 }   - \frac{2 \Re \bigl(   \frac{\partial c_0}{\partial \bar{z}_1}   (X_u)_{j \bar{j}, {1}} \bigr)  }{(1 + \tilde{\lambda}_1)\lambda_1 \cdots \lambda_n \lambda_j} \tag{3.38} \\
&\geq \frac{C}{\lambda_1^2 \cdots \lambda_n^2 h_j}    \Biggl|         \frac{\sum_i h_i h_{ij}       }{  \sum_i h_i^2}    + \frac{   1   }{  \lambda_j    }       \Biggr|^2 + \frac{C \sum_i h_i h_{ij}}{\lambda_1 \cdots \lambda_n \lambda_1 \sum_i h_i^2}  -      \frac{C}{\lambda_1 \cdots \lambda_n \lambda_1 \lambda_j} \\
&\geq C h_j + \frac{C}{\lambda_1} h_j + \frac{C}{\lambda_1} h_j \geq C h_j.
\end{align*}
Similarly, by Theorem~\hyperlink{T:2.3}{2.3} and Lemma~\hyperlink{L:3.5}{3.5}, we have the following inequalities:
\begin{align*}
\label{eq:3.39}
 \biggl | \frac{  \frac{\partial^2 c_0}{\partial z_1 \partial \bar{z}_1 } }{ (1 + \tilde{\lambda}_1) \lambda_1 \cdots \lambda_n } \biggr | &\leq \frac{C}{\lambda_1 \cdots \lambda_n \lambda_1} \leq \frac{C}{\lambda_1 x_0^n }  \leq \frac{C}{\lambda_1 \inf_{c \in S} x_0^n(c) }  \leq \frac{C}{\lambda_1}; \tag{3.39} \\ 
 \label{eq:3.40}
 \frac{  \sum_{i, j} h_i h_j h_{ij}  |  \frac{\partial c_0}{\partial z_1}|^2 }{  (1 + \tilde{\lambda}_1) \lambda_1^2 \cdots \lambda_n^2 ( \sum_i h_i^2)^2}  &\leq  \frac{C\frac{n(n+1)}{\lambda_n^4}}{\lambda_1 \cdots \lambda_n \lambda_1 \frac{1}{\lambda_n^4}  \bigl( 1 -      {\sum_{k=1}^{n-2} c_k \binom{n-1}{k-1} x_0^{k-n} }  \bigr)^2 } \tag{3.40} \\
 &\leq  \frac{C }{ \lambda_1 \inf_{c \in S} x_0^n(c)   \bigl( 1 -   \sup_{c \in S}   {\sum_{k=1}^{n-2} c_k \binom{n-1}{k-1} x_0^{k-n}(c) }  \bigr)^2 } \leq \frac{C}{\lambda_1}.
\end{align*}Notice that here we use the $\Upsilon$-dominance theorem, which implies that $\lambda_1 \cdots \lambda_n \geq x_0^n$. Moreover, $x_0^n(c)$ is a positive continuous function on $\tilde{\mathscr{C}}_n$. Since $S$ is a compact subset of $\tilde{\mathscr{C}}_n$, the continuous function $x_0^n(c)$ attains its minimum on $S$ with minimum value greater than $0$.\smallskip

Similarly, by Lemma~\hyperlink{L:2.7}{2.7}, Lemma~\hyperlink{L:3.1}{3.1}, Lemma~\hyperlink{L:3.5}{3.5}, and above arguments, we have the following
\begingroup
\allowdisplaybreaks
\begin{align*}
\label{eq:3.41}
&\kern-2em h_1 \frac{\bigl |    (  X_u     )_{1 \bar{1},1}    \bigr |^2}{(1 + \tilde{\lambda}_1)^2}     - \frac{2 \sum_{i} h_i h_{1i}  \Re \bigl(  \frac{\partial c_0}{\partial \bar{z}_1}     (X_u)_{1 \bar{1}, 1 }     \bigr)  }{(1 + \tilde{\lambda}_1) \lambda_1  \cdots \lambda_n \sum_l h_l^2 }   -  \frac{2 \Re \bigl(  \frac{\partial c_0}{\partial \bar{z}_1}  (X_u)_{1 \bar{1}, {1}} \bigr)  }{(1 + \tilde{\lambda}_1)\lambda_1 \cdots \lambda_n \lambda_1 } \tag{3.41} \\
&= h_1 \biggl |  \frac{(X_u)_{1 \bar{1}, 1}}{1 + \tilde{\lambda}_1} - \Bigl( \frac{\sum_i h_i h_{1i}}{ \lambda_1 \cdots \lambda_n h_1  \sum_i h_i^2} + \frac{ 1}{ \lambda_1 \cdots \lambda_n \lambda_1 h_1   } \Bigr)   \frac{\partial c_0}{\partial z_1}   \biggr |^2 \\
&\kern2em- h_1 \biggl|   \frac{\partial c_0}{\partial z_1} \biggr|^2 \Bigl( \frac{\sum_i h_i h_{1i}}{ \lambda_1 \cdots \lambda_n h_1  \sum_i h_i^2} + \frac{ 1}{ \lambda_1 \cdots \lambda_n \lambda_1 h_1   } \Bigr)^2 \\
&\geq h_1 \biggl |  \frac{(X_u)_{1 \bar{1}, 1}}{1 + \tilde{\lambda}_1} - \Bigl( \frac{\sum_i h_i h_{1i}}{ \lambda_1 \cdots \lambda_n h_1  \sum_i h_i^2} + \frac{ 1}{ \lambda_1 \cdots \lambda_n \lambda_1 h_1   } \Bigr)   \frac{\partial c_0}{\partial z_1}   \biggr |^2 \\
&\geq 2h_1  \biggl |  \frac{(X_u)_{1 \bar{1}, 1}}{1 + \tilde{\lambda}_1} \biggr|^2 + 2h_1 \biggl| \Bigl( \frac{\sum_i h_i h_{1i}}{ \lambda_1 \cdots \lambda_n h_1  \sum_i h_i^2} + \frac{ 1}{ \lambda_1 \cdots \lambda_n \lambda_1 h_1   } \Bigr)   \frac{\partial c_0}{\partial z_1}   \biggr |^2 \\
&\geq 2h_1  \biggl |  \frac{(X_u)_{1 \bar{1}, 1}}{1 + \tilde{\lambda}_1} \biggr|^2 - C \biggl| \frac{\bigl(\sum_i h_i h_{1i}\bigr)^2}{ \lambda_1^2 \cdots \lambda_n^2 h_1  \bigl( \sum_i h_i^2 \bigr)^2} \biggr| - C \biggl|  \frac{ 1}{ \lambda_1^2 \cdots \lambda_n^2 \lambda_1^2 h_1   }   \biggr | \\
&\geq 2h_1  \biggl |  \frac{(X_u)_{1 \bar{1}, 1}}{1 + \tilde{\lambda}_1} \biggr|^2 - \frac{C}{\lambda_1} = 2A^2h_1 |u_1|^2 - \frac{C}{\lambda_1}.
\end{align*}
\endgroup
Last, by Lemma~\hyperlink{L:3.5}{3.5}, by inequalities (\ref{eq:3.32}), (\ref{eq:3.33}), (\ref{eq:3.38}), (\ref{eq:3.39}), (\ref{eq:3.40}), and (\ref{eq:3.41}), and by choosing $A> 0$ sufficiently large, at the maximum point $q$, when $\lambda_1$ is large, we obtain 
\begin{align*}
0 \geq \mathcal{L}_c    ( \tilde{U}     ) &\geq 
 \bigl( C-A \kappa \bigr)  \sum_i h_i  +   2A^2 h_1 |u_1|^2 - \frac{C}{\lambda_1}   +  C \sum_{j \neq 1} h_j      \\  
%&\geq  \bigl( C-A \kappa \bigr)  \sum_i h_i      +  C \sum_{j \neq 1}   h_j -C \\
% &\kern2em  + h_1 \Bigl |  \frac{(X_u)_{1 \bar{1}, 1}}{1 + \tilde{\lambda}_1} -   \Bigl( \frac{ \sum_l h_l h_{1l}}{ \lambda_1 \lambda_2 \lambda_3 \lambda_4 h_1  \sum_l h_l^2} + \frac{ 1   }{ \lambda_1^2 \lambda_2 \lambda_3 \lambda_4 h_1   } \Bigr)   \frac{\partial c_0}{\partial z_1}  \Bigr |^2 \\
&\geq - (A \kappa - C) \sum_i h_i -  \frac{2 A^2  |u_1|^2 }{\lambda_1}  - \frac{C}{\lambda_1} \\
&\geq \frac{A \kappa  }{2\lambda_n} \Bigl(1- \sup_{c \in S} \sum_{k=1}^{n-2} c_k \mybinom[0.8]{n-1}{k-1} x_0^{k-n}(c) \Bigr) -  \frac{2 A^2  |u_1|^2 }{\lambda_1} \\
&\geq \frac{A \kappa  }{2 \inf_{c \in S} x_0(c)} \Bigl(1- \sup_{c \in S} \sum_{k=1}^{n-2} c_k \mybinom[0.8]{n-1}{k-1} x_0^{k-n}(c) \Bigr) -  \frac{2 A^2  |u_1|^2 }{\lambda_1}.
\end{align*}So, at the maximum point $q$, we get the following estimate for the largest eigenvalue $\lambda_1$:
\begin{align*}
\frac{4A |u_1|^2 \inf_{c \in S} x_0(c)        }{ \kappa \bigl(1 - \sup_{c \in S} \sum_{k=1}^{n-2} c_k \binom{n-1}{k-1} x_0^{k-n}(c) \bigr) } \geq \lambda_1.
\end{align*}
By plugging back into the original test function $U = -Au +G(\Lambda)$, we will obtain a $C^2$ estimate for any $c \in \tilde{\mathcal{P}}^d$ with $c(z) \in S$ and $c_{1}, \cdots, c_{n-2}$ constants. This finishes the proof.
\end{proof}

\subsection{The \texorpdfstring{$C^1$}{} Estimate}
\label{sec:3.2}
Here, we use a blow-up argument proved by Collins--Jacob--Yau \cite{collins20201} to obtain a $C^1$ estimate. One can also check a more general setting considered by Székelyhidi \cite{szekelyhidi2018fully}, or the complex Hessian equation studied by Dinew--Kołodziej \cite{dinew2017liouville}.

\hypertarget{P:4.2}{\begin{fprop}[Collins--Jacob--Yau \cite{collins20201}]
Suppose $u \colon M   \rightarrow \mathbb{R}$ satisfies 
\begin{itemize}
\item[(1)] $X + \sqrt{-1} \partial \bar{\partial} u \geq -K \omega$,
\item[(2)] $\|u\|_{L^\infty(M)} \leq K$,
\item[(3)] $\| \partial \bar{\partial} u \|_{L^\infty(M)} \leq K \bigl( 1 + \sup_M |\nabla u |^2  \bigr)$,
\end{itemize}for a uniform constant $K < \infty$. Then there exists a constant $C$, depending only on $M, \omega$, $X$, and $K$ such that $\sup_M |\nabla u| \leq C$.
\end{fprop}}

\subsection{Higher Order Estimates}
\label{sec:3.3}
The proof follows from Siu \cite{siu2012lectures} and Collins--Jacob--Yau \cite{collins20201}. The equation is elliptic and the solution set is convex, we can exploit the convexity of the solution set to obtain $C^{2, \alpha}$ estimates by a blow-up argument. Furthermore, for higher regularity, we apply the standard Schauder estimates and bootstrapping. \smallskip

The argument here is standard and follows verbatim, so we will not show the details here, the interested reader is referred to the author's previous works \cite{lin2023d, lin2023thesis}.

\section{Continuity Path and the Solvability}
\label{sec:4}
In this section, we consider the solvability of the following equation 
\begin{align*}
\label{eq:4.1}
X^n = d_{n-2} \mybinom[0.8]{n}{n-2} X^{n-2}  \wedge \omega^{2} + \cdots   + d_1 \mybinom[0.8]{n}{1} X  \wedge \omega^{n-1} + d_{0}(z) \mybinom[0.8]{n}{0}   \omega^{n}, \tag{4.1}
\end{align*}which corresponds to the following map $d \colon M \rightarrow \tilde{\mathscr{C}}_n$ with $d_1, \cdots, d_{n-2}$ constants. We always assume a $C$-subsolution exists and by changing representative foe convenience, we say $X$ is this $C$-subsolution. Once we know equation (\ref{eq:4.1}) is solvable, since $M$ is a compact manifold, a representative $X_u \coloneqq X + \sqrt{-1} \partial \bar{\partial} u$ solving equation (\ref{eq:4.1}) is unique by the maximum principle. To obtain the solvability, we apply the method of continuity, the idea is to connect equation (\ref{eq:4.1}) with a well-known solvable equation and see whether we have a priori estimates on this continuity path. With the framework of Section~\ref{sec:3}, we consider the following continuity path first:
\begin{align*}
\label{eq:4.2.t}
X^n = d_{n-2} \mybinom[0.8]{n}{n-2} X^{n-2}  \wedge \omega^{2} + \cdots   + d_1 \mybinom[0.8]{n}{1} X  \wedge \omega^{n-1} + d_{0}(z, t) \mybinom[0.8]{n}{0}   \omega^{n}, \tag{4.2.$t$}
\end{align*}where
\begin{align*}
d_0(z, t) = t d_0(z) + (1-t) \frac{\int_M d_0(z) \omega^n}{\int_M \omega^n} =  t d_0(z) + (1-t) \tilde{d}_0.
\end{align*}Here, we define $\tilde{d}_0 \coloneqq \int_M d_0(z) \omega^n / \int_M \omega^n$. For any $t \in [0, 1]$, we can check that continuity path (\ref{eq:4.2.t}) satisfies the following constraints.
\begin{enumerate}[leftmargin=5.5cm]
	 \setlength\itemsep{0.6em}
\item[Topological constraint:] $\Omega_0 - \sum_{k=1}^{n-2} d_k  \binom{n}{k} \Omega_{n-k} - \int_M d_0(z, t) \omega^n = 0$.
\item[Boundary constraints:] $d_0(z, 1) = d_0(z)$ and $d_0(z, 0) = \tilde{d}_0$.
\item[Positivstellensatz constraint:] $(d_{n-2}, \cdots, d_1, d_0(z, t)) \in  \tilde{\mathscr{C}}_n$.
\item[$\Upsilon$-dominance constraint:] $(d_{n-2}, \cdots, d_1, d_0(z, t)) \in \tilde{\mathcal{P}}^d$.
\end{enumerate}
Here, we denote $\Omega_i \coloneqq \int_M X^{n-i} \wedge \omega^i$. The idea of these constraints follows from the author's work \cite{lin2023d}. First, the topological constraint is natural, along the continuity path (\ref{eq:4.2.t}), the integration over $M$ should be equal. That is, we can check that for any $t \in [0, 1]$, we have
\begin{align*}
\int_M X^n - \sum_{k=1}^{n-2} \int_M d_k \mybinom[0.8]{n}{k}  X^k \wedge \omega^{n-k} - \int_M d_{0} (z, t) \omega^n = \Omega_0 - \sum_{k=1}^{n-2} d_k  \mybinom[0.8]{n}{k} \Omega_{n-k} - \int_M d_0(z) \omega^n = 0.
\end{align*}
Second, the boundary constraint is to make sure that one endpoint is the unsolved equation (\ref{eq:4.1}) and another endpoint is a solvable one. Though at this moment, we do not know whether the equation $X^n = \sum_{k=1}^{n-2} d_k \binom{n}{k} X^{k} \wedge \omega^{n-k} - \tilde{d}_0 \omega^n$ is solvable or not. We will find a second continuity path later and show that the equation $X^n = \sum_{k=1}^{n-2} d_k \binom{n}{k} X^{k} \wedge \omega^{n-k} - \tilde{d}_0 \omega^n$ will be solvable.
\smallskip

Third, the Positivstellensatz constraint guarantees that the equation will not degenerate on the continuity path. From a real algebraic geometry viewpoint, if a general inverse $\sigma_k$ multilinear polynomial is $\Upsilon$-stable but not strictly $\Upsilon$-stable, then the level set degenerates. It is interesting to see whether it is possible to define a solution in this case and the regularity of this solution if such a solution can be defined.
If a general inverse $\sigma_k$ multilinear polynomial is not even $\Upsilon$-stable, then it is not likely that any arbitrary level set is convex and to have these fruitful algebraic properties coming from the $\Upsilon$-stability. If we consider the set $\{ d_0 \colon (d_{n-2}, \cdots, d_1, d_0) \in \tilde{\mathscr{C}}_n \} \subset{\mathbb{R}}$, then this set is a ray and we have $\tilde{d}_0 \in [ \min_{z \in M} d_0(z), \max_{z \in M} d_0(z) ] \subset \{ d_0 \colon (d_{n-2}, \cdots, d_1, d_0) \in \tilde{\mathscr{C}}_n \}$. This implies that $d_0(z, t) = t d_0(z) + (1-t) \tilde{d}_0 \in \{ d_0 \colon (d_{n-2}, \cdots, d_1, d_0) \in \tilde{\mathscr{C}}_n \}$ for any $t \in [0, 1]$.\smallskip

Last, the $\Upsilon$-dominance constraint ensures that we can apply a priori estimates on this continuity path (\ref{eq:4.2.t}). If $(d_{n-2}, \cdots, d_1, d_0(z, t)) \in \tilde{\mathcal{P}}^d$, then by Remark~\hyperlink{R:2.4}{2.4}, $X$ is also a $C$-subsolution to $(d_{n-2}, \cdots, d_1, d_0(z, t))$. Hence, by Theorem~\hyperlink{T:3.1}{3.1}, we have a priori estimates for the solution to $(d_{n-2}, \cdots, d_1, d_0(z, t))$. If we have a priori estimates on this whole continuity path, then we can apply the Arzelà--Ascoli theorem and standard point-set topology argument and say that (\ref{eq:4.2.t}) is solvable for any $t \in [0, 1]$ provided that (\ref{eq:4.2.t}) is solvable when $t = 0$.\smallskip

To see whether equation (\ref{eq:4.2.t}) is solvable when $t = 0$, that is, whether $(d_{n-2}, \cdots, d_1, \tilde{d}_0) \colon M^n \rightarrow  \tilde{\mathscr{C}}^n$ is solvable, similarly we consider a continuity path 
\begin{align*}
\label{eq:4.3.t}
X^n = d_{n-2}(t) \mybinom[0.8]{n}{n-2} X^{n-2}  \wedge \omega^{2} + \cdots   + d_1(t) \mybinom[0.8]{n}{1} X  \wedge \omega^{n-1} + \tilde{d}_{0}(t) \mybinom[0.8]{n}{0}   \omega^{n}, \tag{4.3.$t$}
\end{align*}satisfying the following constraints for any $t \in [0, 1]$.
\begin{enumerate}[leftmargin=5.5cm]
	 \setlength\itemsep{0.6em}
\item[Topological constraint:] $\Omega_0 - \sum_{k=1}^{n-2} d_k(t)  \binom{n}{k} \Omega_{n-k} - \tilde{d}_0(t) \Omega_n =0$.
\item[Boundary constraints:] $ (d_{n-2}(1), \cdots, d_1(1), \tilde{d}_0(1)) = (d_{n-2}, \cdots, d_1, \tilde{d}_0)$  \\
and $ (d_{n-2}(0), \cdots, d_1(0), \tilde{d}_0(0)) \in \CY_n$.
\item[Positivstellensatz constraint:] $(d_{n-2}(t), \cdots, d_1(t), \tilde{d}_0(t)) \in  \tilde{\mathscr{C}}_n$.
\item[$\Upsilon$-dominance constraint:] $(d_{n-2}(t), \cdots, d_1(t), \tilde{d}_0(t)) \in \tilde{\mathcal{P}}^d$.
\end{enumerate}
Here, we denote $\Omega_i \coloneqq \int_M  X^{n-i} \wedge \omega^i$. Same as the arguments above, if we can show that such continuity path (\ref{eq:4.3.t}) exists, then (\ref{eq:4.3.t}) is solvable for any $t \in [0, 1]$ due to the hypothesis that $(d_{n-2}(0), \cdots, d_1(0), \tilde{d}_0(0)) \in \CY_n$. When $t = 0$, (\ref{eq:4.3.t}) becomes the well-known complex Monge--Ampère equation and is solvable due to Yau in \cite{yau1978ricci} provided that a $C$-subsolution exists.

\hypertarget{T:4.1}{\begin{fthm}
The following continuity path $(d_{n-2}(t), \cdots, d_1(t), \tilde{d}_0(t))$ satisfies all the four constraints on (\ref{eq:4.3.t}), where $d_k(t) \coloneqq t^{n-k} d_k$ for $k \in \{1, \cdots, n-2\}$ and 
\begin{align*}
\tilde{d}_0(t) \coloneqq \frac{\Omega_0 - \sum_{k=1}^{n-2} d_k(t)  \binom{n}{k} \Omega_{n-k} }{\Omega_n}.
\end{align*} 
\end{fthm}}

\begin{proof}
Let $d_k(t) \coloneqq t^{n-k} d_k$ for $k \in \{1, \cdots, n-2\}$ and $\tilde{d}_0(t) \coloneqq \bigl(\Omega_0 - \sum_{k=1}^{n-2} d_k(t)  \binom{n}{k} \Omega_{n-k} \bigr)/\Omega_n$. The diagonal restriction of $\lambda_1 \cdots \lambda_n - \sum_{k=1}^{n-2}d_k(t) \sigma_k(\lambda) - \tilde{d}_0(t)$ will be 
\begin{align*}
\label{eq:4.4.t}
x^n - \sum_{k = 1}^{n-2} d_k(t) \mybinom[0.8]{n}{k} x^k - \tilde{d}_0(t) =    x^n - \sum_{k = 1}^{n-2} t^{n-k} d_k  \mybinom[0.8]{n}{k} x^k - \frac{\Omega_0 - \sum_{k=1}^{n-2} t^{n-k} d_k  \binom{n}{k} \Omega_{n-k} }{\Omega_n}. \tag{4.4.$t$}
\end{align*}By fixing $t$, if we denote $x_l(t)$ the largest real root of the $l$-th derivative of (\ref{eq:4.4.t}) with respect to $x$ and $x_l$ the largest real root of the $l$-th derivative of $x^n - \sum_{k = 1}^{n-2} d_k  \binom{n}{k} x^k - \tilde{d}_0$, then we have $x_l(t) = t x_l(1) = tx_l$ for $l \in \{1, \cdots, n-2\}$.\smallskip

First, we show that this choice satisfies the topological constraint. We have
\begin{align*}
\Omega_0 - \sum_{k=1}^{n-2} d_k(t)  \mybinom[0.8]{n}{k} \Omega_{n-k} - \tilde{d}_0(t) \Omega_n = \Omega_0 - \sum_{k=1}^{n-2} d_k(t)  \mybinom[0.8]{n}{k} \Omega_{n-k} - \Bigl( \Omega_0 - \sum_{k=1}^{n-2} d_k(t)  \mybinom[0.8]{n}{k} \Omega_{n-k} \Bigr) = 0.
\end{align*}

Second, we show that this choice satisfies the boundary constraint. When $t = 0$, $d_k(0) = 0$ for $k \in \{1, \cdots, n-2\}$ and $\tilde{d}_0(0) = \Omega_0/\Omega_n > 0$. So $(d_{n-2}(0), \cdots, d_1(0), \tilde{d}_0(0)) \in \CY_n$. When $t = 1$, for $k \in \{1, \cdots, n-2\}$, we have $d_k(1) = 1^{n-k} d_k = d_k$ and 
\begin{align*}
\tilde{d}_0(1) =  \frac{ \Omega_0 - \sum_{k=1}^{n-2} d_k   \binom{n}{k} \Omega_{n-k} }{\Omega_n} = \frac{\tilde{d}_0 \int_M   \omega^n}{\int_M \omega^n} = \tilde{d}_0.
\end{align*}

Third, we show that this choice satisfies the $\Upsilon$-dominance constraint. An elegant way to see this is by Theorem~\hyperlink{T:2.3}{2.3}, because $x_1 \geq   \cdots \geq x_{n-2} \geq 0$, for any $t \in [0, 1]$, we have
\begin{align*}
(x_{n-2}, \cdots, x_1)  \gtrdot (x_{n-2}(t), \cdots, x_1(t)) = (tx_{n-2}, \cdots, tx_1).
\end{align*}Hence, $(d_{n-2}(t), \cdots, d_1(t), \tilde{d}_0(t)) \in \tilde{\mathcal{P}}^d$ for any $t \in [0, 1]$.\smallskip

Last, we show that this choice satisfies the Positivstellensatz constraint. If $x_0(t)$ exists and $x_0(t) > x_1(t) = t x_1$ for all $t \in [0 ,1]$, then we are done. When $t = 1$, we have $x_0(1) = x_0 >x_1(1) = x_1$. When $t = 0$, we have $x_0(0) = \sqrt[n]{\Omega_0/\Omega_n} > x_1(0) = 0$. For any $t \in (0, 1)$, we \hypertarget{claim 1}{\textbf{claim}} that 
\begin{align*}
\Omega_0 - \sum_{k=1}^{n-2} t^{n-k} d_k  \mybinom[0.8]{n}{k} \Omega_{n-k} - t^n \tilde{d}_0 \Omega_n    > 0.
\end{align*}

If the \bhyperlink{claim 1}{\textbf{claim}} is true, then by the fact that $\tilde{d}_0 +  \sum_{k=1}^{n-2} d_k  \binom{n-1}{k} x_1^k > 0$ and by plugging in $x_1(t)$ to equation (\ref{eq:4.4.t}), we get
\begingroup
\allowdisplaybreaks
\begin{align*}
&\kern-2em x_1^n(t) - \sum_{k = 1}^{n-2} d_k(t) \mybinom[0.8]{n}{k} x^k_1(t) - \tilde{d}_0(t) \\
&=  t^n \Bigl( x_1^n - \sum_{k=1}^{n-2} d_k  \mybinom[0.8]{n}{k} x_1^k    \Bigr) -  \frac{\Omega_0 - \sum_{k=1}^{n-2} t^{n-k} d_k  \binom{n}{k} \Omega_{n-k} }{\Omega_n} \\
&=  -t^n  \sum_{k=1}^{n-2} d_k  \mybinom[0.8]{n-1}{k} x_1^k      -  \frac{\Omega_0 - \sum_{k=1}^{n-2} t^{n-k} d_k  \binom{n}{k} \Omega_{n-k} }{\Omega_n}  \\
&= - \frac{1}{\Omega_n} \Bigl( \Omega_0 - \sum_{k=1}^{n-2} t^{n-k} d_k  \mybinom[0.8]{n}{k} \Omega_{n-k} + t^n  \sum_{k=1}^{n-2} d_k  \mybinom[0.8]{n-1}{k} x_1^k \Omega_n  \Bigr) \\
&<  - \frac{1}{\Omega_n} \Bigl( \Omega_0 - \sum_{k=1}^{n-2} t^{n-k} d_k  \mybinom[0.8]{n}{k} \Omega_{n-k} -  t^n  \tilde{d}_0 \Omega_n  \Bigr) < 0.
\end{align*}
\endgroup
This implies that $(d_{n-2}(t), \cdots, d_1(t), \tilde{d}_0(t)) \in  \tilde{\mathscr{C}}_n$ for any $t \in [0, 1]$, which finishes the proof. \smallskip

Now, we prove the \bhyperlink{claim 1}{\textbf{claim}}. By the assumption that $X$ is a $C$-subsolution to $(d_{n-2}, \cdots, d_1, d_0)$, pointwise we always have
\begin{align*}
\label{eq:4.5.l}
X^{n-l} -  \sum_{k = l}^{n-2} d_k \mybinom[0.8]{n-l}{k-l} X^{k-l} \wedge \omega^{n-k} > 0 \tag{4.5.$l$}
\end{align*}as a $(n-l, n-l)$-form for any $l \in \{1, \cdots, n-1\}$. Recall that $X$ is a Kähler form, for any $l \in \{1, \cdots, n-1\}$, by wedging $(n-l, n-l)$-form (\ref{eq:4.5.l}) with $X^l$ and integrating over $M$, we get
\begin{align*}
\label{eq:4.6.l}
\Omega_0 -  \sum_{k=l}^{n-2} d_k \mybinom[0.8]{n-l}{k-l} \Omega_{n-k}  
%&= \int_M X^n - \sum_{k=l}^{n-2} d_k \mybinom[0.8]{n-l}{k-l}   \int_M X^{k} \wedge \omega^{n-k}  \\
&=  \int_M \Bigl(X^{n-l} -  \sum_{k = l}^{n-2} d_k \mybinom[0.8]{n-l}{k-l} X^{k-l} \wedge \omega^{n-k} \Bigr) \wedge X^l  > 0. \tag{4.6.$l$}
\end{align*}
In addition, by the topological constraint, we have
\begin{align*}
\label{eq:4.7}
\Omega_0 -  \sum_{k=1}^{n-2} d_k \mybinom[0.8]{n}{k} \Omega_{n-k} - \tilde{d}_0 \Omega_n  = 0. \tag{4.7}
\end{align*} 

By fixing $t \in (0, 1)$, we consider a sequence $\{a_l\}_{l \in \{0, \cdots, n-1\}}$ satisfying the following recurrence relation with initial $a_0 = t^n$ and
\begin{align*}
a_l = \mybinom[0.8]{n}{l}t^{n-l} - a_0\mybinom[0.8]{n}{l} -   \cdots - a_{l-1}  \mybinom[0.8]{n-l+1}{1} = \mybinom[0.8]{n}{l}t^{n-l} - \sum_{k=0}^{l-1} a_k\mybinom[0.8]{n-k}{l-k}  
\end{align*}for $l  \in \{1, \cdots, n-1\}$. If $a_l > 0$ for any $l \in \{1, \cdots, n-1\}$ and $1 \geq \sum_{l=0}^{n-1} a_l$, then by considering the sum of inequality (\ref{eq:4.6.l}) multiplied by $a_l$ from $l = 1$ to $n-1$, we get
\begingroup
\allowdisplaybreaks
\begin{align*}
0 &< \sum_{l = 1}^{n-1} a_l \Bigl( \Omega_0 -  \sum_{k=l}^{n-2} d_k \mybinom[0.8]{n-l}{k-l} \Omega_{n-k}  \Bigr) = \sum_{l = 1}^{n-1} a_l \Omega_0 - \sum_{l = 1}^{n-1}  \sum_{k=l}^{n-2} a_l d_k \mybinom[0.8]{n-l}{k-l} \Omega_{n-k} \\
&=  \sum_{l = 1}^{n-1} a_l \Omega_0 - \sum_{k = 1}^{n-2}  \sum_{l=1}^{k} a_l  \mybinom[0.8]{n-l}{k-l} d_k\Omega_{n-k} =  \sum_{l = 1}^{n-1} a_l \Omega_0 - \sum_{k = 1}^{n-2}     (t^{n-k} - a_0) d_k\mybinom[0.8]{n}{k} \Omega_{n-k} \\
&=  \sum_{l = 1}^{n-1} a_l \Omega_0  - \sum_{k = 1}^{n-2}   t^{n-k}   d_k\mybinom[0.8]{n}{k} \Omega_{n-k} + t^n \sum_{k = 1}^{n-2} d_k\mybinom[0.8]{n}{k} \Omega_{n-k} \\
&=  \sum_{l = 0}^{n-1} a_l \Omega_0  - \sum_{k = 1}^{n-2}   t^{n-k}   d_k\mybinom[0.8]{n}{k} \Omega_{n-k} - t^n  \tilde{d}_0 \Omega_n \leq \Omega_0 - \sum_{k=1}^{n-2} t^{n-k} d_k  \mybinom[0.8]{n}{k} \Omega_{n-k} - t^n \tilde{d}_0 \Omega_n,
\end{align*}
\endgroup
where the second to last equality is by equation (\ref{eq:4.7}). This proves the \bhyperlink{claim 1}{\textbf{claim}}. So the only remaining things to show are $a_l > 0$ for any $l \in \{1, \cdots, n-1\}$ and $1 \geq \sum_{l=0}^{n-1} a_l$. We use mathematical induction to prove that $a_l = \binom{n}{l}t^{n-l}(1-t)^l$ for any $l \in \{0, \cdots, n-1\}$. When $l = 0$, we have $a_0 = t^n =  \binom{n}{0}t^{n-0}(1-t)^0$. Suppose the statement is true when $l \leq  m$. When $ l = m+1$, we have
\begin{align*}
a_{m+1}  &= \mybinom[0.8]{n}{m+1}t^{n-m-1} - \sum_{k=0}^{m} a_k\mybinom[0.8]{n-k}{m+1-k}  =    \mybinom[0.8]{n}{m+1}t^{n-m-1} - \sum_{k=0}^{m}  \mybinom[0.8]{n-k}{m+1-k}    \mybinom[0.8]{n}{k}t^{n-k}(1-t)^k       \\
&=  \mybinom[0.8]{n}{m+1}t^{n-m-1} \Bigl ( 1  - \sum_{k=0}^{m}     \mybinom[0.8]{m+1}{k}t^{m+1-k}(1-t)^k  \Bigr) =  \mybinom[0.8]{n}{m+1}t^{n-m-1} (1-t)^{m+1},
\end{align*}which is also true. In addition, we can verify that
\begin{align*}
\sum_{l=0}^{n-1} a_l = \sum_{l=0}^{n-1} \mybinom[0.8]{n}{l}t^{n-l}(1-t)^l = 1 - (1-t)^n \leq 1
\end{align*}when $t \in (0, 1)$. This finishes the proof. 
\end{proof}

In conclusion, we have the following.
\begin{fthm}[Solvability]
Suppose that $[X]$ and $[\omega]$ satisfies the following integrability condition: 
\begin{align*}
\int_M X^n = d_{n-2} \int_M X^{n-2} \wedge \omega^2 + \cdots + d_1 \int_M X \wedge \omega^{n-1} + \int_M d_0(z) \omega^n.
\end{align*}with $d = (d_{n-2}, \cdots, d_1, d_0(z)) \colon M \rightarrow \tilde{\mathscr{C}}_n$ for any $z \in M$.
If there exists a $C$-subsolution to $d$, then there exists a unique representative $X_u \coloneqq X + \sqrt{-1} \partial \bar{\partial}u \in [X]$ such that
\begin{align*}
X_u^n = d_{n-2}  X_u^{n-2} \wedge \omega^2 + \cdots + d_1  X_u \wedge \omega^{n-1} +  d_0(z) \omega^n.
\end{align*}
\end{fthm}

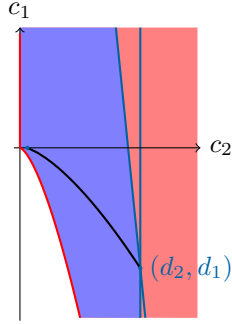
\begin{figure}
\centering
\begin{tikzpicture}[scale=0.4]
\begin{scope}
    %\clip(0, 0) -- (0, 3) -- (4.2325, 3) --  (4.2325, 0);
    \fill[color=red!50,opacity=0.3,thick,domain=0:2,samples=101]
    (0, 4) -- plot (\x,{-2*(\x^(3/2))}) --(5.9,{-2*(2^(3/2))}) --(5.9,4);
  \end{scope}
\begin{scope}
    %\clip(0, 0) -- (0, 3) -- (4.2325, 3) --  (4.2325, 0);
    \fill[color=blue!50,opacity=0.3,thick,domain=0:2,samples=101]
    (0, 4) -- plot (\x,{-2*(\x^(3/2))}) --(4,{-2*(2^(3/2))}) -- (4,-4) --(3.188,4);
  \end{scope}
% \draw [color=red,fill=white] (0, 0) circle[radius= 0.4 em]; 
  \draw[color=UCIB,thick,domain={-2*(2^(3/2))}:4]    plot ({4},\x)             ;
  \draw[color=UCIB,thick,domain={3.188}:4.168]    plot ({\x},{3.2836^3 - 3*3.2836*\x})             ;
  \draw[color=black,thick,domain=0:1]    plot ({(8-7.875*(1-\x))^(2/3)},{-4*(\x)})             ;
 \filldraw [color=UCIB] (4, -4) circle[radius= 0.15 em] node[right] {$(d_2, d_1)$}; 
  \draw[->] (-0.2,0) -- (6,0) node[right] {$c_2$};
  \draw[->] (0,-5.75) -- (0,4) node[above] {$c_1$};
 \filldraw [color=UCIB] (0.25, 0) circle[radius= 0.15 em] ; 
  \draw[color=red,thick,domain=0:3.8]    plot ({0},\x)             ;
  \draw[color=red,thick,domain=0:2]    plot (\x,{-2*(\x^(3/2))})             ;
% \draw[color=UCIB,dashed,domain=0:1]    plot ( {1 + 2*\x} ,{-5.1*\x})            ;
% \draw[color=UCIB,dashed,domain=0:1]    plot ( {3} ,{-5.1*\x})            ;
\end{tikzpicture}
\caption{Continuity path (\ref{eq:4.8})}
\label{fig:4.1}
\end{figure}

\begin{frmk}
The continuity path in Theorem~\hyperlink{T:4.1}{4.1} is inspired by the author's works \cite{lin2023c, lin2023thesis}. When the complex dimension equals four, let $d = (d_2, d_1, d_0) \in \tilde{\mathscr{C}}_4$ with $d_1 < 0$. The author showed that the following continuity path will satisfy all the four constraints: 
\begin{align*}
\label{eq:4.8}
d_{2, \ell} (t) &\coloneqq \bigl( d_2^{3/2} +   {(1-t)\ell d_1}/{2} \bigr)^{2/3}; \quad
d_1(t) \coloneqq td_1;  \quad
d_{0, \ell}(t) \coloneqq \frac{    \Omega_0 -  6d_{2, \ell}(t)      \Omega_2  - 4d_1(t)       \Omega_3}{     \Omega_4}, \tag{4.8}
\end{align*}provided that $\ell = \ell(\Omega_0, \cdots, \Omega_4)$ is sufficiently close to $-2d_2^{3/2}/d_1$. In Figure~\ref{fig:4.1}, we graph the $c_1, c_2$-plane, the pink shaded region is the defining region of $(c_1, c_2)$, the purple region is the cross section of polyhedron $\tilde{\mathcal{P}}^d$ when $c_0 = 0$, and the black curve is the continuity path (\ref{eq:4.8}) projected onto the $c_1, c_2$-plane.
When $\ell = -2d_2^{3/2}/d_1$, then $d_{2, \ell}(t) = t^{2/3}d_2$, $x_{1, \ell}(t) = t^{1/3}x_1$, and $x_{2, \ell}(t) = t^{1/3} x_2$. Here, $x_k$ is the largest real root of the $k$-th derivative of $x^4 - 6d_2 x^2 - 4d_1 x - d_0$ and $x_{k, \ell}$ is the largest real root of the $k$-th derivative of $x^4 - 6d_{2, \ell}(t) x^2 - 4d_1(t) x - d_{0, \ell}(t)$. We could not set $\ell = -2d_2^{3/2}/d_1$ in the previous works because of the lack of a priori estimates crossing different stratas in $\tilde{\mathscr{C}}_4$. We only had a priori estimates on the generic strata of $\tilde{\mathscr{C}}_4$ in \cite{lin2023c, lin2023thesis}, so the path could only be in the generic strata of $\tilde{\mathscr{C}}_4$. This caused the estimates in \cite{lin2023c, lin2023thesis} very complicated and needed to do an equicontinuity argument when $\ell$ is close to $-2d_2^{3/2}/d_1$ but not $-2d_2^{3/2}/d_1$.
\end{frmk}

\Address

\end{document}